\documentclass[final,11pt]{siamart171218}
\usepackage[utf8]{inputenc}
\usepackage[utf8]{inputenc}
\usepackage{url}
\usepackage{xr-hyper}
\usepackage[utf8]{inputenc}
\usepackage{url}
\usepackage{xr-hyper}
\usepackage{amsmath,amsfonts, multirow, graphicx}
\usepackage[caption=false]{subfig} 
\usepackage{algorithm,algpseudocode}
 \usepackage{booktabs}
 \usepackage{changes}

\newtheorem{remark}{Remark}
\usepackage{xcolor}
 \definecolor{mygrn}{RGB}{0,140,70}

\usepackage[mathscr]{eucal}
\usepackage{amsbsy,bm}
\usepackage{amssymb}
\usepackage{hyperref,cleveref}
\usepackage{tikz}
\usepackage{calc}
\usetikzlibrary{calc,arrows}
\usetikzlibrary{decorations.pathreplacing}
\usepackage{changes}
\usepackage[T1]{fontenc}

\newcommand{\T}[2][]{\boldsymbol{#1\mathscr{\MakeUppercase{#2}}}}

\newcommand{\TT}[2]{\T{T}_{\T{#1},#2}} 
\newcommand{\slice}[2]{} 
\newcommand{\TTslice}[2]{ 
\ifx1#2  
	\renewcommand{\slice}[2]{\TT{#1}{1}(i_1,:)}
\else \ifx#2N 
	\renewcommand{\slice}[2]{\TT{#1}{N}(:,i_N)}
\else 
	\renewcommand{\slice}[2]{\TT{#1}{#2}(:,i_{#2},:)}
\fi \fi
\slice{#1}{#2} 
}

\newcommand{\R}{\mathbb{R}}

\makeatletter
\newcommand*{\addFileDependency}[1]{
  \typeout{(#1)}
  \@addtofilelist{#1}
  \IfFileExists{#1}{}{\typeout{No file #1.}}
}
\makeatother

\usepackage{ifthen}
\usepackage{float}
\floatstyle{ruled}
\newfloat{algorithm}{htb}{alg}
\floatname{algorithm}{Algorithm}
\newcounter{algo@row}
\newcounter{algo@rowindent}
\newcommand{\algofont}[1]{\textbf{#1}}
\newcommand{\algonumbersize}[1]{\scriptsize{#1}}
\newcommand{\algopreitem}[1][\arabic{algo@row}]{\texttt{\algonumbersize{#1}}}
\newcommand{\algoitemskip}{\hspace{\value{algo@rowindent}cc}}
\newcommand{\algonewnestedopen}[2]{
	\newcommand{#1}[1][]{%
		\ifthenelse{\equal{##1}{}}{\item}{\item[{\algopreitem[##1]}]}
		\algoitemskip\algofont{#2}%
		\addtocounter{algo@rowindent}{1}%
		\ignorespaces
	}
}
\newcommand{\algonewnestedaux}[2]{
	\newcommand{#1}[1][]{
		\addtocounter{algo@rowindent}{-1}
		\ifthenelse{\equal{##1}{}}{\item}{\item[{\algopreitem[##1]}]}
		\algoitemskip\algofont{#2}%
		\addtocounter{algo@rowindent}{+1}%
		\ignorespaces
	}
}
\newcommand{\algonewnestedclose}[2]{
	\newcommand{#1}[1][]{
		\addtocounter{algo@rowindent}{-1}
		\ifthenelse{\equal{##1}{}}{\item}{\item[{\algopreitem[##1]}]}
		\algoitemskip\algofont{#2}%
		\ignorespaces
	}
}
\newcommand{\algonewcommand}[2]{
	\newcommand{#1}[1][default]{
		\ifthenelse{\equal{##1}{default}}{\item}{\item[{\algopreitem[##1]}]}%
		\algoitemskip\algofont{#2}%
		\ignorespaces
	}%
}
\newcommand{\algonewkeyword}[2]{\newcommand{#1}{\algofont{#2}}}
\algonewcommand{\STATE}{\ignorespaces}
\algonewcommand{\INPUT}{Input: }
\algonewcommand{\pINPUT}{\phantom{Input: }}
\algonewcommand{\COMPUTE}{Compute: }
\algonewcommand{\OUTPUT}{Output: }
\algonewcommand{\pOUTPUT}{\phantom{Output: }}

\algonewnestedopen{\IF}{if }
\algonewnestedaux{\ELSEIF}{else if }
\algonewnestedaux{\ELSE}{else }
\algonewnestedclose{\ENDIF}{end if }
\algonewnestedopen{\FOR}{for }
\algonewnestedclose{\ENDFOR}{end for }
\algonewnestedopen{\WHILE}{while }
\algonewnestedclose{\ENDWHILE}{end while }
\algonewcommand{\BREAK}{break}%

\algonewkeyword{\To}{to }%
\algonewkeyword{\Do}{do }%
\algonewkeyword{\Then}{then }%
\algonewkeyword{\End}{end }%
\algonewkeyword{\AND}{and }%
\algonewkeyword{\True}{true }%
\algonewkeyword{\False}{false }%
\algonewkeyword{\irbleigs}{irbleigs }%
\algonewkeyword{\tridiag}{tridiag}%
\algonewkeyword{\reorth}{reorth}%

\newcommand{\bA}{{\bf A}}
\newcommand{\bB}{{\bf B}}

\newcommand{\bH}{{\bf H}}
\newcommand{\bI}{{\bf I}}

\newcommand{\bL}{{\bf L}}

\newcommand{\bP}{{\bf P}}
\newcommand{\bQ}{{\bf Q}}
\newcommand{\bR}{{\bf R}}
\newcommand{\bS}{{\bf S}}
\newcommand{\bT}{{\bf T}}
\newcommand{\bU}{{\bf U}}
\newcommand{\bV}{{\bf V}}
\newcommand{\bW}{{\bf W}}

\newcommand{\bb}{{\bf b}}

\newcommand{\bd}{{\bf d}}
\newcommand{\be}{{\bf e}}

\newcommand{\br}{{\bf r}}

\newcommand{\bu}{{\bf u}}
\newcommand{\bv}{{\bf v}}
\newcommand{\bw}{{\bf w}}
\newcommand{\bx}{{\bf x}}

\newcommand{\bz}{{\bf z}}

\usepackage{gensymb}

\newcommand{\bSigma}{{\boldsymbol{\Sigma}}}

\newcommand{\bepsilon}{{\boldsymbol{\epsilon}}}

\usepackage{amsmath}
\usepackage{amssymb}
\usepackage[utf8]{inputenc}

\usepackage{optidef}
\usepackage{subfig}
\title{Recycling MMGKS for large-scale dynamic and streaming data}
\author{Mirjeta Pasha \thanks{Department of Mathematics, Tufts University, MA, USA. (\email{mirjeta.pasha@tufts.edu}), (\email{misha.kilmer@tufts.edu})} \and Eric de Sturler \thanks{Department of Mathematics, Virginia Tech, Blacksburg, VA, USA. (\email{sturler@vt.edu})}\and Misha E. Kilmer \footnotemark[1]}
\begin{document}
\maketitle
\begin{abstract}
Reconstructing high-quality images with sharp edges requires the use of edge-preserving constraints in the regularized form of the inverse problem. The use of the $\ell_q$-norm on the gradient of the image is a common such constraint. 
For implementation purposes, the $\ell_q$-norm term is typically replaced with a sequence of $\ell_2$-norm weighted gradient terms with the weights determined from the current solution estimate. 
While (hybrid) Krylov subspace methods can be employed on this sequence, it would require generating a new Krylov subspace for every new two-norm regularized problem. 
The majorization-minimization Krylov subspace method (MM-GKS) addresses this disadvantage by combining norm reweighting with generalized Krylov subspaces (GKS).  After projecting the problem using a small dimensional subspace - one that expands each iteration - the regularization parameter is selected. Basis expansion repeats until a sufficiently accurate solution is found. Unfortunately, for large-scale problems that require many expansion steps to converge, storage and the cost of repeated orthogonalizations presents overwhelming memory and computational requirements. 

In this paper we present a new method, recycled MM-GKS (RMM-GKS), that keeps the memory requirements bounded through recycling the solution subspace. Specifically, our method alternates between enlarging and compressing the GKS subspace, recycling directions that are deemed most important via one of our tailored compression routines.  We further generalize the RMM-GKS approach to handle experiments where the data is either not all available simultaneously, or needs to be treated as such because of the extreme memory requirements.  
Numerical examples from dynamic photoacoustic tomography and streaming X-ray computerized tomography (CT) imaging are used to illustrate the effectiveness of the described methods.
\end{abstract}

\section{Introduction}
A linear inverse problem is one in which the forward model that maps the input to the measured data is assumed to be linear.  If $\b{d} \in \mathbb{R}^{m} $ denotes the measured data, then the linear forward model in the presence of noise is  
\begin{equation}\label{eq:LinSys}
\bA \bx_{\rm true}
+ \be = \bd
\end{equation}
where $\bx_{\rm true} \in \mathbb{R}^n$ is the vectorized representation of the pixel or voxel values corresponding to the medium of interest. The vector $\be \in \R^{m}$ represents additive noise that stems from discretization, measurement and rounding errors.  The inverse problem refers to the recovery of an approximation of $\bx_{\rm true}$ given knowledge of the operator $\bA\in \R^{m\times n}$, measured data $\bd$, and some assumptions on the additive noise $\b{e}$.  

The forward problem (\ref{eq:LinSys}) may be underdetermined, in which case the least squares solution to (\ref{eq:LinSys}) will not be unique.  Moreover, even in the full rank case, the singular values of the forward operator in such problems decay rapidly to zero.  The small singular values cause least squares solutions
to be highly sensitive to the noise present in the measured data.  
A method for mitigating the effects of noise is known as \emph{regularization}, where the original least squares problem is replaced by a well-posed problem whose solution better approximates the desired solution. 
If Gaussian noise is assumed, one standard regularization technique is Tikhonov regularization 
\begin{equation} \label{eq:ellone} 
\min_{\bx} \| \bA \bx - \bd\|_2^2 + \lambda \|\Psi \bx \|_q^q  
\end{equation}
where $0 < q \leq 2$, $\Psi$ is a regularization matrix, and $\lambda$ denotes the regularization parameter whose value determines
the balance of the constraint term (second term) vs. the data misfit term (first term).  Frequently, a choice of $q = 1$ is made to ensure adequate representation of edges and that the constraint term defines a valid norm \cite{buccini2021linearized, chung2019flexible, chung2022}. We will therefore use $q = 1$ throughout the paper, but note that {\it the machinery we develop here can easily be adapted to other choices of $q$}. 

In this paper, we are concerned with the fast, memory efficient, and accurate solution to edge-regularized linear inverse problems of the form (\ref{eq:ellone}) in the context of dynamic or streaming data.  Dynamic linear inverse problems are exceedingly large because of the need to recover multiple images simultaneously and the fact that the regularization couples all these large systems together.  The streaming data case arises when either a) the data may not be available all at once, due to the measurement system setup or b) the data may be available but is so large that memory constraints require we can process only chunks of data at a time.

Due to the non-differentiability when $q=1$, it is common to replace the constraint term in (\ref{eq:ellone}) with a differentiable approximation, and instead minimize
\begin{equation} \label{eq:Je}
\min_{\bx} \mathcal{J}_{\epsilon,\lambda}(\bx) = \min_{\bx} 
\| \bA \bx - \bd\|_2^2 +
 \lambda \sum_{j=1}^n \phi_{\epsilon}( (\Psi \bx_{j} )
 \end{equation}
where $\phi_{\epsilon}(t) = \sqrt{t^2 + \epsilon^2}$ is a smoothed approximation to the $\ell_1$ norm for small $\epsilon$.  

The majorization-minimization approach \cite{hunter2004tutorial, lange2016mm, rodriguez2008efficient} could in theory be employed to solve (\ref{eq:Je}) (see Section \ref{sec:background} for details).  The core of the algorithm requires solution of a sequence of regularized linear least squares problems
\begin{equation} \label{eq:reg2}  \bx^{(k+1)} := \min_{\bx} \|\bA \bx - \bd \|_2^2 + \lambda\| \bP^{(k)}\Psi \bx \|_2^2, \end{equation}
where $\bP^{(k)}$ denotes a diagonal weighting matrix whose entries are determined from the current solution estimate $\bx^{(k)}$.  For a {\it fixed $k$} and known $\lambda$, (\ref{eq:reg2}) could be solved by a Krylov subspace method.  If $\lambda$ is not known,  both $\lambda$ and 
the solution to this can be computed simultaneously using hybrid Krylov subspace methods \cite{gazzola2020inner, chung2019flexible, gazzola2018ir}.

However, each new $k$ requires an entirely new Krylov subspace be built to accommodate the changing regularization term. This is the computational bottleneck, and can also constitute a memory concern.   
To resolve such issues, a majorization-minimization Krylov subspace method (MM-GKS) was proposed \cite{lanza2015generalized, huang2017majorization} that combines this norm reweighting with generalized Krylov subspaces (GKS) \cite{lampe2012large} to solve the reweighted problem.

Nevertheless, computational experiments \cite{huang2017majorization, pasha2021efficient} illustrate that when MM-GKS is used for solving massive problems which require hundreds of iterations to converge, the computational cost can become prohibitive, and the memory requirements can easily exceed the capacity.  
Hence, there is a tremendous need to develop alternatives which will reduce the amount of storage that may be required by MM-GKS without sacrificing the reconstruction quality. Recently, a restarted version of MM-GKS \cite{buccini2023limited} has been proposed that addressees the potentially high memory requirements of the MM-GKS. However, completely ignoring the previous solution subspace increases the risk of throwing away relevant information, thereby necessitating extra work to reintroduce the information into the subspace. 

Thus, in this work we focus on developing a more effective strategy to build a solution subspace through \emph{recycling}.  Our method, called recycling MM-GKS (RMM-GKS), alternately enlarges and compresses the solution subspace while keeping the relevant information throughout the process. 
The idea of recycling Krylov methods for the solution of linear systems and least squares problems is not new \cite{soodhalter2014krylov, parks2006recycling, wang2007large, kilmer2006recycling, ahuja2015recycling, soodhalter2016block, keuchel2016combination}. 
However, subspace recycling in the context of MM-GKS is entirely new.  Moreover, we adapt our RMM-GKS to handle streaming data, such as would arise when either data arrives asynchronously or must be processed as such due to strict memory requirements. 

Our main contributions include the following: 
   \begin{itemize}
    \item a new procedure for the initialization of the solution basis such that edge information is encoded sooner into the process;
    \item {a new recycling variant of MM-GKS in which we}  
   \begin{itemize}
       \item expand by a modest number of vectors to a most $k_{\rm max}$ columns;
       \item compress the basis dimension $k_{\rm min}$ using techniques that we show are designed to retain important information;
   \end{itemize}
\item {a new RMM-GKS approach for the streaming data problem, which can handle situations in which}
    \begin{itemize}
        \item chunks of data (all pertaining to the same solution) are available only sequentially and
        \item storage is at a premium, so we can only store some rows of the system.
     \end{itemize}
  \item rigorous numerical comparisons on several applications illustrating the performance superiority of our method with regard to both storage and quality against competing methods. 
 \end{itemize}   
   \medskip
Our paper is organized as follows.  In Section \ref{sec:background}, we give the background on the MM-GKS background and make the case for more memory and computationally efficient variants. 
We introduce our recycled MM-GKS algorithm in Section \ref{sec: recycle}, and give details on its memory and computational efficiency.
We show, in Section \ref{sec: streaming}, how we can also adapt our RMM-GKS to handle streaming data. 
Extensive numerical experiments are presented in
Section \ref{sec: NumericalExperiments}. We  
provide a final summary in Section ~\ref{sec: conclusion}.  

\section{Background} \label{sec:background}
To make the paper self-contained, we first provide some background on the MM-GKS methods.
First proposed in \cite{lanza2015generalized}, the MM-GKS method computes a stationary point of the functional $\mathcal{J}_{\epsilon}(\bx)$ by using a majorization-minimization approach which constructs a sequence of iterates $\bx^{(k)}$ that converge to a stationary point of $\mathcal{J}_{\epsilon,\lambda}(\bx)$. The functional $\mathcal{J}_{\epsilon,\lambda}(\bx)$ is majorized at each iteration by a quadratic function $\bx \rightarrow \mathcal{Q}(\bx, \bx^{(k)})$. We briefly outline the MM-GKS method to make the paper self contained for the reader. 

\subsection{The majorization step}
Huang et al. \cite{huang2017majorization} describe two approaches to construct a quadratic 
tangent majorant for (\ref{eq:Je}) at an available approximate solution $\bx=\bx^{(k)}$. The 
majorants considered in \cite{huang2017majorization} are referred to as adaptive or fixed 
quadratic majorants. The latter are cheaper to compute, but may give slower convergence. 
In this paper we focus on the adaptive quadratic majorants, but all results also 
hold for fixed quadratic majorants. First we give the following definition.
\begin{definition}[\cite{huang2017majorization}]\label{def: 1}
The functional $\bx\mapsto\mathcal{Q}(\bx,\bv):\R^{n}\rightarrow \R$ is a 
quadratic tangent majorant for $\mathcal{J}_{\epsilon,\lambda}(\bx)$ 
at $\bx=\bv\in\R^n$, if for all $\bv \in \R^n$,
\begin{enumerate}
\item $\bx\mapsto\mathcal{Q}(\bx,\bv)$ is quadratic,
\item $\mathcal{Q}(\bv,\bv)=\mathcal{J}_{\epsilon,\lambda}(\bv)$,
\item $\bigtriangledown_{\bx}\mathcal{Q}(\bv,\bv)=\bigtriangledown_{\bx}\mathcal{J}_{\epsilon,\lambda}(\bv)$,
\item $\mathcal{Q}(\bx,\bv)\geq \mathcal{J}_{\epsilon,\lambda}(\bx)\quad\forall \bx\in \R^{n}$.
\end{enumerate}
\end{definition}
Let $\bx^{(k)}$ be an available approximate solution of \eqref{eq:Je}. We define the vector $\bu^{(k)}=\Psi\bx^{(k)}$ and  the weights
\begin{equation}\label{eq: weightslplq}
\bw_{\epsilon}^{(k)}=\left((\bu^{(k)})^2+\epsilon^2\right)^{-1/2},
\end{equation}
where all operations in the expressions on the right-hand sides, including squaring, are 
element-wise.
We define a weighting matrix 
\begin{equation}\label{eq: weightingMatrix}
    \bP_{\epsilon}^{(k)} = ( \text{diag}(\bw^{(k)}_{\epsilon}))^{1/2}
    \end{equation}
and consider the quadratic tangent majorant for the functional $\mathcal{J}_\epsilon(\bx,\lambda)$ as
\begin{equation}\label{eq: QuadraticMajorantQ}
\begin{array}{rcl}
\mathcal{Q}(\bx, \bx^{(k)})  &=&
\|\bA\bx-\bd\|^{2}_{2}
+\lambda \left(\|\bP_{\epsilon}^{(k)}\Psi\bx\|^{2}_{2}\right)+c,
\end{array}
\end{equation}
where $c$ is a suitable constant\footnote{A non-negative value of $c$ is technically required for the second condition in the definition to hold. However, since the minimizer is independent of the value of $c$, we will not discuss it further here.} that is independent of $\bx^{(k)}$.
We refer the reader to \cite{lanza2015generalized} for the derivation of the weights $\bw_{\epsilon}^{(k)}$.

\subsection{The minimization step}\label{sec: MM_GKS}
We briefly describe here the computation of the approximate solution $\bx^{(k+1)}$ by considering the quadratic majorant formulation \eqref{eq: QuadraticMajorantQ}.
An approximate solution at iteration $k+1$, $\bx^{(k+1)}$ can be determined as the zero of the gradient of \eqref{eq: QuadraticMajorantQ} by solving the normal equations
\begin{equation}\label{eq: normaleqQuadMajorant}
(\bA^T\bA + \lambda  \Psi^{T} (\bP_{\epsilon}^{(k)})^2\Psi)\bx = \bA^T\bd.
\end{equation}
The system \eqref{eq: normaleqQuadMajorant} has a unique solution if
\begin{equation}\label{eq: NANL}
\mathcal{N}(\bA^T\bA)\cap \mathcal{N}(\Psi^{T} (\bP_{\epsilon}^{(k)})^2\Psi)=\{0\},
\end{equation}
is satisfied. 
This condition typically holds in practice, so we will not concern ourselves with it further. The solution $\bx^{(k+1)}$ of \eqref{eq: normaleqQuadMajorant} is the unique minimizer\footnote{In making this statement, we are assuming that $\lambda$ is fixed and known. The case when it is not fixed, which is usually the case in practice, will be discussed in the Appendix.} of the quadratic tangent majorant function $\bx \rightarrow \mathcal{Q}(\bx, \bx^{(k)})$, hence a solution method is well defined.

Unfortunately, solving \eqref{eq: normaleqQuadMajorant} for large matrices $\bA$ and $\Psi$ may be computationally demanding or even prohibitive.  Therefore, the idea presented in earlier work \cite{huang2017majorization} is to compute an approximation by projections onto a smaller dimensional problems.

Their generalized Golub-Kahan method first determines an initial reduction of $\bA$ to a small bidiagonal matrix 
by applying $1\leq\ell\ll\min\{m,n\}$ steps of Golub--Kahan bidiagonalization to $\bA$ with 
initial vector $\bd$. This gives a decomposition 
\begin{equation}\label{bdiag}
\bA\bV_{\ell}=\bU_{\ell}\bB_{\ell},
\end{equation}
where the matrix $\bV_{\ell}\in\R^{n\times\ell}$ has orthonormal columns that span the Krylov
subspace $K_\ell(\bA^T\bA,\bA^T\bd)={\rm span}\{\bA^T\bd,(\bA^T\bA)\bA^T\bd,\ldots,(\bA^T\bA)^{\ell-1}\bA^T\bd\}$, the
matrix $\bU_{\ell}\in\R^{m\times(\ell+1)}$ has orthonormal columns. 
The matrix
$\bB_{\ell}\in\R^{(\ell+1)\times\ell}$ is lower bidiagonal. 

Let $k=0$, let $\bx^{(0)}$ be given and assume $\bP^{(k)}$ is obtained for this initial solution. 
We compute the QR factorizations 
\begin{align}\label{eq: QR}
\bA\bV_{\ell} = \bQ_{\bA}\bR_{\bA}, \quad
\bP^{(k)}\Psi\bV_{\ell} = \bQ_{\Psi}\bR_{\Psi}, 
\end{align}
where
$\bQ_{\bA}\in\R^{m\times\ell}$ and $\bQ_{\Psi}\in\R^{s\times\ell}$ have orthonormal columns and 
$\bR_{\bA}\in\R^{\ell\times\ell}$ and $\bR_{\Psi}\in\R^{\ell\times\ell}$ are upper triangular matrices. 

If we now constrain our next solution, $\bx^{(k+1)}$, to live in the space $\bV_{\ell}$, the minimization problem \eqref{eq: normaleqQuadMajorant} will simplify to finding
\begin{equation}\label{eq: minKryov2}
\bz^{(k+1)}=\arg\min_{\bz\in\R^\ell}\left\|\left[\begin{array}{c} \bR_{\bA}\\ \sqrt{\lambda} \bR_{\Psi}\end{array}\right]
\bz-\left[\begin{array}{c} \bQ_{\bA}^T\bd\\ 
0\end{array}\right]\right\|_2^2,
\end{equation}
followed by the assignment $\bx^{(k+1)} = \bV_{\ell}\bz^{(k+1)}$.

Now we can compute the residual vector corresponding to the current normal equations:
\begin{equation*}\label{eq: residual}
\br^{(k+1)}=\bA^T(\bA\bV_{\ell}\bz^{(k+1)} -\bd)+\lambda  \Psi^T(\bP^{(k)})^2\Psi\bV_{\ell}\bz^{(k+1)}.
\end{equation*}
If this is not suitably small, we expand the solution space by dimension one.
Specifically, we use the normalized residual vector $\bv_{\rm new}=\br^{(k+1)}/\|\br^{(k+1)}\|_2$ to expand the solution subspace. We define the matrix 
\begin{equation}\label{eq: enlargeMMGKS}
\bV_{\ell+1}=[\bV_{\ell},\bv_{\rm new}]\in\R^{n\times(\ell+1)},\end{equation} whose columns form an orthonormal basis 
for the expanded solution subspace. We note here that in exact arithmetic  $\bv_{\rm new}$ is orthogonal to the columns of $\bV_{\ell}$, but in computer arithmetic the vectors may lose orthogonality, hence reorthogonalization of the $\bV_{\ell+1}$ is needed.
Now, $k \leftarrow k+1$,  $\bP^{(k)}$ is recomputed for the current solution estimate, and the process in steps \eqref{eq: QR} - \eqref{eq: enlargeMMGKS} is repeated, expanding the initial $\ell$-dimensional solution space by 1 each iteration until a suitable solution is reached.  

\begin{remark}
Up to now, we have ignored the issue of parameter selection.  The process outlined above should produce a sequence of iterates $\bx^{(k)}$ that converge to the minimizer of $\mathcal{J}_{\epsilon,\lambda}$.  However, in reality, the best choice of $\lambda$ to define $\mathcal{J}_{\epsilon,\lambda}$ is not known a priori.  On the other hand, a suitable regularization parameter $\lambda^{(k)}$ for each individual projected problem (\ref{eq: minKryov2}) can easily be determined by known 
heuristics, such as generalized cross validation (GCV), since the dimension of this problem is small. Thus, $\bx^{(k+1)}$ is taken to be the 
estimated solution to the minimizer of the quadratic tangent majorant for $\mathcal{J}_{\epsilon,\lambda^{(k)}}$ for 
the $\lambda^{(k)}$ that has been selected at the current iteration.   
\end{remark}
We summarize the MM-GKS method in Algorithm \ref{Alg: MM-GKS}.

\begin{algorithm}[!ht]
\caption{MM-GKS \cite{lanza2015generalized}}%
\label{Alg: MM-GKS}
\begin{algorithmic}[1]
	\Require{$\bA, \Psi, \bd, \bx^{(0)}, \epsilon$} 
	\Ensure{An approximate solution $\bx^{(k+1)}$}
	\Function{$\bx^{(k+1)} = $ MM-GKS }{$\bA, \Psi, \bd, \bx^{(0)}, \epsilon$}
	\State	Generate the initial subspace basis $\bV_{\ell}\in \R^{n\times \ell}$ such that $\bV_{\ell}^{T}\bV_{\ell}=\bI$\;
\FOR {$k=0,1,2,\ldots$ \text{until convergence}}{
\State	 $\bu^{(k)}=\Psi \bx^{(k)}$\;
\State $\bw^{(k)}_{\epsilon}=\left((\bu^{(k)})^2+\epsilon^2\right)^{-1/2}$\;	
\State $\bP_{\epsilon}^{(k)} = (\rm diag(\bw_{\epsilon}^{(k)}))^{1/2}$
\State $\bA \bV_{\ell+k}$ and $\bP_{\epsilon}^{(k)}\Psi \bV_{\ell+k}$\;
\State	$\bA\bV_{\ell+k} = \bQ_{\bA}\bR_{\bA}$ and $\bP_{\epsilon}^{(k)}\Psi \bV_{\ell+k}=\bQ_{\Psi}\bR_{\Psi}$\;\Comment{Compute/Update the QR}
\State Select $\lambda^{(k)}$ by heuristic (e.g. GCV) on (\ref{eq: minKryov2}). 
\State $\bz^{(k+1)}$ to satisfy (\ref{eq: minKryov2}) with selected $\lambda^{(k)}$
\State $\bx^{(k+1)}=\bV_{\ell+k}\bz^{(k+1)}$\;
\State $\br^{(k+1)}=\bA^T(\bA\bV_{\ell+k} \bz^{(k+1)} -\bd)+\lambda^{(k)} \Psi^T\bP_{\epsilon}^{(k)}\Psi\bV_{\ell+k}\bz^{(k+1)}$\;
\State		 $\br^{(k+1)} = \br^{(k+1)} - \bV_{\ell+k}\bV_{\ell+k}^{T}\br^{(k+1)}$\; \Comment{Reorthogonalize, if needed}
\State		 $\bv_{\rm new}=\frac{\br^{(k+1)}}{\|\br^{(k+1)}\|_{2}}$: $\bV_{\ell+k+1}=[\bV_{\ell+k}, \bv_{\rm new}]$\;\Comment{Enlarge the solution subspace}}
\ENDFOR
	\EndFunction
\end{algorithmic}
\end{algorithm}

\section{Majorization minimization generalized Krylov subspace with recycling (RMM-GKS)}\label{sec: recycle}
MM-GKS methods have been widely used to efficiently solve large-scale inverse problems \cite{buccini2020modulus, pasha2021efficient}. However, for large-scale problems requiring many basis expansion steps (i.e., large $k$) to converge, storing the necessary solution basis vectors can easily exceed the memory capacities or require an excessive amount of computational time. 
The method we propose, called RMM-GKS, modifies MM-GKS by a) initializing with a better, initial skinny search space and then by
b) keeping the memory requirements constant without sacrificing the reconstruction quality.  Indeed, numerical results show that we can improve the reconstruction quality over typical MM-GKS. The regularization parameter will still be cheaply and automatically determined on small subspaces at each iteration. We alternate between two main steps, enlarging and compressing until a desired reconstructed solution is obtained. 

\subsection{Initialization}
As noted in Section \ref{sec: MM_GKS}, MM-GKS is initialized with $\bV_{\ell}$.  The columns of the matrix $\bV_{\ell}$ form (in exact arithmetic) an orthonormal basis for the Krylov subspace $\mathcal{K}_{\ell}(\bA^T \bA,\bA^T \bd)$. However, this basis is known~\cite{golub2013matrix} to be fairly smooth for small $\ell$, thereby delaying convergence to solutions with edge information.  
Therefore, we insert another stage into the initialization of the solution space.
Our method begins as if we were performing the initialization phase and one step of MM-GKS: 
we determine $\bV_{\ell}$ and an initial solution $\bx^{(0)}$ as before. We then use $\bx^{(0)}$ to define $\bP_{\epsilon}^{(0)}$, and QR factor
both $\bA \bV_{\ell}$, $\bP_{\epsilon}^{(0)} \Psi \bV_{\ell}$.  We find the optimal regularization
parameter $\lambda^{(0)}$, and corresponding
solutions $\bz^{(0)}$, $\bx^{(0)}$, and overwrite $\bP_{\epsilon}^{(0)}$ using our new
$\bx^{(0)}$.  We do one step of MM-GKS expansion to return $\lambda^{(1)}$, a new solution estimate $\bx^{(1)}$, from which we get an updated weighting matrix $\bP^{(1)}_{\epsilon}$.   

Thus far, we have matched the initialization phase of MM-GKS and completed one MM-GKS expansion step.  However, instead of continuing to perform MM-GKS expansions from this point, we use what we now have to generate {\it a new initial seed basis that is encoded with the edge information we have obtained so far}.  Specifically, we generate a matrix $\bV_{k_{\rm min}}$ whose columns are an orthonormal basis for the Krylov subspace
\[ \mathcal{K}_{k_{\rm min}}(\bA^T \bA + \lambda^{(1)} \Psi^T \bP_\epsilon^{(1)} \Psi, \bA^T \bd).\]
Here, $k_{\rm min}$ is the desired minimal basis dimension that we will use throughout.  

\begin{remark}
In comparison with the generation of the  initialization space $\bV_{\ell}$ used in MM-GKS, to get our initialization space $\bV_{k_{\rm min}}$, we incur the additional cost of 1 MM-GKS expansion step, and $k_{\rm min}$ extra matvecs with $\bA$, $\bP_{\epsilon}^{(1)} \Psi$ and their respective transposes. Here, we do not assume that $k_{\rm min} = \ell$.  
\end{remark}

\subsection{Recycling}
\comment{
We assume that an initial approximation of the desired solution $\bx^{(0)}$ and a matrix $\bV_{\ell_0} \in \R^{m \times \ell_0}$ with orthonormal columns that span the search space in
which we wish to find the solution are given. Considering the initial approximation $\bx^{(0)}$, we compute the weighting matrix as in \eqref{eq: weightingMatrix}. We compute the QR factorization of the skinny matrices 
\begin{equation}\label{eq: startingQR}
\bA\bV_{\ell_0} = \bQ_{\bA} \bR_{\bA}\\
\bP^{(0)}_{p, \epsilon}\bA\bV_{\ell_0} = \bQ_{\Psi} \bR_{\Psi}
\end{equation}
Substituting \eqref{eq: startingQR} into \eqref{eq: normaleqQuadMajorant} allows us to compute a better approximation of the desired solution by solving 
$\bz^{(1)}=( \bR_{\bA}^T\bR_{\bA} + \lambda^{(0)} \bR_{\Psi}^T\bR_{L})^{-1}\bR_{\bA}^T\bQ_{\bA}\bd $\ and projecting to the original subspace, $\bx^{(1)}=\bV_{\ell_0}\bz^{(1)}$. Choosing a dimension $\ell$ relatively small and using the regularization parameter $\lambda^{(0)}$ and the weighting matrix $\bP^{(0)}_{p, \bepsilon}$ allows us to compute a new basis for the solution subspace that encodes information from the regularization operator $\bP^{(0)}_{p, \bepsilon}\Psi$. Notice that such subspace is generated only once. }

Our recycling technique is an iterative process that consists of two main phases:  a basis enlargement phase and a basis compression phase.  In the compression phase, we use information we have to determine a subspace of small dimension that it is appropriate to recycle as a new solution space. We now describe each phase in turn.

\subsubsection{Basis Expansion}
This phase essentially mimics the MM-GKS expansion phase.  We assume that our solution basis is not permitted to hold more than $k_{\rm max}$ columns.  Thus, from an initial solution space of $k_{\rm min}$ columns, we can expand at most $s = k_{\rm max} - k_{\rm min}$ columns. For the first expansion step, we begin with the solution space $\bV_{k_{\rm min}}$ that is initialized as described above. To keep consistent with the numbering starting from 0, we overwrite $\bx^{(0)}, \bP^{(0)}$ with the $\bx^{(1)}, \bP^{(1)}$ that were the estimates which we used to generate $\bV_{k_{\rm min}}$. For subsequent expansion steps that follow a compression step, we use the output of Algorithm \ref{Alg: Compress} to give us the initial solution space and guesses. Then we repeat steps 4-14 of Alg 2.1 for
$j = 0,\ldots,s-1$ or until a convergence criteria is met.  If we have not detected convergence by the end of the $s$ steps, 
we will have a solution subspace $\bV_{k_{\rm max}}$
and a solution estimate $\bx^{(k_{\rm max})}$.
A summary of the enlarging phase is presented in Algorithm \ref{Alg: Enlarge}.

\begin{algorithm}[!ht]
\caption{Enlarge}%
\label{Alg: Enlarge}
\begin{algorithmic}[1]
	\Require{$\bA, \Psi, \bV_{k_{\rm min}}, \bd, \bx^{(0)}, \bepsilon$} 
	\Ensure{An approximate solution and regularization parameter, the enlarged subspace and reduced matrices, $[\bx^{(k_{\rm max})}, \lambda^{(k_{\rm max})}, \bV_{k_{\rm max}}, \bR_{\bA}, \bR_{\Psi}]$}.
	\Function{$[\bx^{(k_{\rm max})}, \lambda^{(k_{\rm max})}, \bV_{k_{\rm max}}, \bR_{\bA}, \bR_{\Psi}] = $ Enlarge }{$\bA, \Psi, \bV_{k_{\rm min}}, \bd, \bx^{(0)}, \epsilon$, $s$, $tol_1$}
\FOR {$k=0,1,2,\ldots, s-1$}
	\State	 $\bu^{(k)}=\Psi \bx^{(k)}$\;
	\State	 $\bw^{(k)}_{\epsilon}=\left((\bu^{(k)})^2+\epsilon^2\right)^{-1/2}$\;	
	\State Compute $\bP_{\epsilon}^{(k)} = (\rm diag(\bw_{\epsilon}^{(k)}))^{1/2}$
	\State $\bA \bV_{k_{\rm min}+k}$ and $\bP_{\epsilon}^{(k)}\Psi \bV_{k_{\rm  min}+k}$\;
	\State
          $\bA\bV_{k_{\rm min}+k} = \bQ_{\bA}\bR_{\bA}$ and $\bP_{\epsilon}^{(k)}\Psi \bV_{k_{\rm min}+k}=\bQ_{\Psi}\bR_{\Psi}$\; \Comment{Compute/Update the QR}
          \State $\lambda^{(k)}  = \min_{\lambda}\Theta\left(\lambda\right)$\; \Comment{Define the regularization parameter}
	\State $\bz^{(k+1)}=( \bR_{\bA}^T\bR_{\bA} + \lambda^{(k)} \bR_{\Psi}^T\bR_{\Psi})^{-1}\bR_{\bA}^T\bQ_{\bA}^T\bd $\;\Comment{Solve the small problem}
	\State $\bx^{(k+1)}=\bV_{k_{\rm min}+k}\bz^{(k+1)}$\;
	\State 
	$\br^{(k+1)}=\bA^T(\bA\bV_{k_{\rm min}+k} \bz^{(k+1)} -\bd)+\lambda^{(k)} \Psi^T(\bP_{\epsilon}^{(k)})^2\Psi\bV_{k_{\rm min}+k}\bz^{(k+1)}
	$\;
	\State	 $\br^{(k+1)} = \br^{(k+1)} - \bV_{k_{\rm min}+k}\bV_{k_{\rm min}+k}^{T}\br^{(k+1)}$\;
    \Comment{Reorthogonalize if needed}
	\State	 $\bv_{\rm new}=\frac{\br^{(k+1)}}{\|\br^{(k+1)}\|_{2}}$
 \State $\bV_{k_{\rm min}+k+1}=[\bV_{k_{\rm min}+k}, \bv_{\rm new}]$\;\Comment{Enlarge the solution subspace}
 
\IF{$\|\bx^{(k+1)}-\bx^{(k)}\|_2/\|\bx^{(k)}\|_2\leq tol_1$}
\STATE \text{break;}
\ENDIF
\ENDFOR
\State $\bx^{(k_{\rm max})} = \bx^{(k+1)}$\;
\State $\bV_{k_{\rm max}} = \bV_{k_{\rm min}+k+1}$\;
\State $\bP^{(k_{\rm max})}_{\epsilon} = \left((\Psi \bx^{(k_{\rm max})})^2+\epsilon^2\right)^{-1/4}$
\EndFunction
\end{algorithmic}
\end{algorithm}

\subsubsection{Basis Compression}\label{subsec: basiscompression}
If the solution obtained at the end of the current expansion phase is not sufficiently accurate, we want to reduce the dimension of the solution space to $k_{\rm min}$, but do it in a way that {\it does not lose what we have learned about the solution thus far}. 
While many compression methods are possible (see for instance \cite{jiang2021hybrid}), in this section we describe the basis compression step in such a way that it is independent on the compression technique used. Then, in Section \ref{sec: compression}, we give details on specific compression routines that give the desired properties.  

Assume the basis built thus far is $\bV_{k_{\rm max}}$, that we have the associated triangular factors $\bR_{\bA}, \bR_{\Psi}$ of size $k_{\rm max} \times k_{\rm max}$ 
and that we have a current regularization parameter estimate, 
 $\lambda_{\rm curr}$.  
 We pass these to our compression subroutine, which we denote by $\chi$, to obtain a mixing matrix of size $k_{\rm max} \times k_{\rm min}-1$, $\bW$:

\begin{equation}\label{eq: xi}
    \bW = \chi\left( \bR_{\bA}, \bR_{\Psi},  \bQ_{\bA}, \bd, \lambda_{\rm curr}\right).
\end{equation}

\emph{The most relevant} $k_{\rm min}-1$ columns of $\bV_{k_{\rm max}}$ are defined according to
\begin{equation}\label{eq: chi}
\widetilde{\bV} := \bV_{k_{\rm max}} \bW,
\end{equation}
so that $\widetilde{\bV}$ 
will have $k_{\rm min}-1$ columns.

Now $range(\widetilde{\bV})$ has encoded within it the most relevant information about previous solutions. But the current solution, $\bx_{\rm curr}$, is not necessarily contained in the space. Therefore, we computed the orthogonal projection of the current solution onto $range(\widetilde{\bV})^{\perp}$: 
 \begin{equation}\label{eq: rangePerp}
     \bx_{\rm new} : = \mbox{prj}_{\widetilde{\bV}^{\perp}}( \bx_{\rm curr}), 
 \end{equation} 
and then we add the normalized projection as the last column of our new solution space as
\begin{equation}\label{eq: addSol}
    \bx_{\rm new} \leftarrow \bx_{\rm new}/\| \bx_{\rm new} \|_2, \qquad \bV_{k_{\rm min}} := [\widetilde{\bV},\bx_{\rm new} ].\end{equation}

\comment{
we repeat refining the solution subspace by removing the less important search directions during the compression and by adding new vectors during the enlarge phase. Once the solution subspace reaches the maximum capacity $\bV_{\ell_{\rm max}}$, we seek to select a new subspace $\mathcal{R}\left(\bV^{(\rm new)}_{\ell_{\rm min}-1}\right) \subset \mathcal{R}\left(\bV_{\ell_{\rm max}}\right)$, where $\bV^{(\rm new)}_{\ell_{\rm min}} = \left[\bV^{(\rm new)}_{\ell_{\rm min}-1}, \bx_{\rm new} \right]$. The compressed subspace $\mathcal{R}\left(\bV^{(\rm new)}_{\ell_{\rm min}-1}\right)$ can be constructed by wisely selecting the most important components through a matrix $\bW_{\ell_{\rm min}-1}$ with orthonormal columns. We can construct such matrix from an SVD (or other compression approaches) of the relatively small matrix $\begin{bmatrix}
    \bR_{\bA}\\
    \lambda^{1/2} \bR_{\Psi}
\end{bmatrix}$, i.e.,  $\begin{bmatrix}
    \bR_{\bA}\\
    \lambda^{1/2} \bR_{\Psi}
\end{bmatrix} = \bU\bSigma \bW^T$ and select the first $\ell_{\rm min}-1$ columns of $\bW$. If $\bx^{(i+1)}  \not\in \bV_{\ell_{\rm min}}$, we consider the solution obtained before the compression and add it to the solution subspace after normalizing it, i.e., 
$\bar{\bx} = \frac{\bar{\bx}}{\|\bar{\bx}\|_2}$, with  \\$\bar{\bx} = \left(\bx^{(i+1)} - \bV_{\ell_{\rm min }-1}\bV_{\ell_{\rm min }-1}^{T}\bx^{(i+1)}\right)/\|\left(\bx - \bV_{\ell_{\rm min }-1}\bV_{\ell_{\rm min }-1}^{T}\bx^{(i+1)})\|_2\right)$. Otherwise, $\bV_{\rm \ell_{\rm min}} = \bV_{\ell_{\rm min}-1}$. 
}

The compression procedure is summarized in Algorithm \ref{Alg: Compress},
and the complete RMM-GKS is in Algorithm \ref{Alg: RMMGKS}.

\begin{algorithm}[!ht]
\caption{Compress}%
\label{Alg: Compress}
\begin{algorithmic}[1]
	\Require{$\bV_{k_{\rm max}}, \bR_{\bA}, \bR_{\Psi}, \bx, \lambda$} 
	\Ensure{The compressed subspace $\bV_{k_{\rm min}}$}
	\Function{$\bV_{k_{\rm min}} = $ Compress }{$\bV_{k_{\rm max}}, \bR_{\bA}, \bR_{\Psi}, \bd, \bx, \bQ_{\bA} , k_{\rm min}, \lambda$}
 \State $ \bW = \chi\left( \bR_{\bA}, \bR_{\Psi},  \bQ_{\bA}, \bd, \lambda_{\rm curr}\right)$ \Comment{see Section \ref{subsec: basiscompression}}
\State $\widetilde{\bV} = \bV_{k_{\rm max}} \bW$

\State $\bx_{\rm new} = \left(\bx - \widetilde{\bV} \widetilde{\bV}^{T}\bx\right)/\|\left(\bx - \widetilde{\bV} \widetilde{\bV}^{T}\bx\|_2\right)$ \;
\State $\bV_{k_{\rm min}} = [ \widetilde{\bV}, \bx_{\rm new}] $ \Comment{Expand the solution subspace}
\EndFunction
\end{algorithmic}
\end{algorithm}

\begin{algorithm}[!ht]
\caption{RMM-GKS}%
\label{Alg: RMMGKS}
\begin{algorithmic}[1]
	\Require{ $\bA, \Psi, \bV, \bd, \bx^{(0)}, \epsilon$} 
	\Ensure{An approximate solution $\bx^{(i+1)}$ and basis $\bV_{k_{\rm min}}$}
	\Function{$[\bx, \bV_{k_{\rm min}}]= $ RMM-GKS}{$\bA, \Psi, \bV, \bd, \bx^{(0)}, \epsilon$, $k_{\rm min}$, $k_{\rm max}$, $tol_1$} 
    \STATE  $s = k_{\rm max}-k_{\rm min}$\;
     \If{$\bV$, $\bP_{\epsilon}$ provided} 
            \State $\bV_{\ell_{0}} = \bV$
            \State $\bu^{(0)}=\Psi \bx^{(0)}$\;
	        \State	 $\bw^{(0)}_{\epsilon}=\left((\bu^{(0)})^2+\epsilon^2\right)^{-1/2}$\;	
	        \State $\bP_{\epsilon}^{(0)} = (\rm diag(\bw_{\epsilon}^{(0)}))^{1/2}$ \Comment{Compute the weights}
       \Else
	        \State Generate the initial subspace basis $\bV_{\ell_0}\in \R^{n\times \ell_0}$ such that $\bV_{\ell_0}^{T}\bV_{\ell_0}=\bI$\;
             \State $\bu^{(0)}=\Psi \bx^{(0)}$\;
	        \State	 $\bw^{(0)}_{\epsilon}=\left((\bu^{(0)})^2+\epsilon^2\right)^{-1/2}$\;	
	        \State $\bP_{\epsilon}^{(0)} = (\rm diag(\bw_{\epsilon}^{(0)}))^{1/2}$ \Comment{Compute the weights}
   \EndIf
	
	\State $\bA\bV_{\ell_0} = \bQ_{\bA}\bR_{\bA}$ and $\bP_{\epsilon}^{(0)}\Psi \bV_{\ell_0}=\bQ_{\Psi}\bR_{\Psi}$\; \Comment{Compute/Update the QR factorizations}
          \State $\lambda^{(0)}  = \min_{\lambda}\Theta\left(\lambda\right)$ \Comment{Define the regularization parameter}
	\State $\bz^{(1)}=( \bR_{\bA}^T\bR_{\bA} + \lambda^{(0)} \bR_{\Psi}^T\bR_{\Psi})^{-1}\bR_{\bA}^T\bQ_{\bA}^T\bd $\;
	\State $\bx^{(1)}=\bV_{\ell_0}\bz^{(1)}$\;
 \State $\bV_{\ell} \in \R^{n \times \ell}$, $\mathcal{K}_{\ell}\left(\bA^T\bA + \lambda^{(0)}\Psi^T\bP^{(0)}_{\epsilon}\Psi, \bA^T\bd\right)$ 
 \Comment{Generate a better basis}
 \State $\bV_{k_{\rm min}} = \bV_{\ell}$
		\FOR {$i= 1,2,\ldots, i_{\rm max}$}{ 
\State $[\bx^{(i+1)}, \lambda^{(i)}, \bV_{k_{\rm max}}, \bR_{\bA}, \bR_{\Psi}] =$ \Call{Enlarge}{$\bA, \Psi, \bV_{k_{\rm min}}, \bd, \bx^{(i)}, \epsilon$,$s$,$tol_1$}
\State $\bV_{k_{\rm min}} =$ \Call{Compress}{$\bV_{k_{\rm max}}, \bR_{\bA}, \bR_{\Psi}, 
\bd, \bx^{(i+1)}, \bQ_{\bA}, k_{\rm min}, \lambda^{(i)}$}  \Comment{Compress the subspace $\bV_{k_{\rm max}}$ to $\bV_{k_{\rm min}}$}
}\ENDFOR
\EndFunction
\end{algorithmic}
\end{algorithm}

\subsection{Computational Cost}
Aside from the overhead involved in initialization of $\bV_{k_{\rm min}}$, 
the major costs are inside loop, which repeats $j=0, \dots, s-1$, with $s = k_{\rm max} - k_{\rm min}$, times or until convergence, are associated with matrix vector products and reorthogonalizations.  Remember that because the weighting matrix $\bP_{\epsilon}^{(j)}$ changes at every iteration, the QR factorization of $\bP_{\epsilon}^{(j)} \Psi \bV_{k_{\rm min}+j}$ must be computed from scratch for each $j$, whereas the QR factorization of $\bA \bV_{k_{\rm min}+j}$ is computed once for $j=0$ and then is easily updated for each new column added to the solution space.  In summary, at the end of the $s$ steps we have incurred the following major costs: 

\medskip
\noindent{\underline{Matvecs}}:  
   \begin{itemize} 
   \item $s$ with $\bA$.
   \item $s$ with $\Psi$.
  \end{itemize}
\medskip
\noindent{\underline{(Re)orthogonalizations}}:  
   \begin{itemize}
       \item $k_{\rm min}\! + \!s$  to first compute, then update, $\bQ_{\bA}\bR_{\bA}$
       \item  $ \sum_{j=1}^s k_{\rm min}\!+\!j$ to compute $\bQ_\Psi\bR_\Psi$ 
       \item  Reorthogonalize residual to current $\bV_{k_{\rm min}+j}$, $j=1,\ldots,s$.    
\end{itemize}
\medskip
\noindent{\underline{Overhead in Compress Routine}}: 
  \begin{itemize}
      \item Compute matrix $\bW$ using one of the methods in Section \ref{sec: compression} (e.g., method 1 requires truncated SVD of a $2 k_{\rm max} \times k_{\rm max}$ matrix).
      \item Compute product $\bV_{k_{\rm max}} \bW$. 
       \item Orthogonalize current solution against $\widetilde{\bV}$.  
  \end{itemize}
\section{RMM-GKS for Streaming Data (s-RMM-GKS)}
\label{sec: streaming}
In this section, we still assume a linear forward model as in~\eqref{eq:LinSys} but now we assume one of two additional conditions hold: 
  \begin{itemize}
      \item  Either there is only a portion of the data that is available at any given time (i.e. it is streamed) or
      \item The problem is massive  so we can only deal with subsets of rows at a time (i.e. we must treat the problem as if the data is streamed).
      \end{itemize}
      If either of these conditions hold, we will refer to the situation as the {\it streaming data case.}

Ideally, we would like to be able to solve
 \begin{equation}
    \label{eq:rtomoproblemall}
    \min_\bx \left\|{\begin{bmatrix}\bA_1\\ \vdots \\ \bA_{n_t}\end{bmatrix} \bx - \begin{bmatrix} \bd_1 \\ \vdots \\ \bd_{n_t} \end{bmatrix}}\right\|_2^2 + \lambda^2 \|\Psi{\bx}\|_q^q, \quad \text{where} \quad q = 1.
  \end{equation}
However, in the streaming data case, estimating the solution to this large regularized problem is not feasible. 

To adapt RMM-GKS for the streaming data case, we first 
partition $\bA, \bd$.  Then, we imagine solving $n_t$ minimization problems in succession:
\begin{align}
 \label{eq:rtomo1}
  \min_{\bx} \|\bA_1 \bx - \bd_1\|_2^2 & +  \lambda_1 \|\Psi\bx\|_q^q   \\
  \label{eq:rtomoi}
\min_{\bx} \|\bA_2 \bx - \bd_2\|_2^2 & +  \lambda_2 \|\Psi\bx\|_q^q  \\ 
 & \vdots & \nonumber \\
\label{eq:rtomont}
\min_{\bx} \|\bA_{n_t} \bx - \bd_{n_t}\|_2^2 & + \lambda_{n_t} \|\Psi\bx\|_q^q, 
\end{align}
with $q=1$
but {\it constrain each of the solution spaces} of systems $2$ through $n_t$ using RMM-GKS previous information so that $\bx^{(n_t)} \approx \bx,$ where $\bx$ refers to the solution of 
(\ref{eq:rtomoi}).

\comment{
For comparison, we provide the results of some of the methods we consider with automatic regularization parameter selection on the entire problem,
  \begin{equation}
    \label{eq:rtomoproblemall}
    \min_\bx \left\|{\begin{bmatrix}\bA_1\\ \vdots \\ \bA_{n_t}\end{bmatrix} \bx - \begin{bmatrix} \bd_1 \\ \vdots \\ \bd_{n_t} \end{bmatrix}}\right\|_2^2 + \lambda^2 \|\Psi{\bx}\|_q^q, \quad \text{where} \quad q = 1,2.
  \end{equation}
}

Assume that we wish to solve $n_t$ systems (\ref{eq:rtomo1})-(\ref{eq:rtomont}) where $\bA_i$, $\bd_i$, $i = 1,2, \dots, n_t$ arrive in a streaming fashion, i.e., once we have received/processed $\bA_i$ and $\bd_i$ we may not be able to store them (for instance due to memory limitations) before $\bA_{i+1}$, $\bb_{i+1}$ arrive. We can solve \eqref{eq:rtomoi} with RMMGKS and we obtain a solution subspace of dimension $k_{\max}$, called $\bV^{(i)}_{k_{\rm max}}$, and an approximate solution that we denote by $\bx_{(i)}$. When the first system is solved, $\bV^{(0)}_{k_{\rm min}}$ (line 2 of Algorithm \ref{Alg: sRMMGKS}) is typically not given, but computed as in line 8 of Algorithm \ref{Alg: RMMGKS}. 
To avoid keeping redundant basis information, we compress the subspace through the compress procedure described in the previous section to obtain $\widetilde{\bV} := \bV_{k_{\max}}^{(i)} \bW$.
Similar to the non-streamed case, we
include the information about the most recent solution into the updated basis. 
We compute 
$$\bar{\bx} = \left(\bx_{(i)} - \widetilde{\bV} \widetilde{\bV}^{T}\bx_{(i)}\right)/\|\left(\bx_{(i)} - \widetilde{\bV} \widetilde{\bV}^{T}\bx_{(i)})\|_2\right),$$
and initialize $\bV^{(i+1)}_{k_{\rm min}} = [\widetilde{\bV},\bar{\bx}]$.
We discard all other information except for $\bV_{k_{\rm min}}^{(i+1)}$ and $\bx_{(i)}$. 
We use now RMM-GKS to solve the $(i+1)$-th system where $\bx_{(i)}$ and $\bV_{k_{\rm min}^{(i+1)}}$ are used as initial approximate solution and initial solution subspace, respectively. The process is repeated for the remaining systems. We summarize the process of solving a streaming problem in Algorithm~\ref{Alg: sRMMGKS}. We illustrate the performance of s-RMM-GKS in the numerical examples section.

\begin{algorithm}[!ht]
\caption{s-RMM-GKS}%
 \label{Alg: sRMMGKS}
\begin{algorithmic}[1]
	\Require{${\bA_1}, {\bA_2}, \cdots, \bA_{n_t}$, $\Psi$, $\bd_1, \bd_2, \cdots, \bd_{n_t}$, $\bx^{(0)}$, $\bepsilon$, $ k_{\rm min}$, $k_{\rm max}$, $tol_1$} 
	\Ensure{An approximate solution $\bx$ of the streaming problem }
	\Function{$\bx = $ s-RMM-GKS}{${\bA_1}, \cdots, \bA_{n_t}, \Psi, \bd_1, \cdots, \bd_{n_t}, \bx^{(0)}, \bepsilon, k_{\rm min}, k_{\rm max}, tol_1$}
	\State $[\bx_{(1)}, \bV^{(1)}_{k_{\rm min}}] =$ RMM-GKS($\bA_1, \Psi, \bV^{(0)}_{k_{\rm min}}, \bd_1, \bx^{(0)}, \bepsilon, k_{\rm min}, k_{\rm max}, tol_1$)
  \For{$j=2$ to $n_t$}
       \State $[\bx_{(j)}, \bV^{(j)}_{k_{\rm min}}] =$ RMM-GKS($\bA_{j}, \Psi, \bV^{(j-1)}_{k_{\rm min}}, \bd_{j}, \bx_{(j-1)}, \bepsilon,  k_{\rm min}, k_{\rm max}, tol_1$) 
  \EndFor
\State $\bx = \bx_{(n_t)}$
	\EndFunction
\end{algorithmic}
\end{algorithm}

\begin{remark}
   Throughout the paper, for simplicity we used $q = 1$. It is straightforward to extend them to a general $0<q\leq2$. In particular, to obtain a general $\ell_q$ constrained RMM-GKS method, we substitute line 8 of Algorithm \ref{Alg: RMMGKS} with $\bw^{(0)}_{p,\epsilon}=\left((\bu^{(0)})^2+\epsilon^2\right)^{p/2-1}$ and line 4 of Algorithm \ref{Alg: Enlarge} with $\bw^{(k)}_{p,\epsilon}=\left((\bu^{(k)})^2+\epsilon^2\right)^{p/2-1}$.
\end{remark}

\section{Compression approaches}\label{sec: compression}
\comment{
Compression allows reducing the total number of the solution vectors without significantly affecting the accuracy of the resulting reconstructed solution. Compression is crucial when the memory capacity is reached (we can not store more basis vectors) without the method converging to the desired solution. In particular, once the memory capacity of storing $k_{\rm max}$ vectors is reached by storing the current set of basis vectors $\bV_{k_{\rm max}}$, we seek to compress the subspace to $\bV_{k_{\rm min}-1}$ through $\bW_{k_{\rm min}-1}$ and by augmenting one additional basis vector into $\bV_{k_{\rm min}-1}$ obtained by an approximate solution at the previous iteration (see for instance lines 3 and 4 in Algorithm \ref{Alg: Compress}). 
}
Previously, the discussion of our recycling algorithm has been agnostic to the routine, $\chi$ used to return the mixing matrix, $\bW$, given the inputs $\bR_{\bA}, \bR_{\Psi},  \bQ_{\bA}, \bd, \lambda_{\rm curr}.$
In this section, we describe four different strategies for constructing the $\bW$ we need
in Algorithm \ref{Alg: Compress}.  
These choices are specific to the fact that we are solving a regularized-projected ill-posed problem.
The first two strategies -- truncated SVD and reduced basis decomposition (RBD) -- rely on properties of the stacked matrix
\begin{equation}\label{eq: SVDStacked}
\bar{\bH}_{k_{\rm max}} = \begin{bmatrix} \bR_{\bA} \\ \sqrt{\lambda_{\rm curr}} \bR_{\Psi} \end{bmatrix}.
\end{equation}

The other two methods rely on properties of a regularized solution to the small projected problem. We call these techniques 
solution-oriented compression and sparsity enforcing compression, respectively.

\paragraph{\bf 1. Truncated SVD (tSVD)}
 We compute the truncated SVD of 
 $\bar{\bH}_{k_{\rm max}}$, truncating to 
 $k_{\rm min}-1$ terms.  

That is, 
$\bar{\bH}_{k_{\rm max}}  \approx \bU \bS \bW^{T},$
where $\bU \in \R^{2k_{\rm max}\times k_{\rm min}-1}$, $\bW\in \R^{k_{\rm max}\times k_{\rm min}-1}$ are matrices with orthogonal columns  
and $\bS\in \R^{k_{\rm min}-1\times k_{\rm min}-1}$ is a diagonal matrix that contains the largest $k_{\rm min}-1$ singular values $\sigma_i$ of the matrix $\bar{\bH}_{k_{\rm max}}$. We then return $\bW$.

\paragraph{\bf 2. Reduced basis decomposition (RBD)} 
As a second compression approach we explore the reduced basis decomposition from reduced order modeling \cite{chen2015reduced} to obtain a compressed representation of a data matrix.
Used in the context of our problem, this is a greedy strategy that determines an approximate factorization for $\bar{\bH}_{k_{\rm max}}^T$ as  
\begin{equation}\label{eq: RBD}
\bar{\bH}_{k_{\rm max}}^T \approx \bW_{j}  \bT_{j},
\end{equation}
where the matrix $\bW_{j} \in \R^{k_{\rm max}\times j}$ has orthonormal columns and $\bT \in \R^{j \times 2k_{\rm max}}$ is the transformation matrix. 
The RBD alogrithm needs a tolerance and a max dimension in order to return the approximate factorization. Let $1\leq j \leq d$, where $d$ is the largest number of basis vectors we wish to keep after compression and we define $$\epsilon_j = \rm max_{1\leq k \leq 2k_{\rm max}}\|\bar{\bH}^T_{k_{\rm max}}\left(:, k\right) - \bW_{j} \bT_{j}\left(:, k\right)\|_2.$$ We let $k_{\rm min} - 1 = j$ if $\epsilon_d < \epsilon_{tol}$ with $j$ being the largest index such that $\epsilon_{j} \geq \epsilon_{tol}$. Otherwise, we set $k_{\rm min} - 1 = d$.
The chosen $\bW$ will be returned as the output
to $\chi$ .

\paragraph{\bf 3. Solution-oriented compression (SOC)}

This technique involves the solution of the regularized projected problem \eqref{eq: minKryov2}.  $\bW$ will be determined to be columns of an identity matrix, so that when multiplying against $\bV_{\rm max}$, from the left, only certain columns of $\bV_{\rm max}$ will remain.  

Let $\bz = [z_1, z_2, \dots, z_{k_{\rm max}}]$ be the solution of the projected problem obtained as in \eqref{eq: minKryov2}. 

Define the sets of indexes such that
\begin{equation}\label{eq: Im}
I = \{i: |z_i|>tol\}
\end{equation}
\begin{equation}\label{eq: Jm}
J = \{i: |z_i|\quad  \text{are the largest} \quad d \quad \text{components}\}
\end{equation}
Choose the index set $K = \{k_1, k_2, \dots, k_{\rm min}-1\} \subseteq \left(I \cap J\right)$. Then if $\bf I$ is the identity matrix of size $k_{\rm max}$, we let $\bW = {\bf I}_{:,K}$. 

\paragraph{\bf 4. Sparsity-enforcing compression (SEC)}

This technique differs from the last only in the way the projected problem is regularized.  We solve an auxiliary problem  

\begin{equation}\label{eq: Sparsity}
\bz^{**} = \arg\min_{\bz \in \R^{k_{\rm max}}}\|\bR_{\bA}\bz - \bQ^T \bd\|^2_2 + \lambda\|\bR_{\Psi}\bz\|^1_1.
\end{equation}
which will enforce sparsity of the solution $\bz$.  We solve this small regularized problem for $\bz$ by MM \cite{lange2016mm}.
Note that we use this $\bz$ only inside the compression function to select the index sets in (\ref{eq: Im}) and (\ref{eq: Jm}), and not elsewhere in the recycling
algorithm as a whole. 
 
A detailed numerical comparison of the compression approaches is presented in the first numerical example in Section \ref{sec: NumericalExperiments}.

\section{Numerical Results} \label{sec: NumericalExperiments}
We illustrate the performance of RMM-GKS and s-RMM-GKS and compare their performance with existing methods on applications in  
image deblurring, dynamic photoacoustic tomography, and computerized tomography. An illustration of RMM-GKS applied on real CT data is given in Appendix~\ref{app:mm}.
In every scenario on each application, we observe that our recycling based approaches
are successful in providing high quality reconstructions with very limited memory requirements. All computations were carried out in MATLAB R2020b with about 15 significant decimal digits running on a laptop computer with core CPU Intel(R) Core(TM)i7-8750H @2.20GHz with 16GB of RAM. When the size of the problem allows for us to compare to full MM-GKS, our reconstructions are improved over the MM-GKS reconstructions.  In other scenarios, the problem is too large to store the MM-GKS iterates beyond a certain iteration, while our RMM-GKS can operate under the memory limitations and provide a high-quality solution. 

\comment{
\paragraph{Discussion on the choice of numerical examples}
In our first set of numerical experiments, we compare reconstruction results between MM-GKS and our RMM-GKS on an image deblurring problem. Since RMM-GKS performance will depend on the choice of compression routine we use, we also provide comparisons of the four techniques from Section \ref{sec: compression}. 
We next consider two scenarios from computerized tomography where multiple linear systems are aimed to solve to reconstruct a medium of interest. The focus of the first scenario is on comparing other existing hybrid methods (with or without recycling) with our RMM-GKS. Such example illustrates that including edge information on the recycled subspace helps to enhance the reconstruction quality (see for instance the comparison between HyBR and MM-GKS methods). Further, we test the robustness of the RMM-GKS with respect to noise level. The second scenario considers a larger test problem where the focus is on illustrating the performance of RMM-GKS and on dynamically selecting the number of iterations on the RMM-GKS. In addition, we test the performance of our s-RMM-GKS when the data are streamed randomly. 
Our third example displays a large-scale example from dynamic photoacoustic tomography where our RMM-GKS shows promising results by keeping the memory limited while MM-GKS can not converge due to reaching the memory capacity. We further compare versus a recenlty proposed restarted MM-GKS that we denote by MM-GKS$_{\rm res}$ \cite{restartMMGKS}. For all the above test problems we have the true solution which is used to display reconstruction quality measures that we discuss below. In the last example displayed in the Supplementary Materials file we treat a large-scale computerized tomography and dynamic problem with real data.
}

For all the examples we perturb the measurements with white Gaussian noise, i.e., the noise vector $\be$ has mean zero and a rescaled identity covariance matrix; 
we refer to the ratio $\sigma=\|\be\|_{2}/\|\bA\bx\|_{2}$ as the noise level. 

To assess the quality of the reconstructed solution we consider several quality measures. We compute the Relative Reconstruction Errors (RRE). That is, for some recovered $\bx^{(k)}$ at the $k$-th iteration, the RRE is defined as
\begin{equation} \label{eq: RRE} 
{\rm RRE}:={\rm RRE}(\bx^{(k)},\bx_{\mathrm{true}}) = \frac{||\bx^{(k)}-\bx_{\mathrm{true}}||_2}
{\|\bx_{\mathrm{true}}\|_2}. 
\end{equation}
In addition, we use Structural SIMilarity index (SSIM) between $\bx^{(k)}$ and $\bx_{\mathrm{true}}$ to measure the quality of the computed approximate solutions. The definition of the SSIM is involved and we refer to \cite{SSIM} for details. Here we recall that the SSIM measures how well the overall structure of the image is recovered; the higher the index, the better the reconstruction. The highest achievable value is $1$. 
Our work is concerned with reconstructing solutions that preserve the edges in the medium of interest, hence we use HaarPSI \footnote{Throughout the paper we use HP when we report the HaarPSI measure on Tables.}, a recently proposed measure that stands for Haar wavelet-based perceptual similarity index \cite{reisenhofer2018haar}.

The iterations in the enlarge routine (see Algorithm ~\ref{Alg: Enlarge}), unless otherwise stated, are terminated as soon as the maximum number of iterations is reached or $\rm T_1$$ \leq tol_1$, where
\begin{equation} \label{eq: TOL1} 
{\rm T_1}:={\rm T_1}(\bx^{(k+1)},\bx^{(k)}) = \frac{||\bx^{(k)} - \bx^{(k-1)}||_2}
{||\bx^{(k-1)}||_2}, k = 1, \cdots.
\end{equation}
We choose $tol_1 = 10^{-3}$.

\subsection{Image deblurring}
\begin{figure}[h!]
\begin{center}
\begin{minipage}{0.32\textwidth}
		\includegraphics[width=\textwidth]{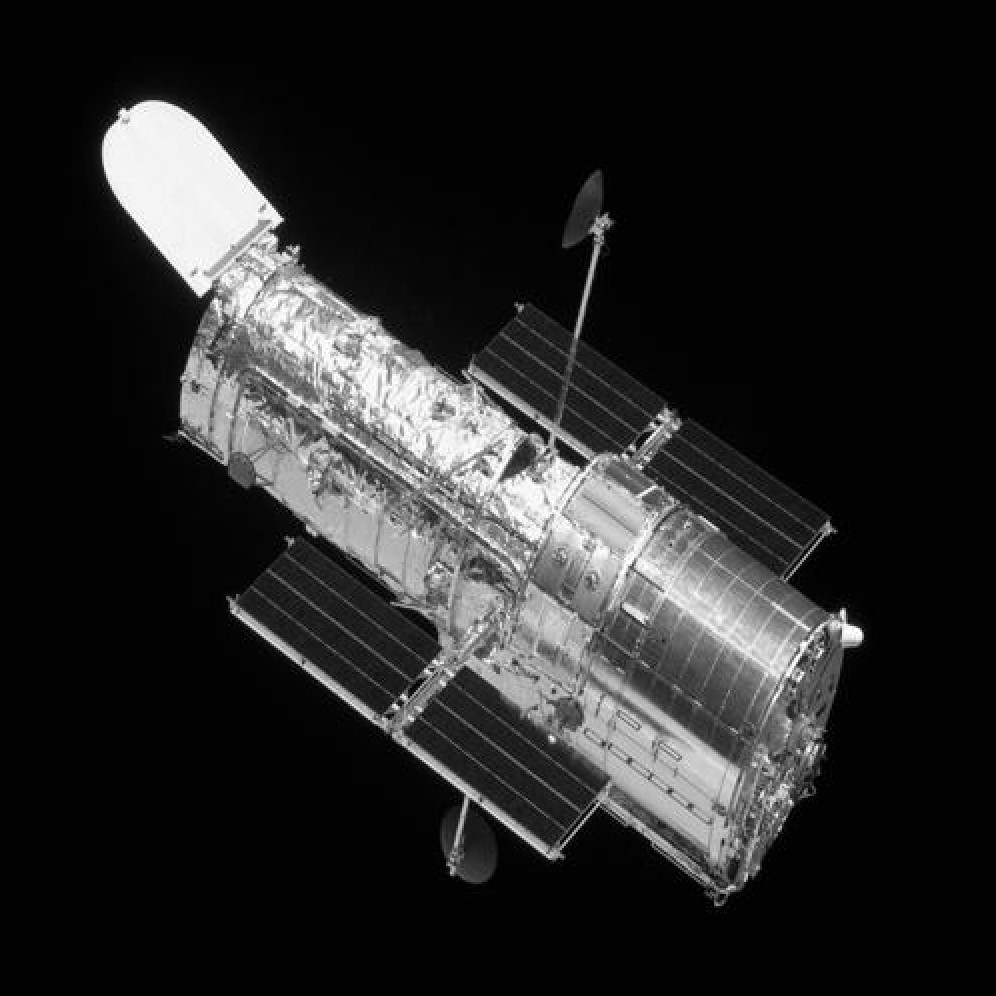}
	\end{minipage}
	\begin{minipage}{0.32\textwidth}
		\includegraphics[width=\textwidth]{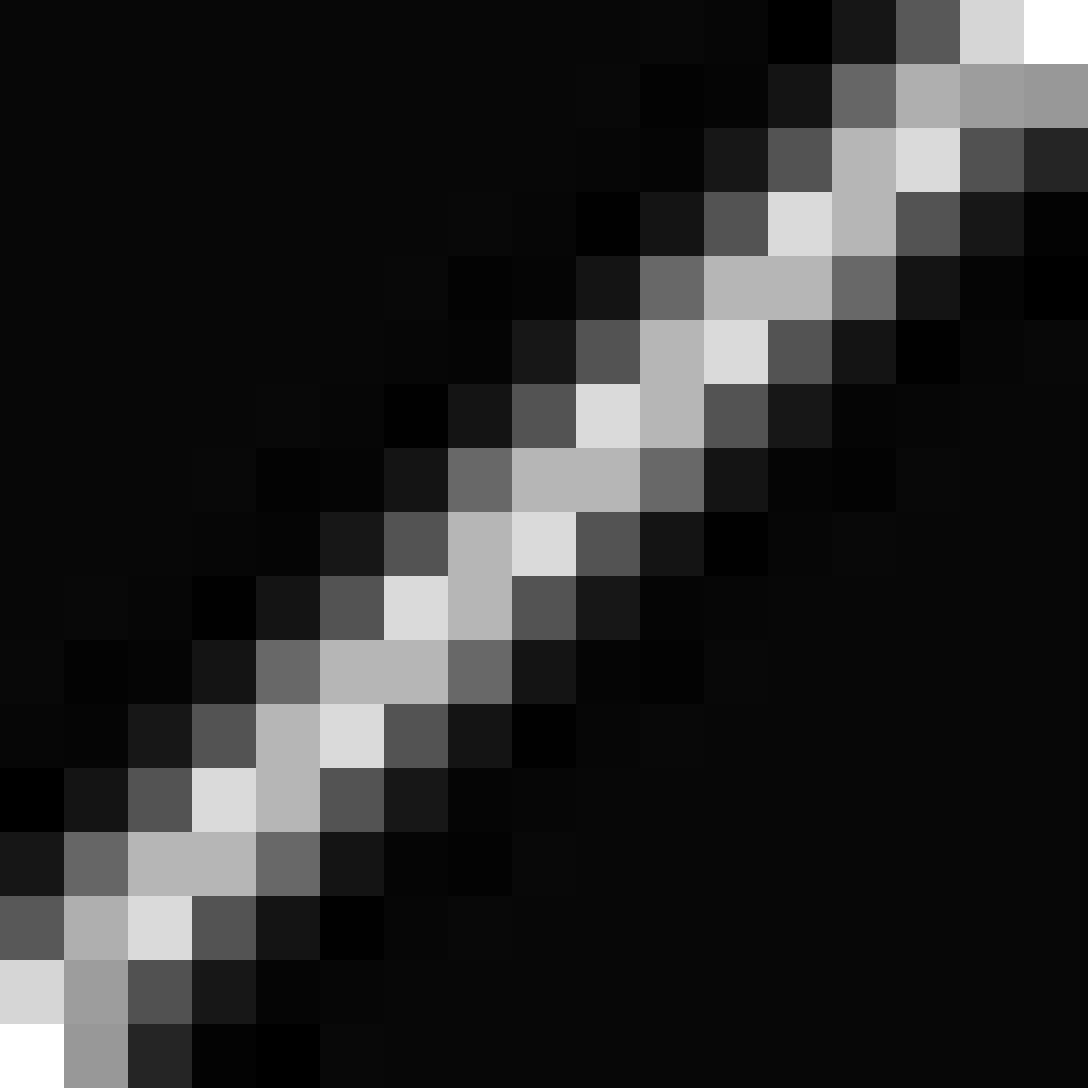}
	\end{minipage}
	\begin{minipage}{0.32\textwidth}
		\includegraphics[width=\textwidth]{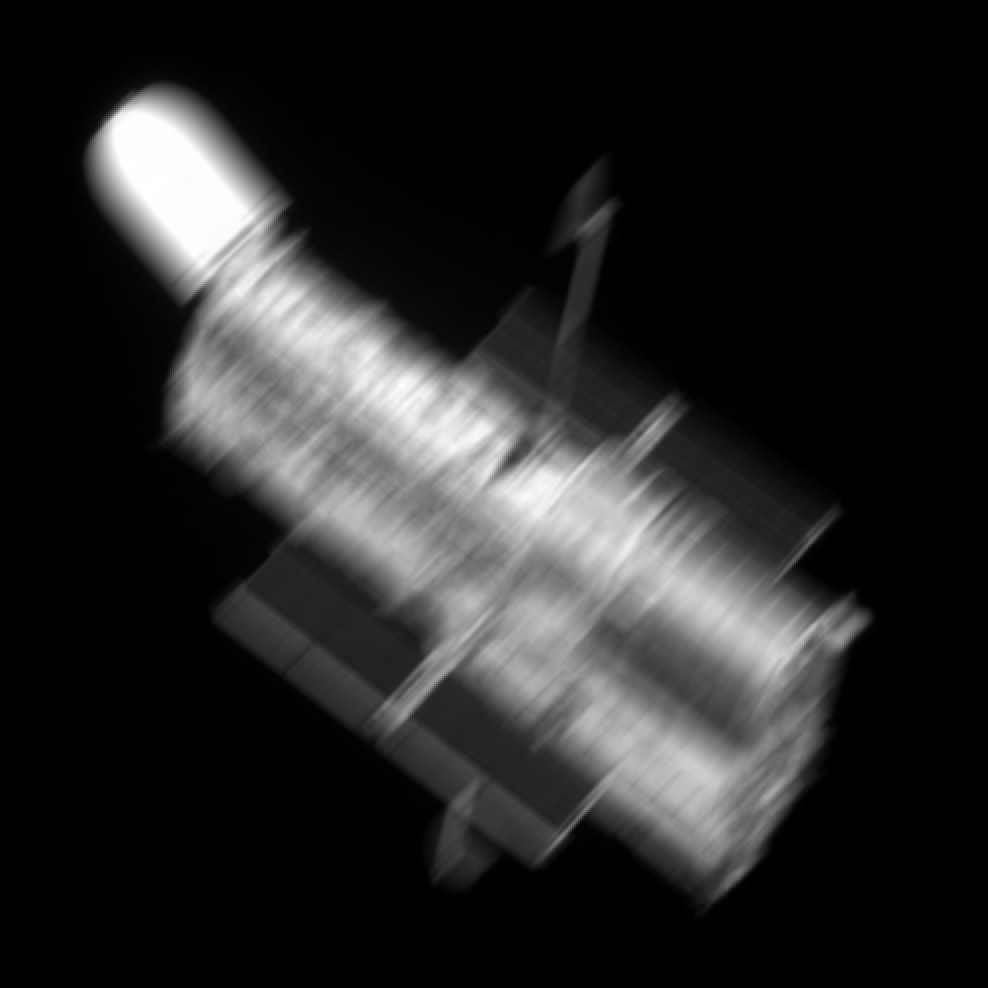}
	\end{minipage}
 \end{center}
\caption{Telescope test problem. (a) True image of $500 \times 500$ pixels. (b) Motion PSF ($17 \times 17$ pixels), (c) Blurred and noisy image with $0.1\%$ Gaussian noise.}
\label{Fig: telescopeImages}
\end{figure}

\paragraph{Telescope test problem}
Here we consider an  
image deblurring problem where the image has been corrupted by motion blur. The aim is to reconstruct an approximation of the true, $500 \times 500$ telescope image shown in Figure \ref{Fig: telescopeImages} given the observed blurred and noisy image with 0.1$\%$ Gaussian noise. 
The motion point-spread function (PSF) of size $17\times 17$ pixels 
and the blurred and noisy telescope image are shown in Figure \ref{Fig: telescopeImages}(b) and \ref{Fig: telescopeImages}(c), respectively. 
The PSF implicitly determines the $250,000 \times 250,000$ blurring operator, $\bA$.  

The goals are 1) to investigate the reconstruction quality when several compression approaches are used, 2) to compare the results with the MM-GKS when the memory capacity is limited, 3) to investigate several choices on the number of solution basis vectors to be kept in the memory after compression. To achieve these goals, we set the memory capacity $k_{max}=25$, varied $k_{min}$ from 5 to 15 and computed the results for each of the four compression routines (also, the compression tolerance on RBD was set to $\epsilon_{tol} = 10^{-5}$ and for SOC the tolerance was set to $tol = 1$). 
All methods are run until the maximum number of iterations (200) is reached or the relative error of two consecutive reconstructions falls below a tolerance $\epsilon_{rel} = 10^{-5}$.

\begin{figure}[ht!]
\centering
  \begin{tabular}{cccc}
   MM-GKS (25) &  TSVD &  RBD & SEC \\
  \includegraphics[height = 0.18\textwidth, width = .2\textwidth]{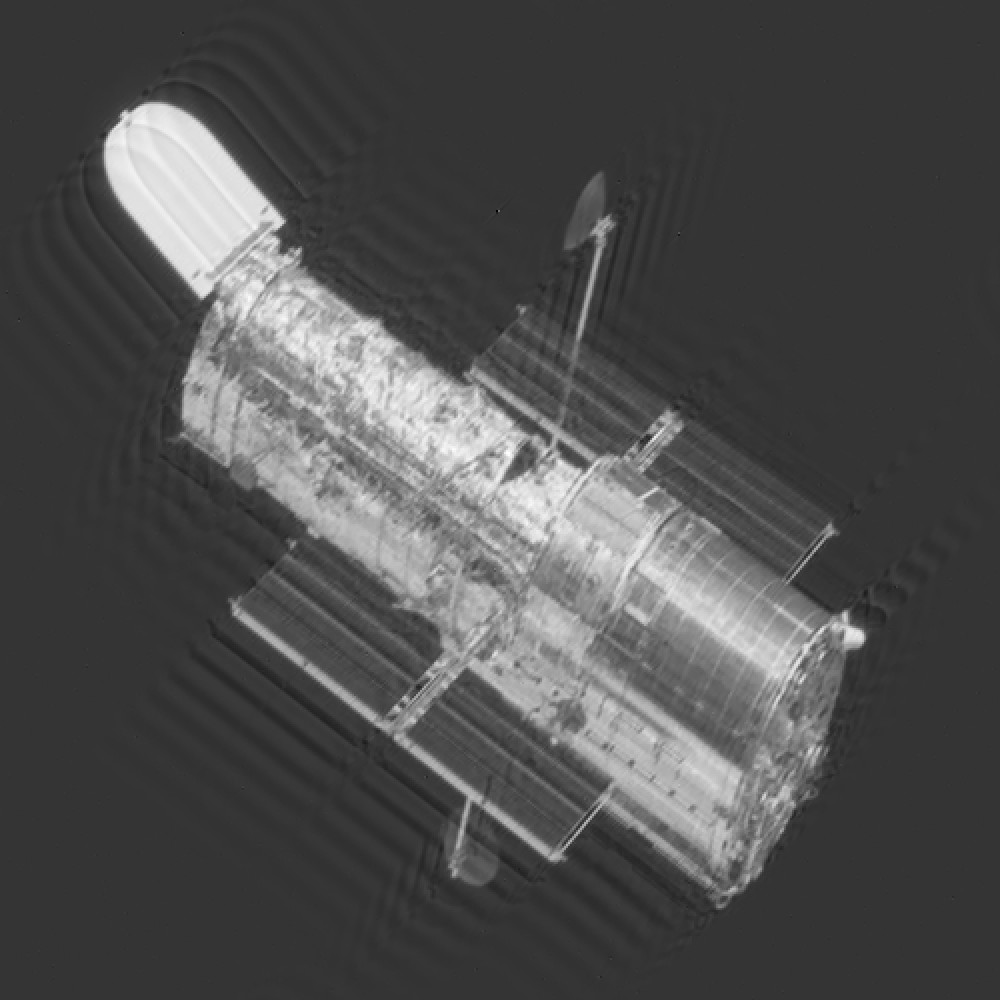} &
  \includegraphics[height = 0.18\textwidth, width = .2\textwidth]{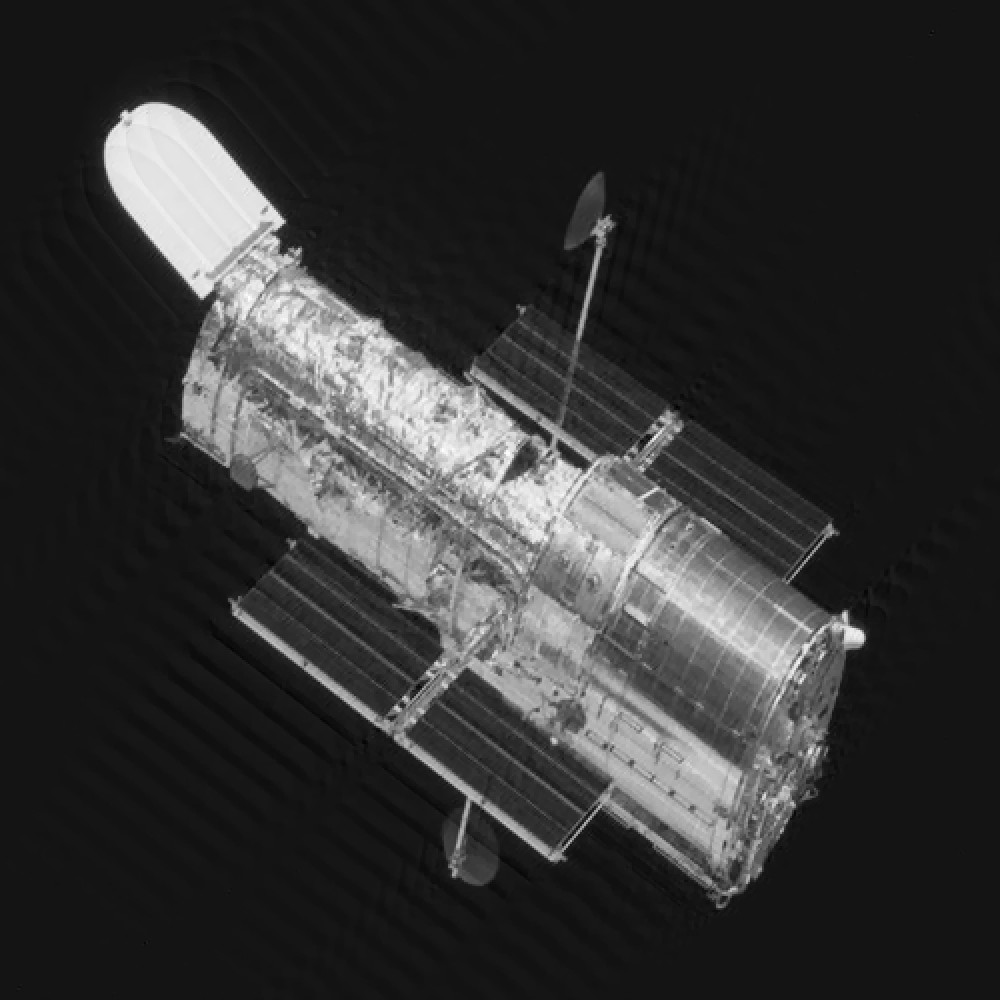} &
  \includegraphics[height = 0.18\textwidth, width = .2\textwidth]{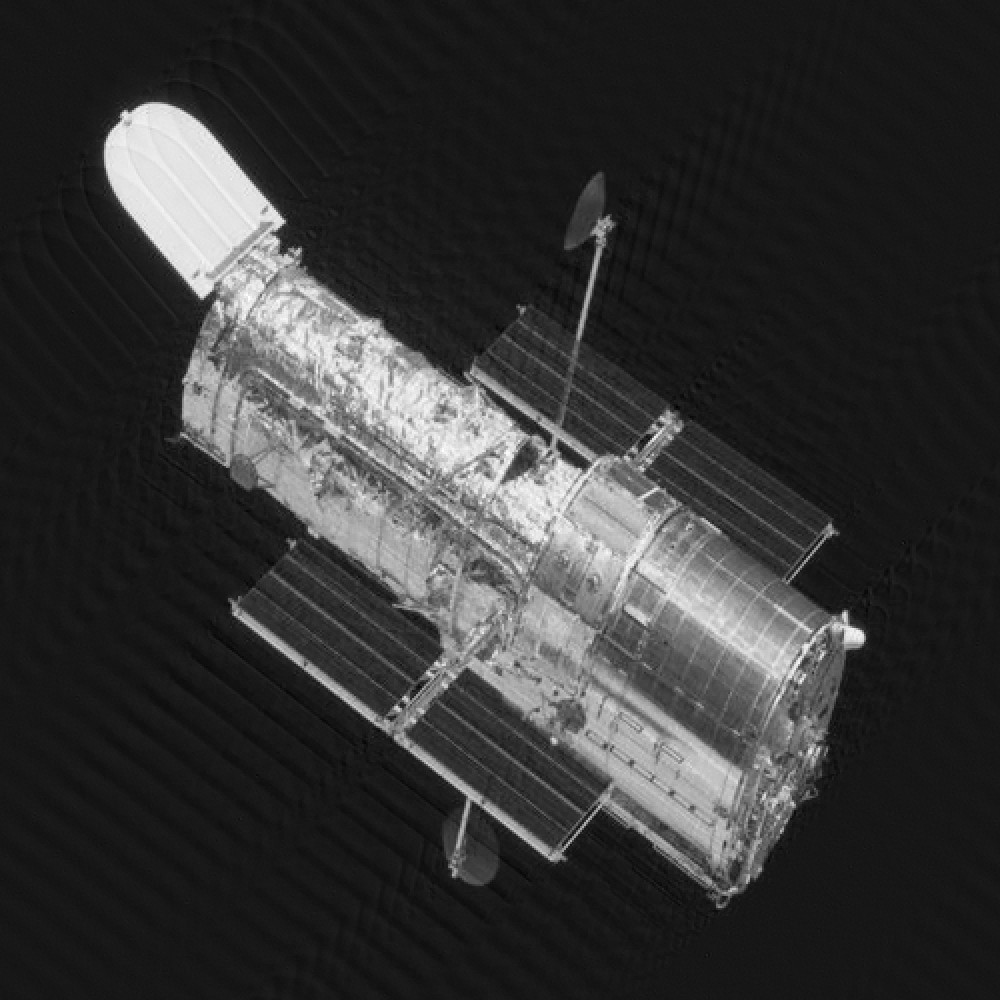}& 
    \includegraphics[height = 0.18\textwidth, width = .2\textwidth]{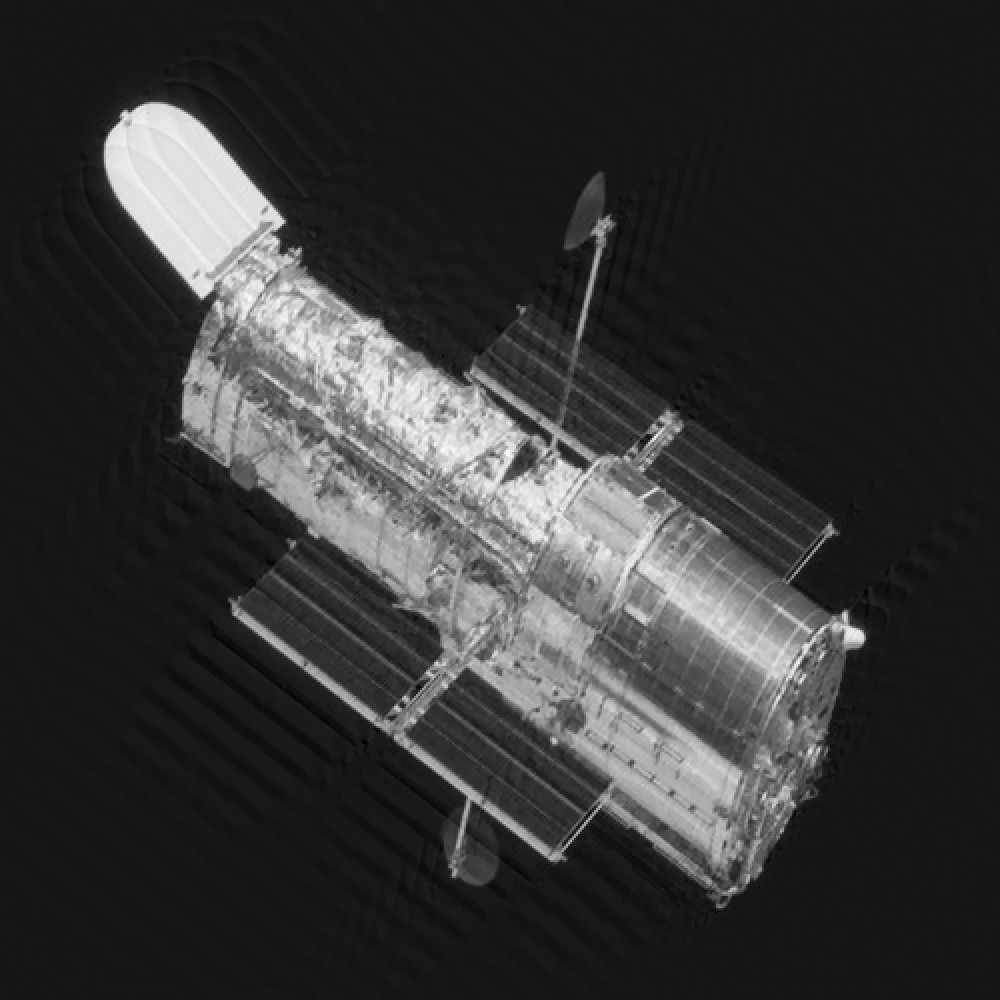} 
\end{tabular}
\begin{tabular}{cccc}
  \includegraphics[height = 0.18\textwidth,  width = .2\textwidth]{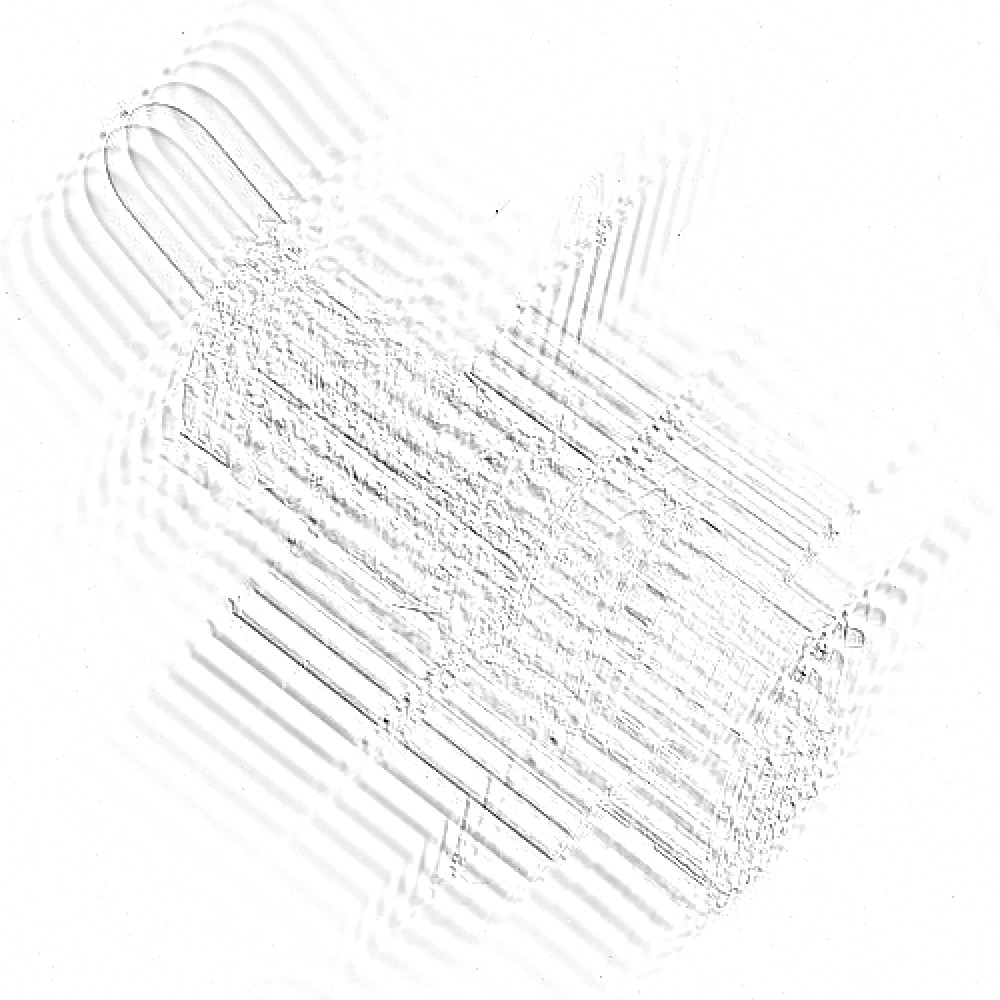}
    &
  \includegraphics[height = 0.18\textwidth,  width = .2\textwidth]{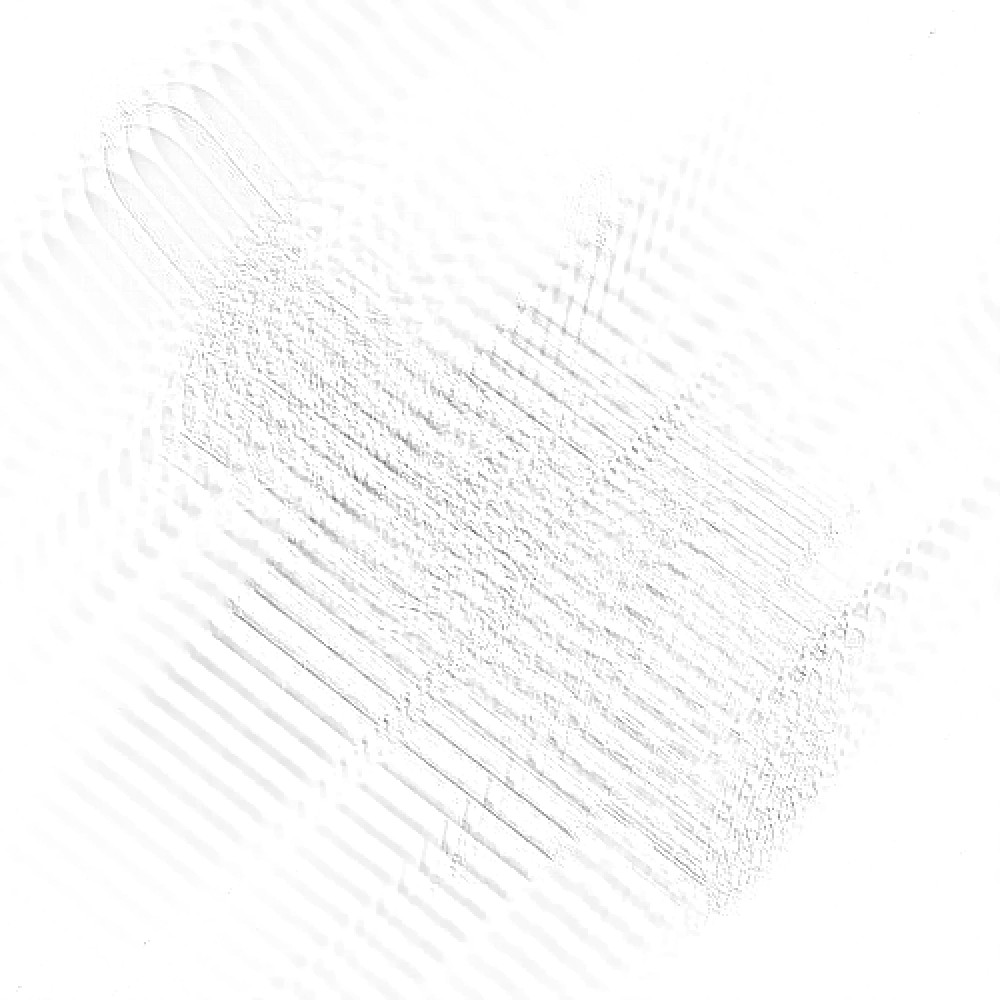}
    &
  \includegraphics[height = 0.18\textwidth, width = .2\textwidth]{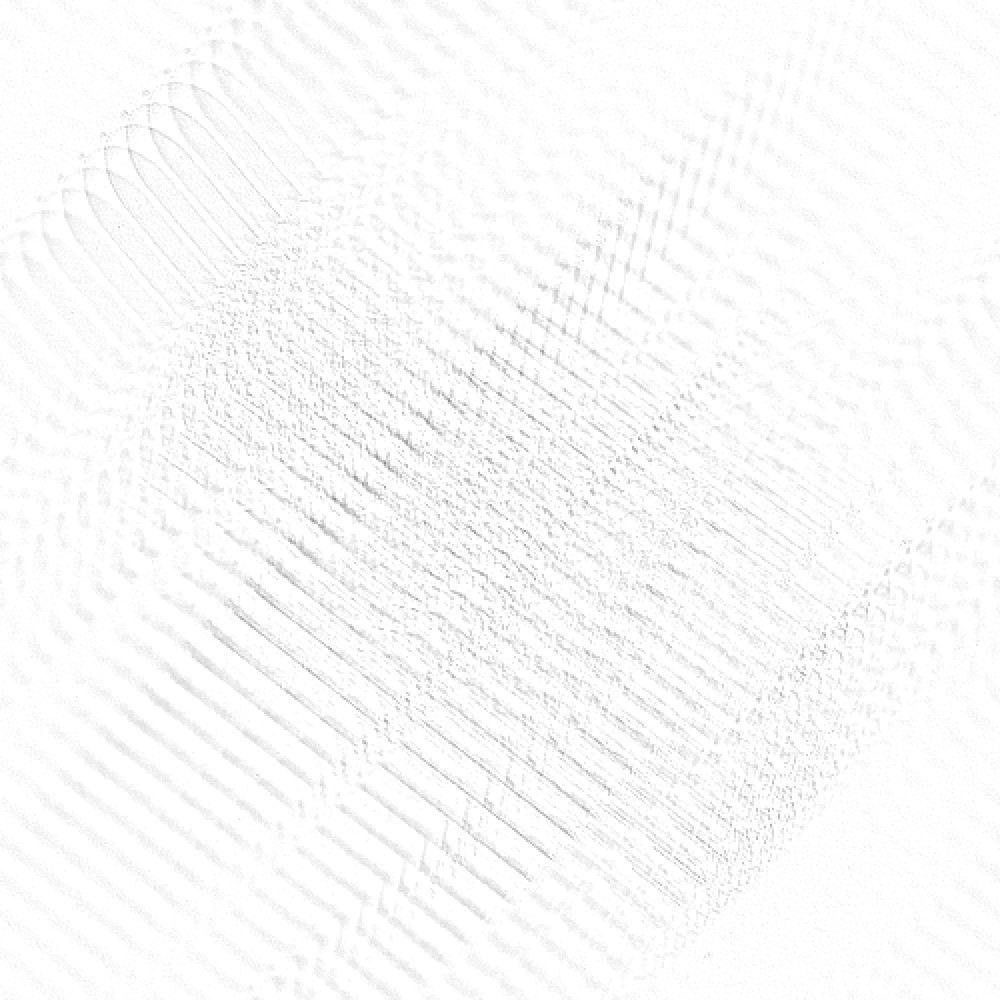}
     &
  \includegraphics[height = 0.18\textwidth, width = .2\textwidth]{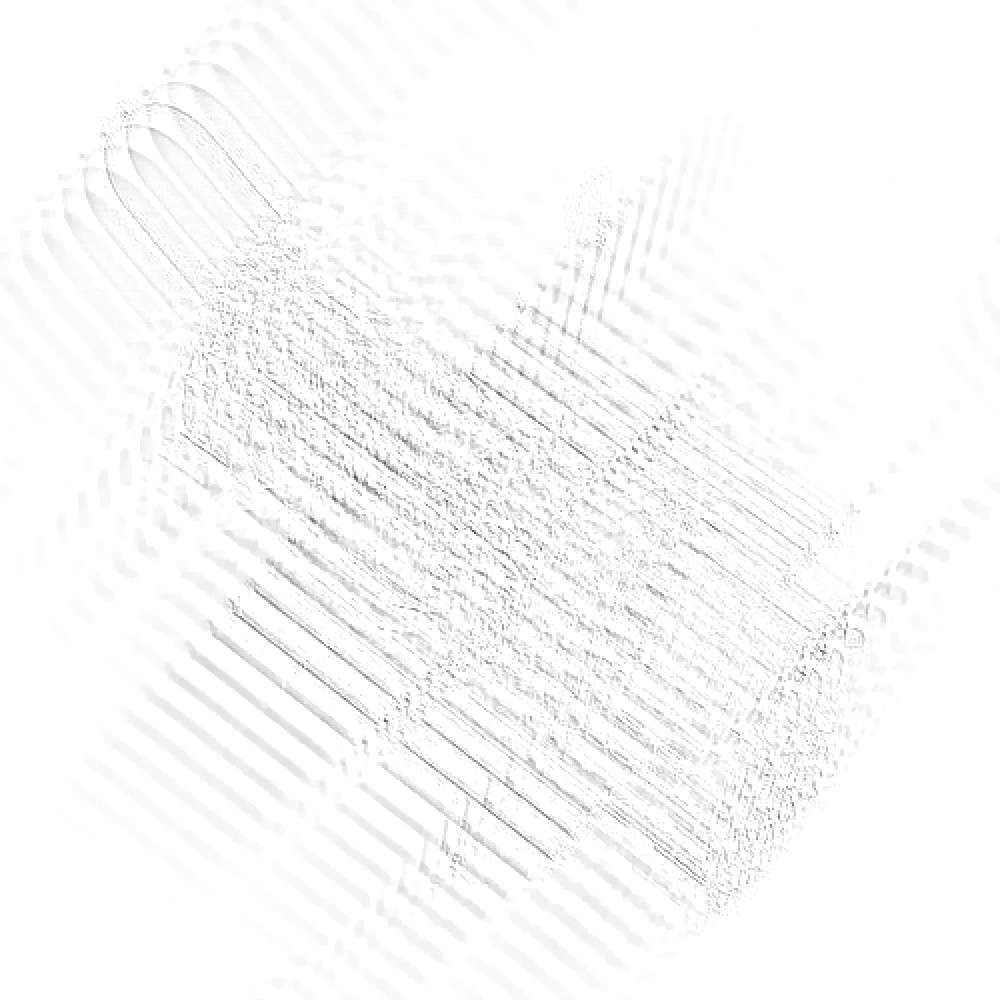}
\end{tabular}
\caption{Telescope test problem: The first row shows the reconstructed images and the second row shows the error images in the inverted colormap by MM-GKS (25 iterations), and RMM-GKS with compression approach being TSVD, RBD, and SEC, respectively, from left to right.}
\label{Figure: Telescope_errors}.
\end{figure}

Results for RMM-GKS with all the compression approaches we consider are shown in Table \ref{Table: RRE_telescope}. 
 For comparison, running MM-GKS up to the hypothetical memory capacity (i.e. the solution subspace maximum was 25) yields an RRE = 0.106 and HaarPSI = 0.949.
 As the results in the table show, in every case the RRE and HaarPSI scores of our RMM-GKS are much improved over the MM-GKS values. We show the reconstructed images along with the error images in the reverted map in the first and second rows of Figure \ref{Figure: Telescope_errors}.

\begin{table}[ht!]
\centering
\fontsize{10pt}{12pt}\selectfont
\begin{tabular}{|l|l|l|l|l|l|l|l|l|}
\hline
\hline
\multirow{2}{*}{$k_{\rm min}$} & \multicolumn{2}{c|}{TSVD}                                             & \multicolumn{2}{c|}{RBD}                           & \multicolumn{2}{c|}{SOC}                           & \multicolumn{2}{c|}{SEC}         \\ \cline{2-9} 
                            & \multicolumn{1}{l|}{RRE} & \multicolumn{1}{l|}{HP} & \multicolumn{1}{l|}{RRE} & \multicolumn{1}{l|}{HP} & \multicolumn{1}{l|}{RRE} & \multicolumn{1}{l|}{HP} & \multicolumn{1}{l|}{RRE} & HP    \\ \hline 
5                           & 0.064                    & 0.981                                     & 0.056                    & 0.984                   & 0.074                    & 0.974                   & 0.074                    & 0.975 \\ \hline
10                          & 0.071                    & 0.977                                   & 0.067                    & 0.977                   & 0.077                    & 0.970                   & 0.077                    & 0.972 \\ \hline
15                          & 0.072                    & 0.975                             & 0.065                    & 0.978                   & 0.081                    & 0.970                   & 0.082                    & 0.925 \\ \hline \hline
\end{tabular}
   \caption{Telescope test problem: RRE and HaarPSI for $k_{\rm min}$ = 5, 10, 15. Algorithm setting: Maximum memory limit is $k_{\rm max}= 25$, Image: Telescope. Size:  $500\times 500$ pixels. All methods are run until the maximum number of iterations (200) is reached or the relative error of two consecutive reconstructions falls below a tolerance $\epsilon_{\rm rel} = 10^{-5}$.}
   \label{Table: RRE_telescope}
\end{table}

\subsection{Computerized tomography}
Our goal in this section is to test the streaming version of RMM-GKS described in Section \ref{sec: streaming}.  Recall that in the ideal 
case, we would want to solve (\ref{eq:rtomoproblemall}). 

We consider two scenarios when this is not possible: in the first set of experiments (labeled Test 1), we assume that this cannot be solved because only a fraction of the sinogram data is available for processing at any given time. In the second set of experiments (labeled Test 2), we mimic the scenario where the data has all been collected but the system is too large to fit in memory. Thus, we use s-RMM-GKS to approximately solve instead the $n_t$
minimization problems (\ref{eq:rtomo1})-(\ref{eq:rtomont}).  

\subsubsection{Test 1} 
\begin{figure}[htbp]
\centering
  \begin{tabular}{cccc}
$\bx_{\rm true}$ &  $\bd_1$: $0^{\circ}-44^{\circ}$ &  $\bd_2$: $45^{\circ}-89^{\circ}$ & $\bd_2$: $90^{\circ}-179^{\circ}$\\
  \includegraphics[height = 0.15\textwidth, width = .18\textwidth]{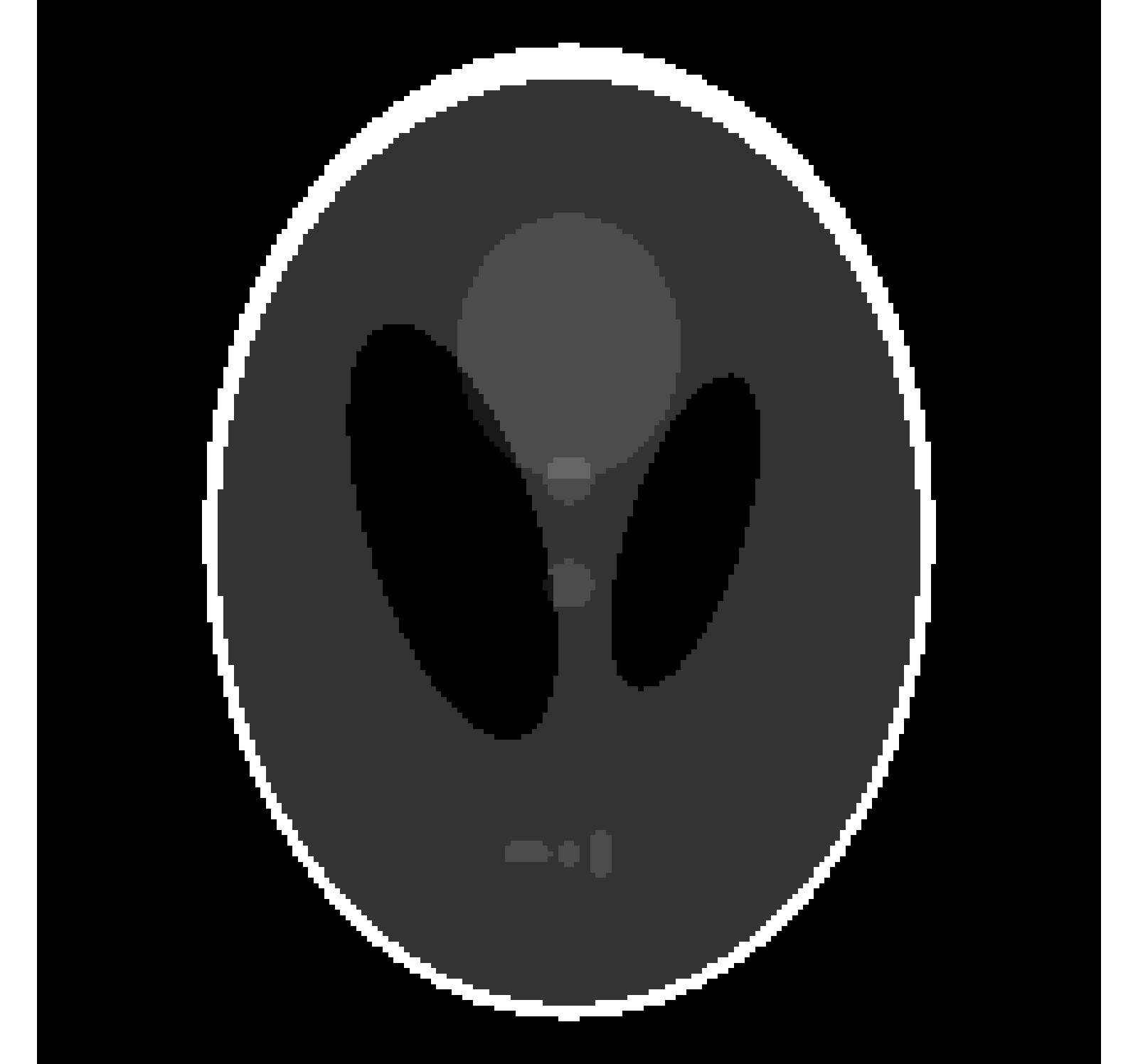} &
  \includegraphics[angle = 90,height = 0.1\textwidth, width = .22\textwidth]{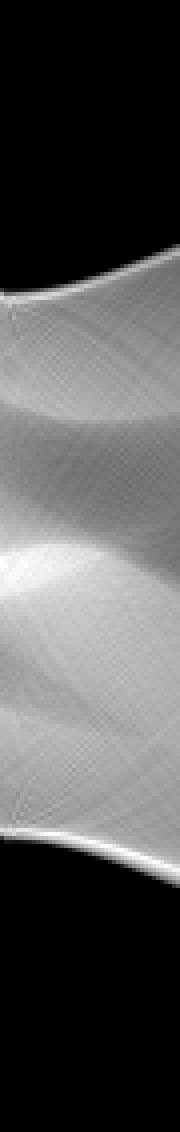} &
  \includegraphics[angle = 90,height = 0.1\textwidth, width = .22\textwidth]{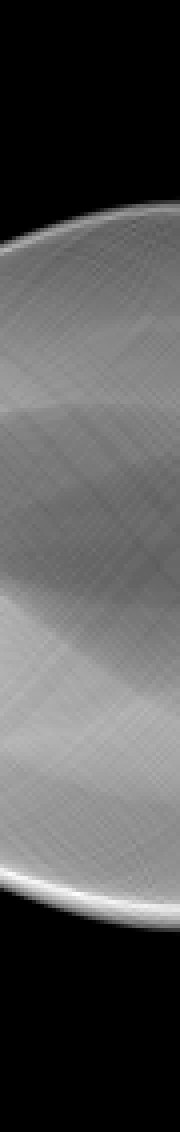} 
   &
  \includegraphics[angle = 90,height = 0.1\textwidth, width = .22\textwidth]{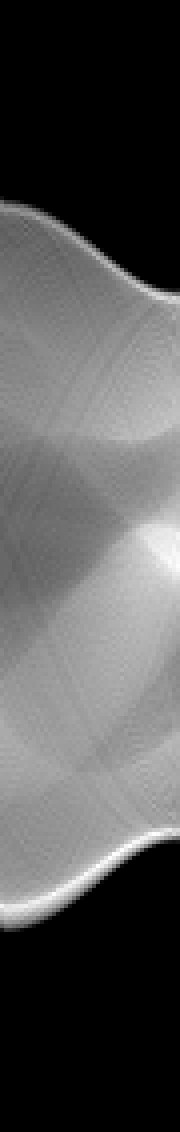}\\
(a) & (b)  & (c)  & (d)
\end{tabular}
  \caption{Test 1: Streaming tomography example. The true image is provided in (a), along with three observed sinograms (b) $\bd_1$, (c) $\bd_2$, and (d) $\bd_3$ corresponding to projections taken at intervals from $\bd_1$: $0^{\circ}-44^{\circ}$,  $\bd_2$: $45^{\circ}-89^{\circ}$, and $90^{\circ}-179^{\circ}$ respectively.}
  \label{fig: Streaming3prob_true}
\end{figure}

For the first experiment we consider the parallel tomography example from IRTools package \cite{gazzola2018ir} where we generate the true phantom to be a Shepp-Logan phantom of size $n_x, n_y = 500$ resulting in $\bx_{\rm true} \in \R^{500^2}$. 
Assume that during the data acquisition process, the data is streamed in three groups, so that $n_t = 3$ in \eqref{eq:rtomont}. The first and the second problems correspond to projection angles from $0^{\circ}-44^{\circ}$ and  $45^{\circ}-89^{\circ}$, respectively, with angle gap $1^{\circ}$. The third problem corresponds to data from $45$ angles from  $90^{\circ}$ to $179^{\circ}$ with angle gap $2^{\circ}$. This setup produces forward operators and the observations $\bA_{1}, \bA_{2}, \bA_3 \in \R^{ 45 \cdot 707 \times 500^2}$ and $\bd_1, \bd_2, \bd_3 \in \R^{45\cdot 707}$, respectively. We artificially add $0.1\%$ noise to each system.
The sinograms and the true image are provided in Figure \ref{fig: Streaming3prob_true}. We set $k_{max}=40$ and $k_{min}=10$.

Here, we compare to the recently proposed hybrid projection method with recycling, HyBR-recycle, from \cite{jiang2021hybrid}.
It is similar to ours in that the
regularization parameter is selected by working with the projected problem and subspace information is recycled.   
However, this algorithm corresponds to $\Psi = \bI$ and a two-norm constraint $q=2$. Thus the method differs from ours in how the subspaces are generated and what is recycled.   
We expect that while the method may be competitive in terms of memory required, it will not be qualitatively competitive because it does not effectively enforce edge-constraints the way s-RMM-GKS does.  

We describe here the numerical comparisons we make. In every algorithm, the regularization parameter is automatically selected using GCV on the projected problem. 

\begin{enumerate}
    \item Run HyBR (version with no recycling) \cite{jiang2021hybrid} on the first problem \eqref{eq:rtomo1} (all HyBR based codes assume $\Psi=\bI; q=2$). Then run HyBR on the full data problem. 
    We refer to these as HyBR 1st dataset and HyBR all data, respectively. 
    \item Run HyBR-recycle:  the first problem with HyBR is run to obtain an approximate solution and an initial subspace. This information is used in solving the second system with HyBR-recycle. The subspace obtained from solving the second system is then used to solve the third system with HyBR-recycle. We refer to the resulted solution as HyBR-recycle.
    \item Run MM-GKS on the first problem and on the full data problem. We refer to their results as MM-GKS 1st and MM-GKS all data, respectively. 
    \item We run RMM-GKS on the first problem (i.e. solve \eqref{eq:rtomo1} with $q=1$) and on all the data (i.e. solve \eqref{eq:rtomoproblemall}), that we refer to as RMM-GKS 1st and RMM-GKS all, respectively. 
   \item We further consider the case when the data are not available or can not be processed at once, hence, we employ our streaming version of RMM-GKS. First we run s-RMM-GKS for a fixed number of iterations for each problem (we choose $i_{max}$ in Algorithm ~\ref{Alg: RMMGKS} such that the total number of new basis vectors added overall in the subspace for each problem is 200). Further we decide to stop the iterations in the enlarge routine if the relative change of two consecutive computed approximations falls below some given tolerance $tol_1$. 
\end{enumerate}

For all MM-GKS and R-MMGKS variations we choose  \begin{equation}\label{eq: L}
    \Psi = \begin{bmatrix}
    \bI_{n_x} \otimes \bL_{x} \\ \bL_{y} \otimes \bI_{n_y}
\end{bmatrix},
\end{equation} where $d = x, y$ and $q=1$. 

\begin{figure}[ht!]
\centering
  \begin{tabular}{cccc}
   HyBR 1st &  MM-GKS 1st &  RMM-GKS 1st & s-RMM-GKS \\
  \includegraphics[height = 0.18\textwidth, width = .2\textwidth]{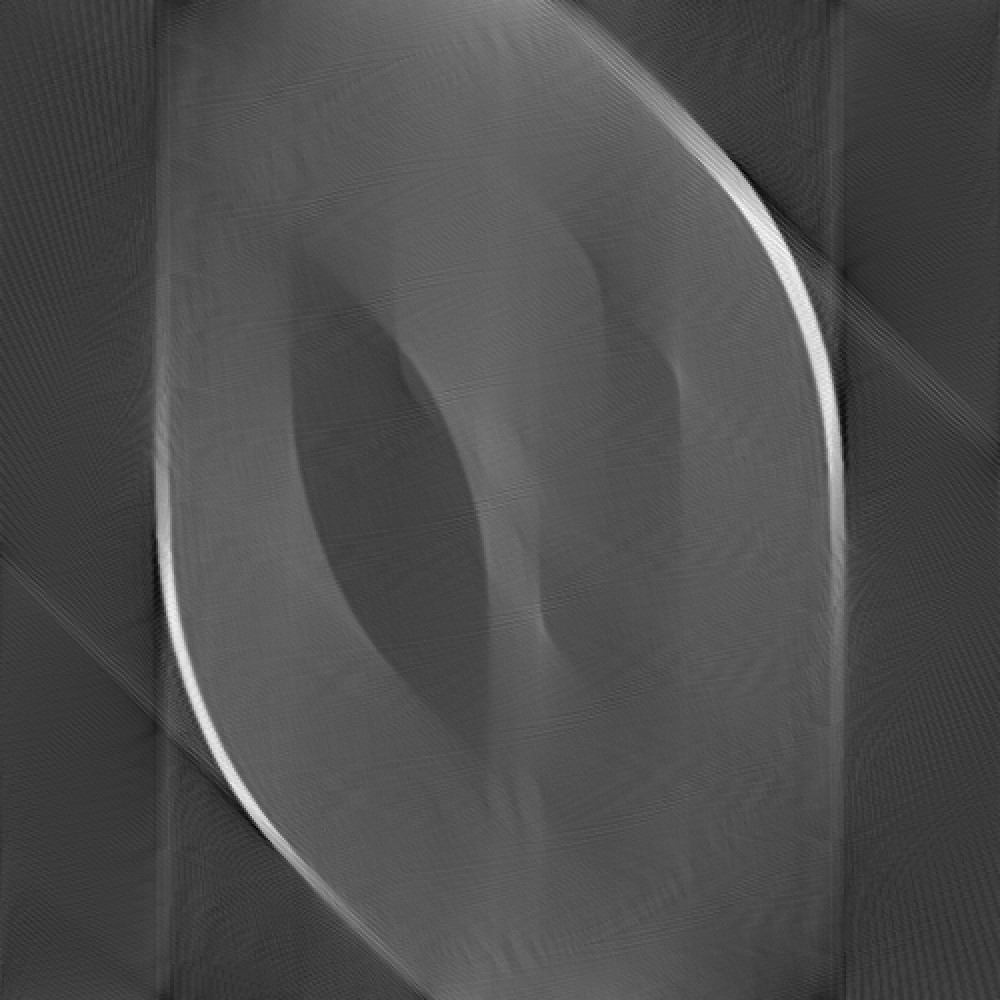} &
  \includegraphics[height = 0.18\textwidth, width = .2\textwidth]{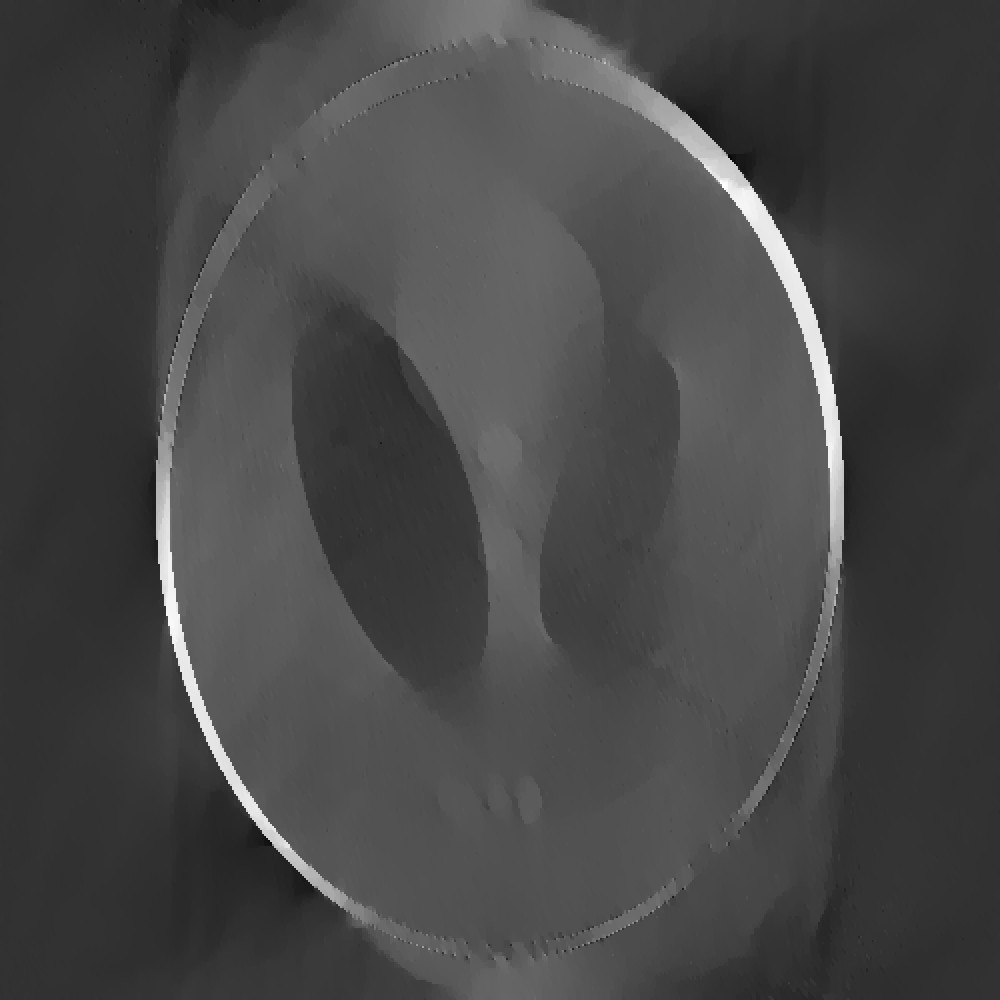} &
  \includegraphics[height = 0.18\textwidth, width = .2\textwidth]{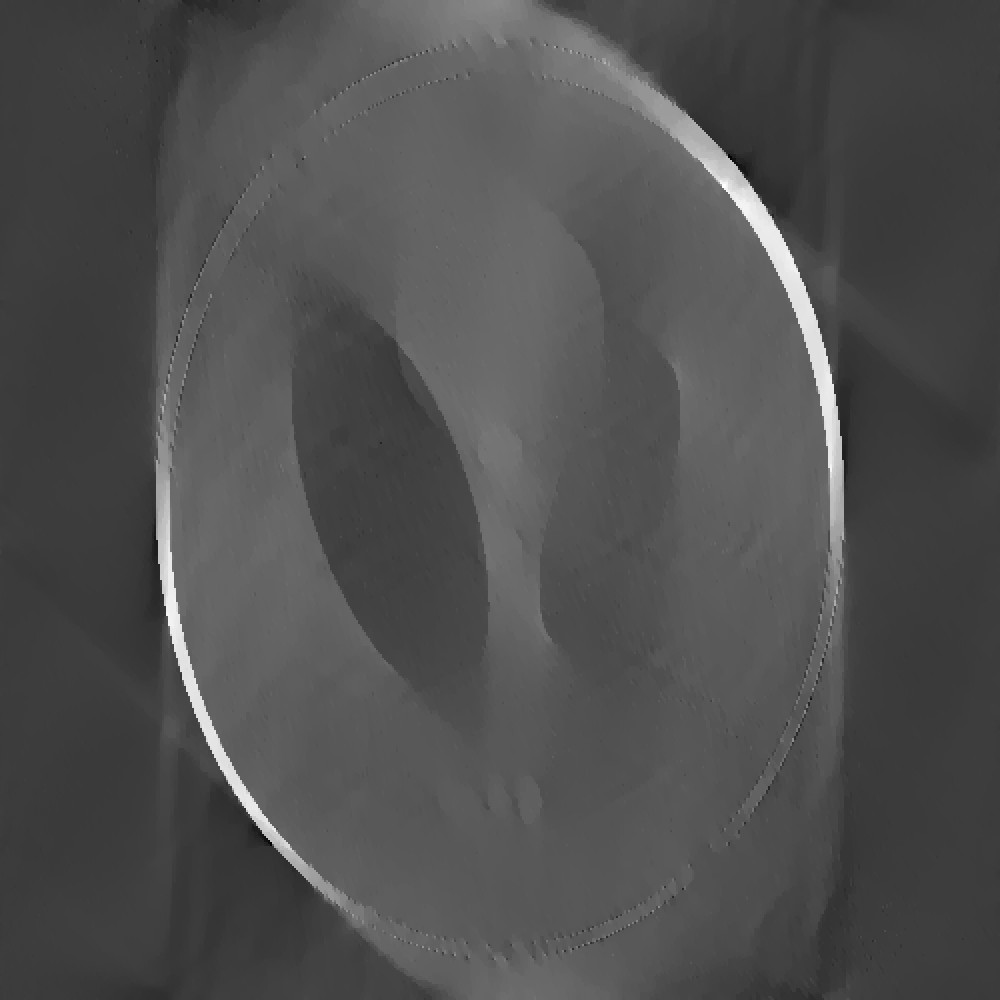}& 
    \includegraphics[height = 0.18\textwidth, width = .2\textwidth]{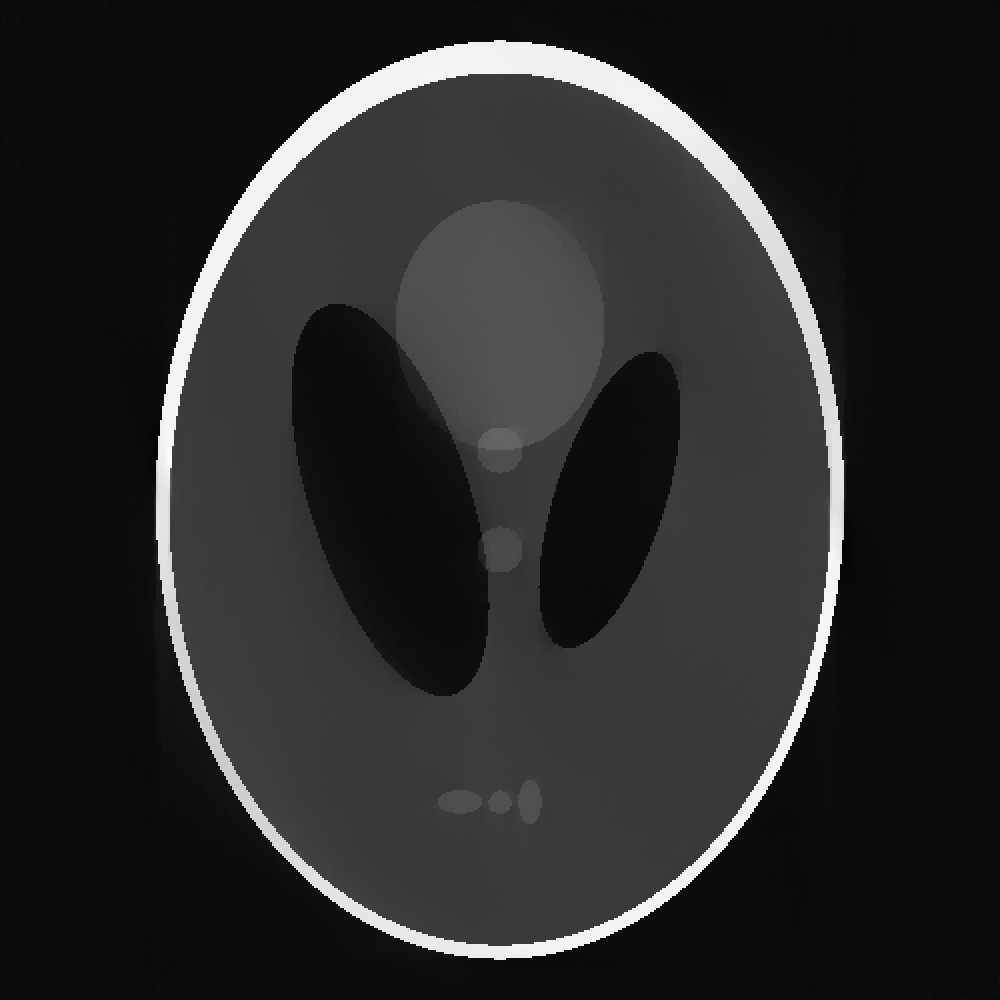} 
\end{tabular}
\begin{tabular}{cccc}
  \includegraphics[height = 0.18\textwidth,  width = .2\textwidth]{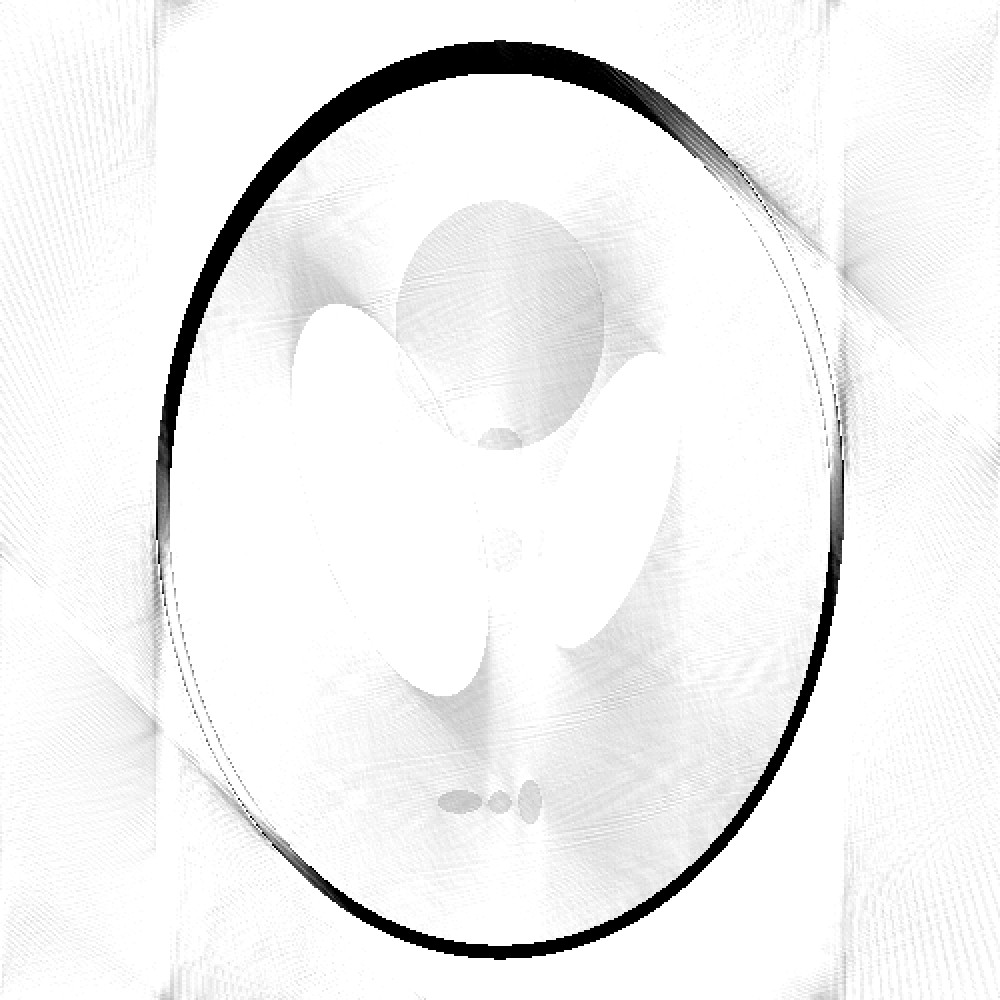}
    &
  \includegraphics[height = 0.18\textwidth,  width = .2\textwidth]{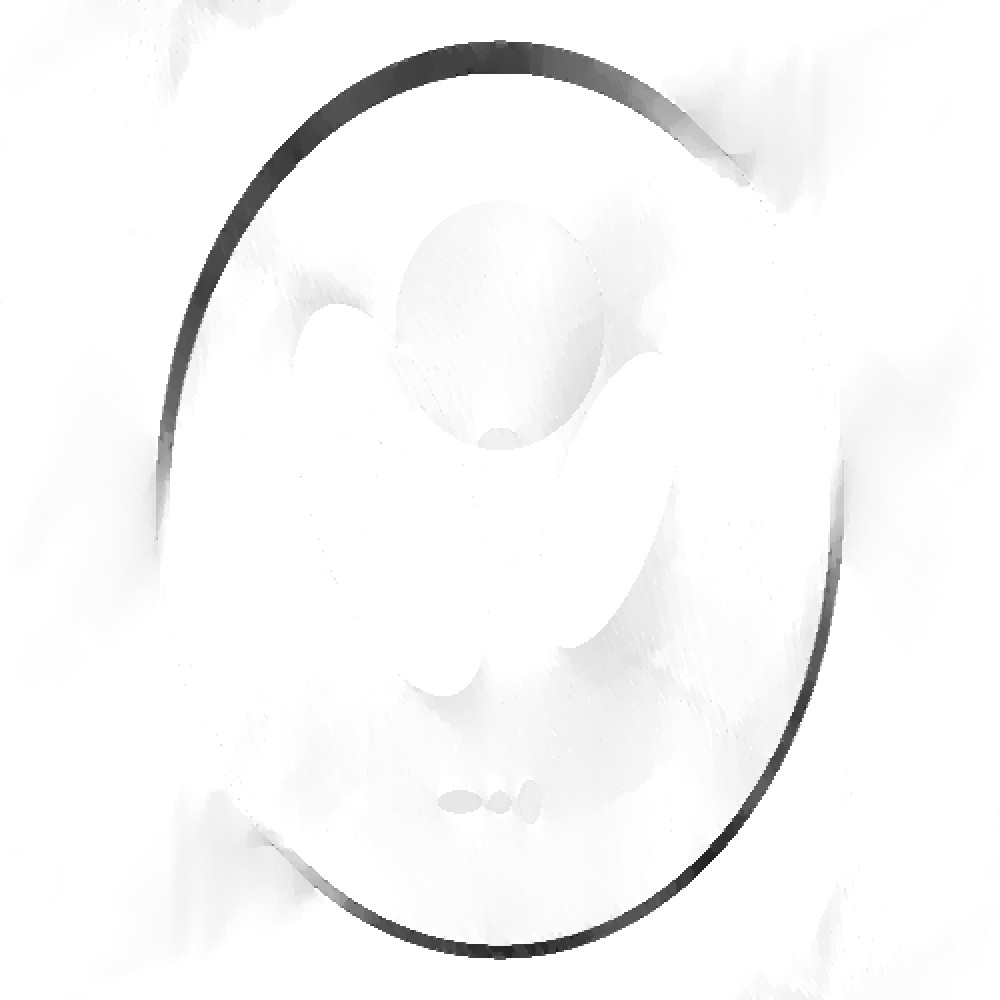}
    &
  \includegraphics[height = 0.18\textwidth, width = .2\textwidth]{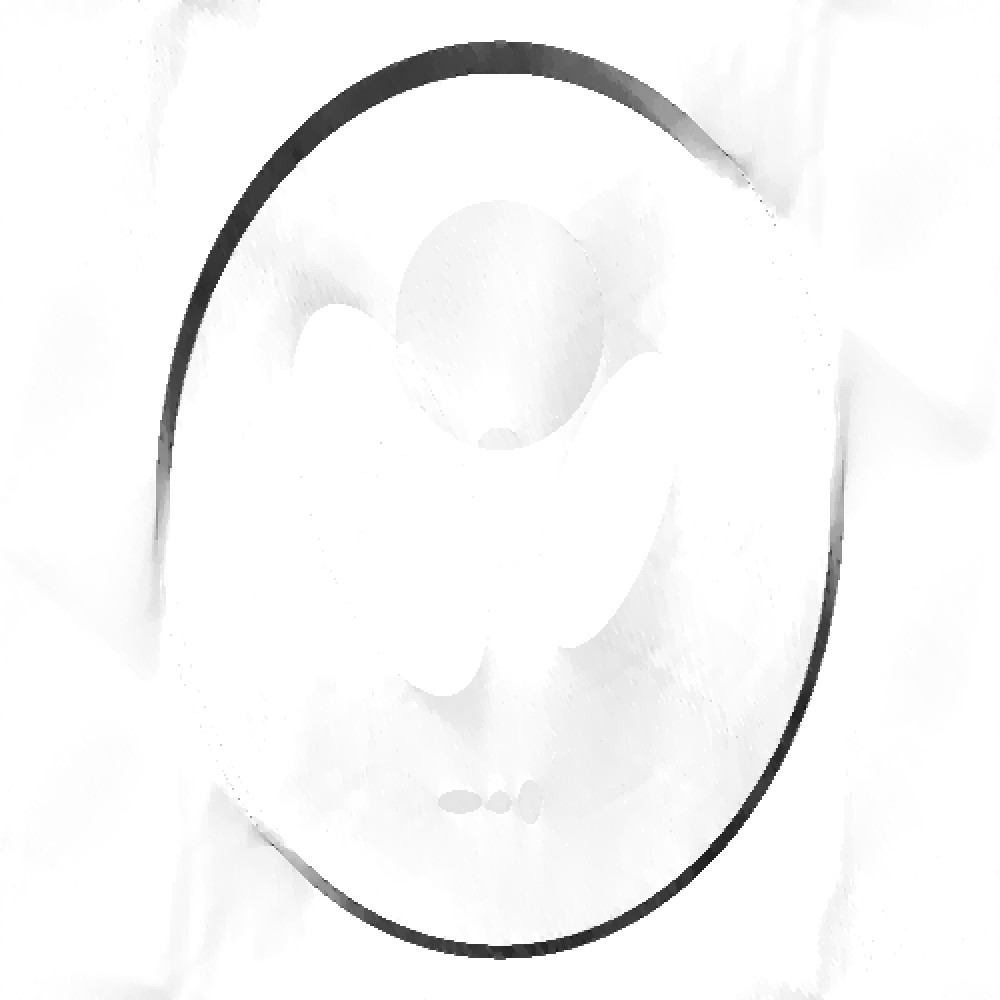}
     &
  \includegraphics[height = 0.18\textwidth, width = .2\textwidth]{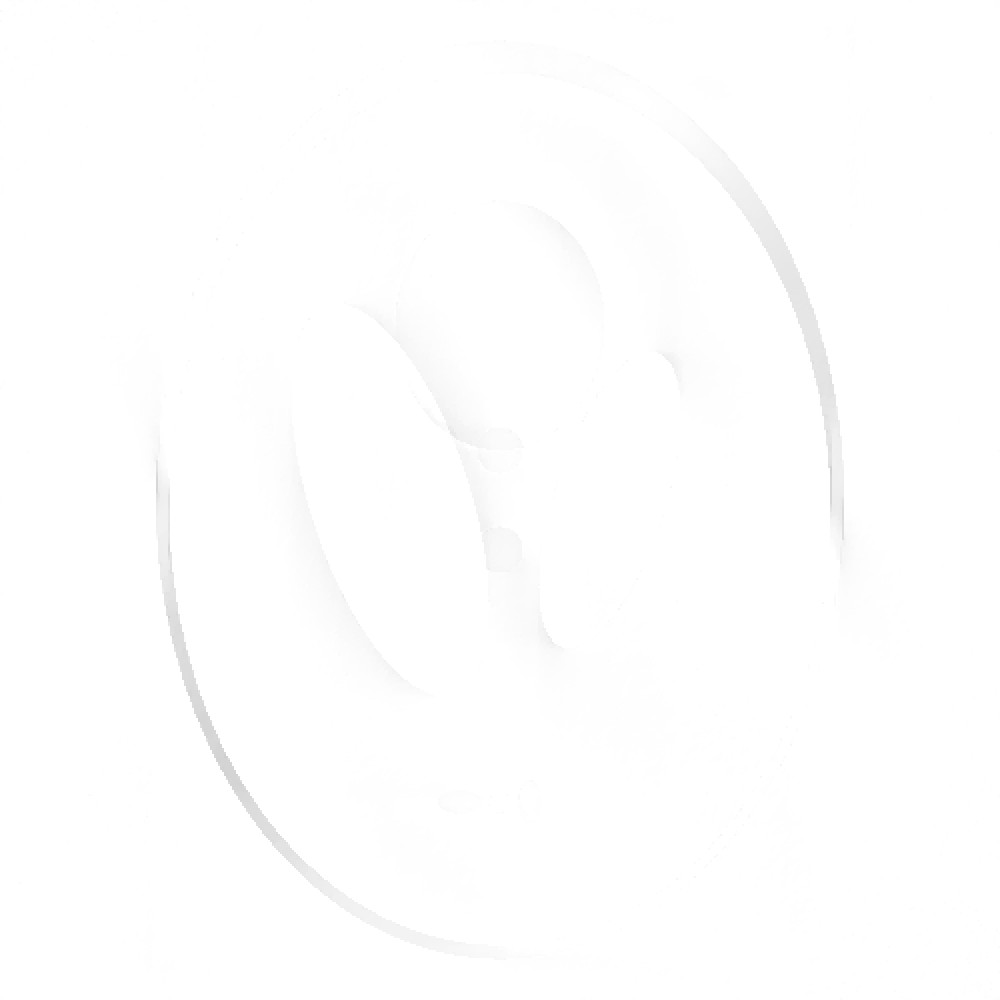}
\end{tabular}
\begin{tabular}{cccc}
HyBR all& HyBR-rec & MM-GKS all & RMM-GKS all \\
  \includegraphics[height = 0.18\textwidth,  width = .2\textwidth]{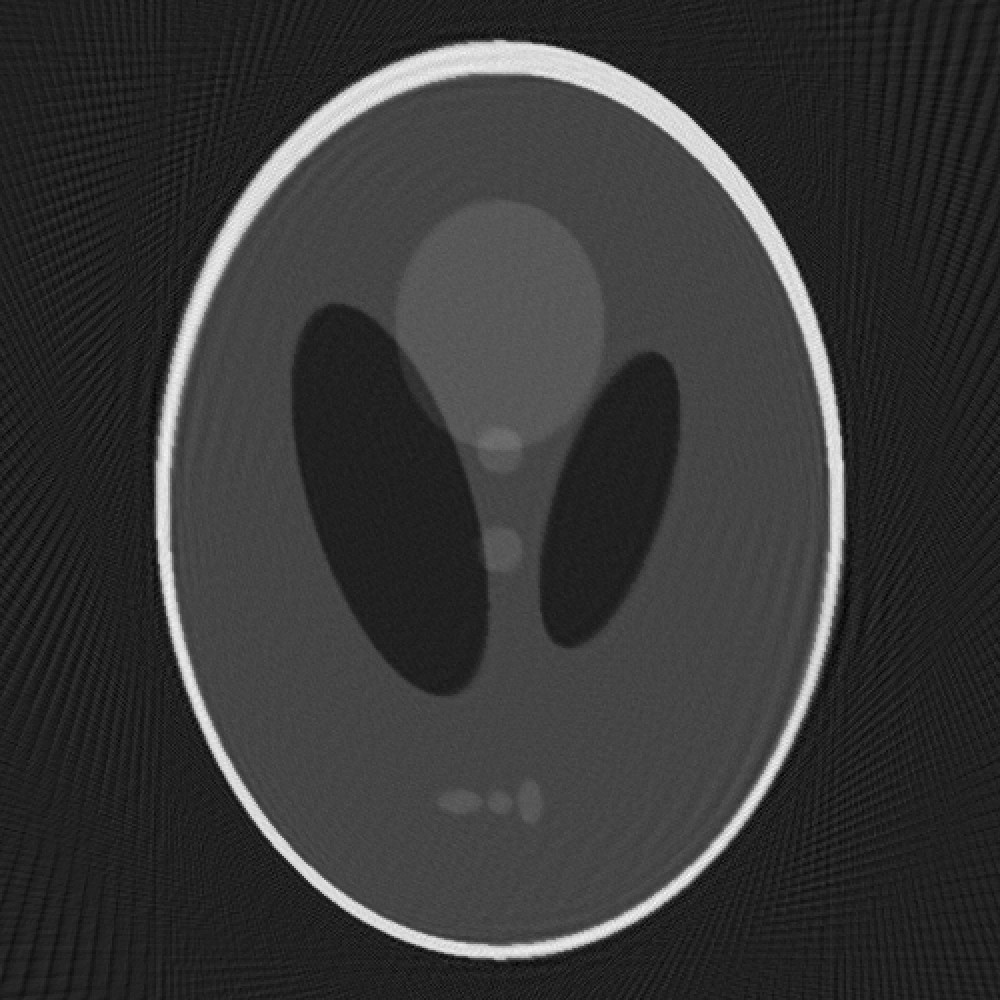}&
  \includegraphics[height = 0.18\textwidth,  width = .2\textwidth]{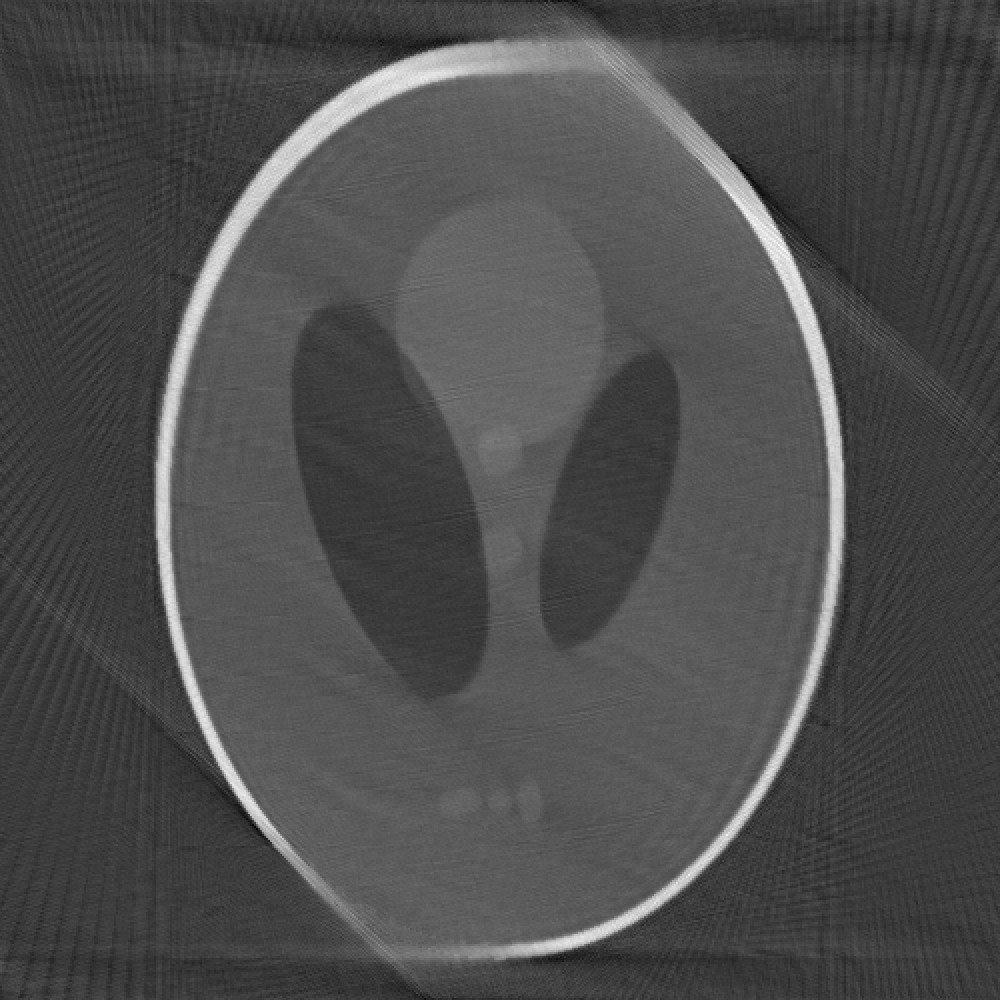}
    &
  \includegraphics[height = 0.18\textwidth,  width = .2\textwidth]{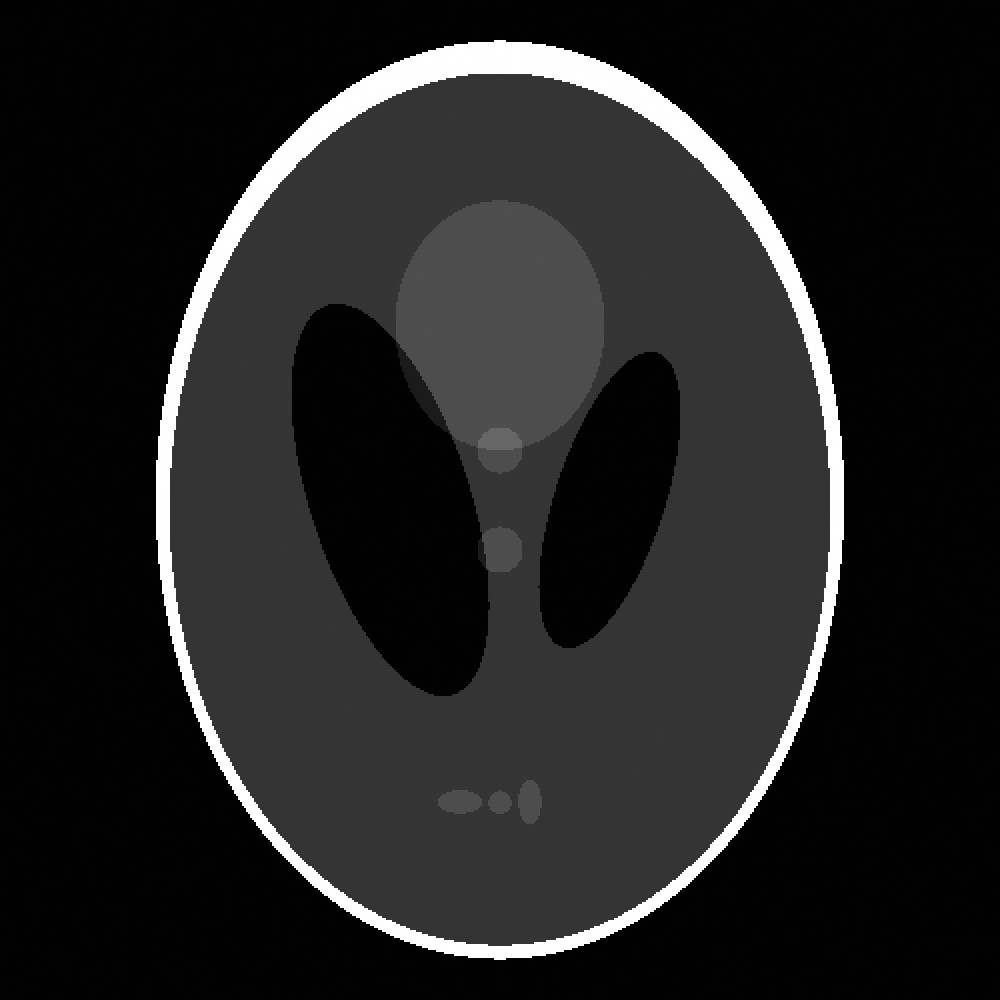}
    &
  \includegraphics[height = 0.18\textwidth,  width = .2\textwidth]{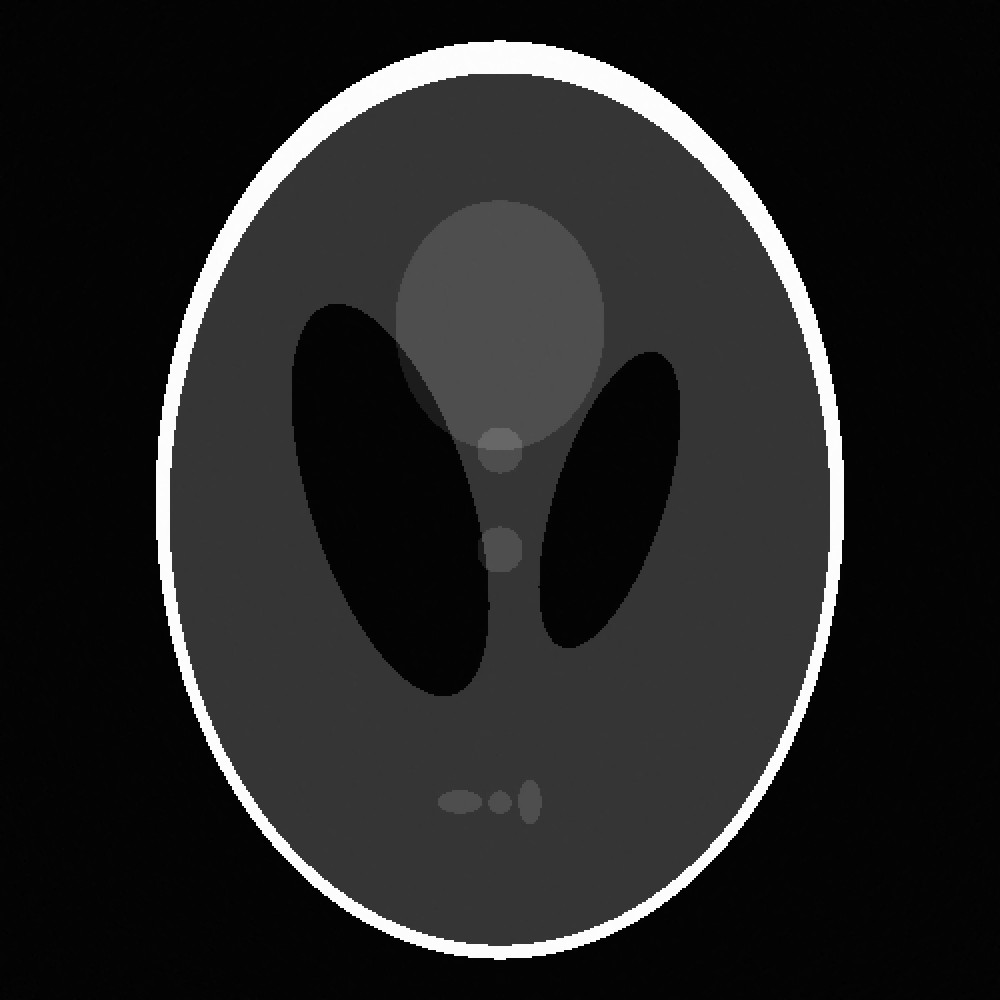}
\end{tabular}
\begin{tabular}{cccc}
  \includegraphics[height = 0.18\textwidth,  width = .2\textwidth]{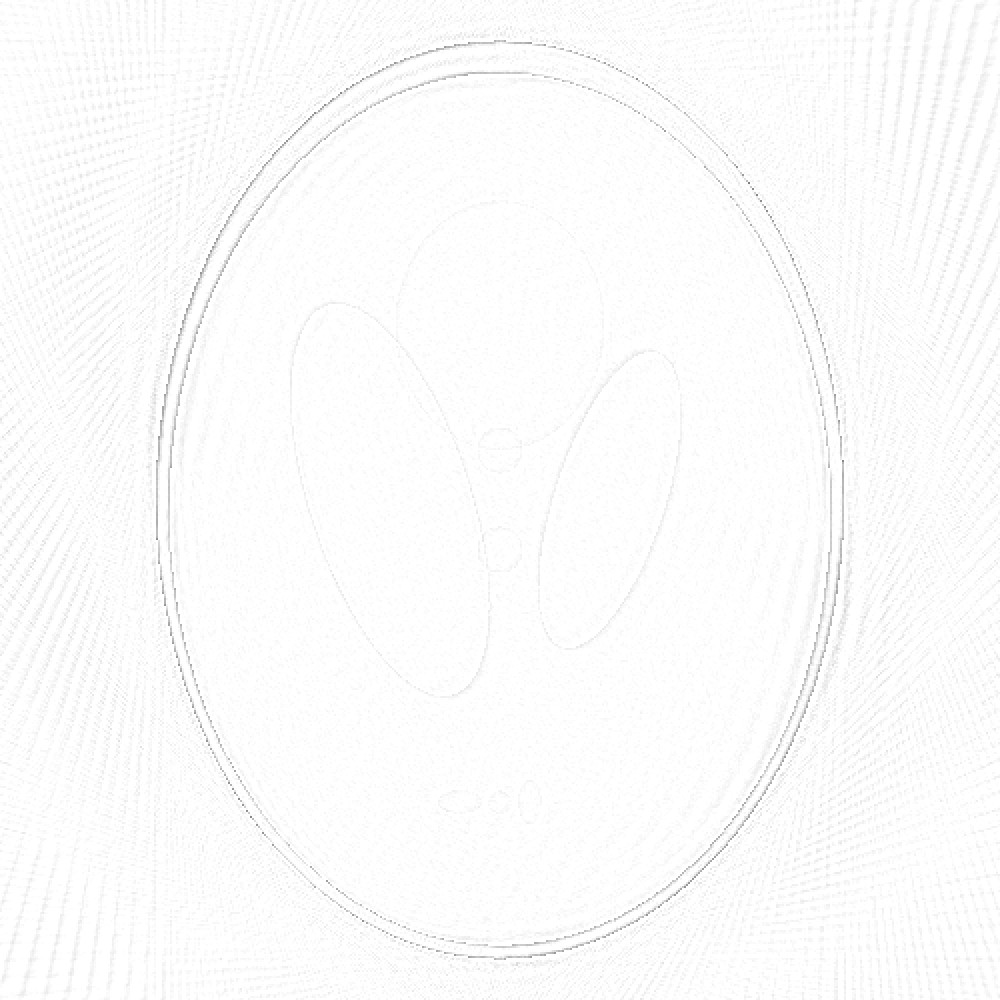}
    &
  \includegraphics[height = 0.18\textwidth,  width = .2\textwidth]{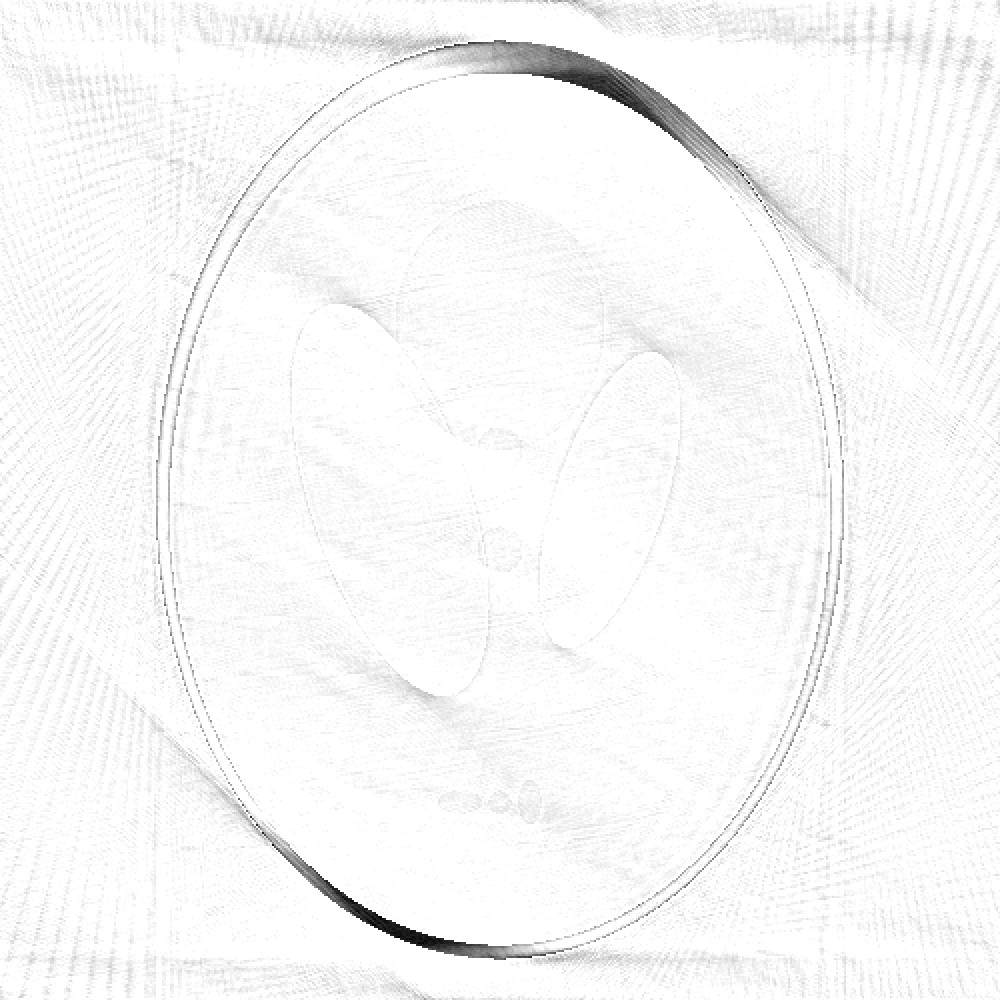}
    &
  \includegraphics[height = 0.18\textwidth,  width = .2\textwidth]{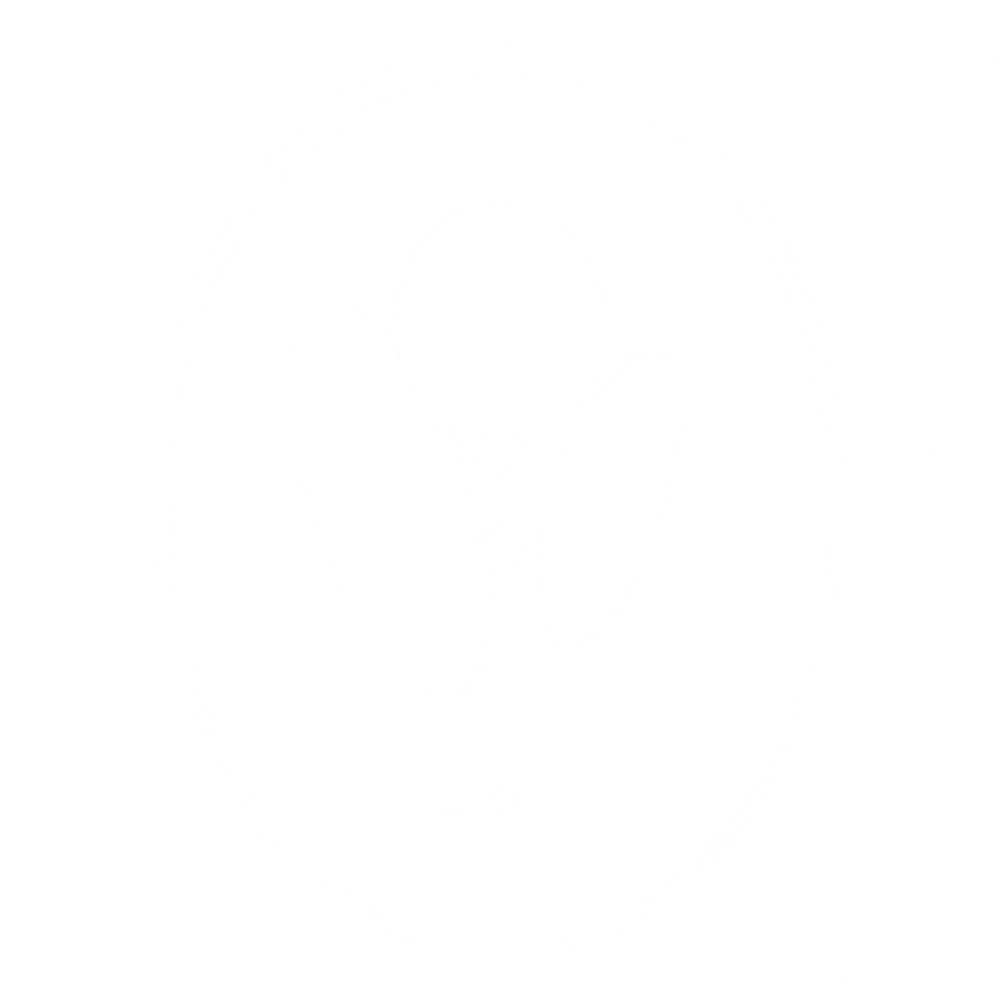}
     &
  \includegraphics[height = 0.18\textwidth,  width = .2\textwidth]{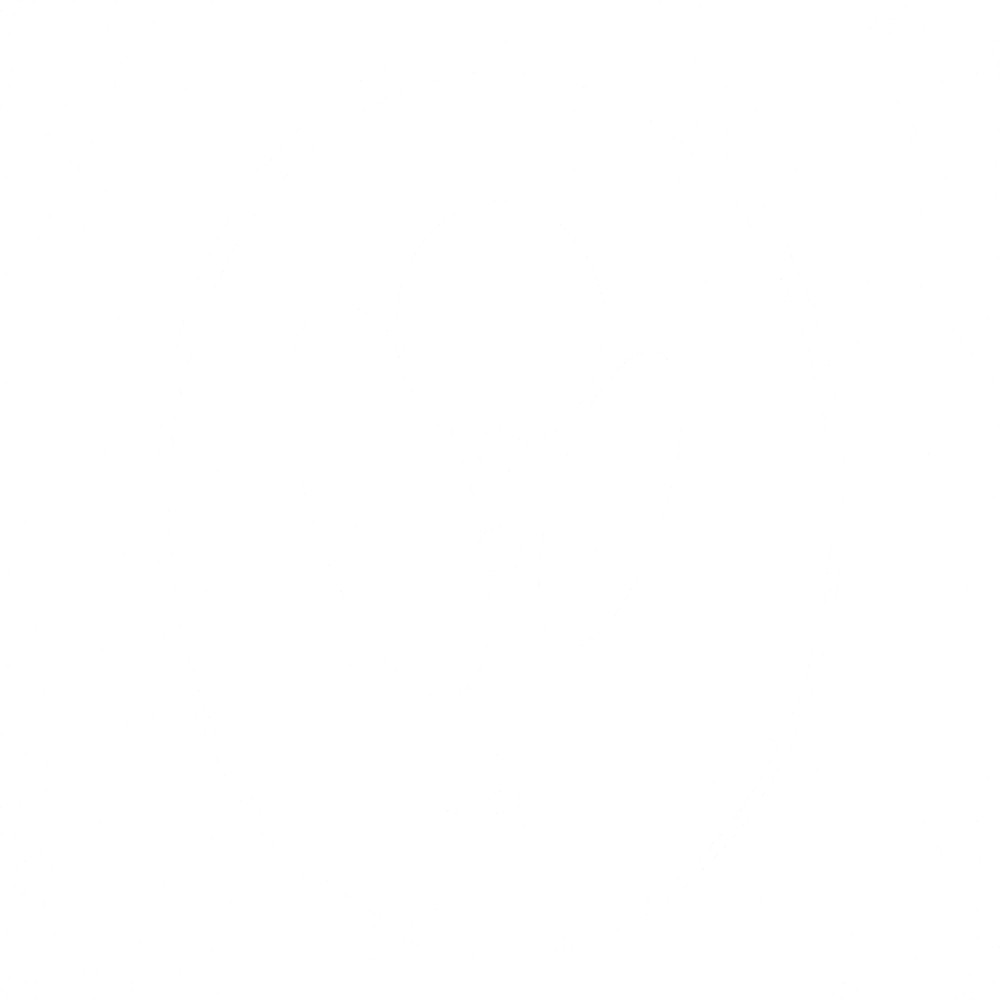}
\end{tabular}
\caption{Test 1: Streaming tomography example. Reconstructed images with the methods considered in the first and the third row. The second and the fourth row show the corresponding error images in the inverted colormap. The last two figures of the last row are shown as white due to the small values of the error images (see the corresponding RREs in Table~\ref{Table: tomo_stream_3prob_errors_noiselevles_differ}).}
  \label{fig: Streaming3prob_reconstructions}
\end{figure}

The reconstructed images along with the error images (i.e., the reconstructed image minus the true image), in inverted color map (we keep the unified scaling (-1, 0) when plotting all the error images) are shown in Figure \ref{fig: Streaming3prob_reconstructions}. 
We observe that both MM-GKS and RMM-GKS when used to solve all the problems together produce the best reconstruction results with RRE $0.39\%$ and $0.55\%$, respectively that are smaller than all the other methods considered. Nevertheless, MM-GKS requires to store all the solution subspace vectors, i.e., at iteration $200$ we have to store and operate with $200$ vectors with $95445$ entries, while for RMM-GKS we store at most $40$ solution vectors of size $95445$ and we achieve comparable results with MM-GKS. For scenarios when we can not solve all the problems at once (as MM-GKS all data and RMM-GKS all data do), we observe that solving only the first system produces relatively low quality reconstructions. Despite that some methods yield reconstructions with higher quality, we remark that solving only the first problem \eqref{eq:rtomo1} is a limited angle tomography problem, and to obtain a more accurate reconstruction we need more information collected from other angles, i.e., obtain information in the other systems \eqref{eq:rtomoi} and \eqref{eq:rtomont}. We do so by either solving all problems together or by recycling information from the previous systems. 

\begin{figure}
\centering
\includegraphics[width=0.99\textwidth, height = 0.3\textheight]{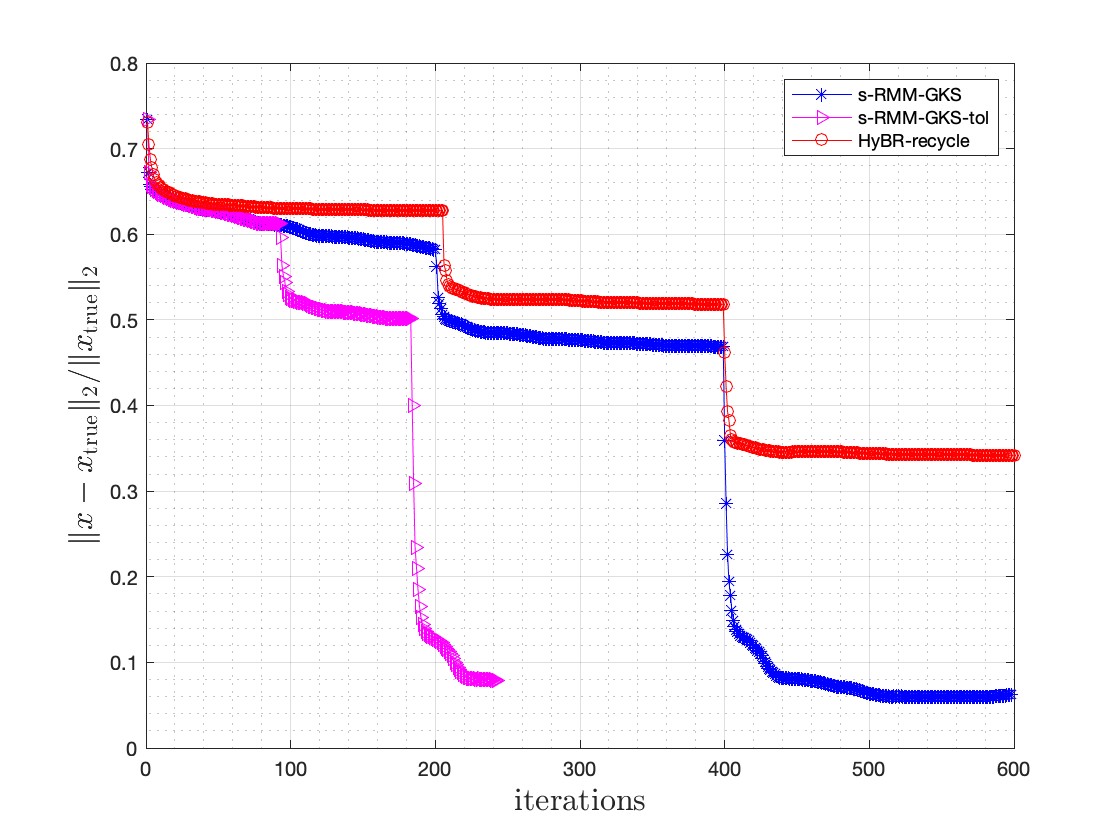}
\caption{Test 1: Streaming tomography example for 3 problems: RRE for s-RMM-GKS (blue-star curve) and HyBR-recycle (red-circle curve) for a fixed number or iterations. RRE for s-RMM-GKS with $\rm T_1 \leq 10^{-3}$.}
\label{Fig: sRMMGKSvsHyBRrecycle}
\end{figure}

Further, we show that the s-RMM-GKS method yields high quality reconstructions by only processing each system one by one and recycling from the previous system a compressed subspace and the computed approximate solution.
A comparison of the reconstruction RREs of our s-RMM-GKS with a previously developed recycling approach HyBR-recycle is shown in Figure~\ref{Fig: sRMMGKSvsHyBRrecycle}. The blue stared line shows the RRE for s-RMM-GKS and the red circled line shows the RRE for HyBR-recycle. For both methods we set $k_{min} = 10$, $k_{max} = 50$ and solved each problem by fixing the number of iterations to 200. The magenta triangled line shows the RRE for s-RMM-GKS when $\rm T_1 \leq 10^{-3}$ is used as a stopping criteria in addition to the maximum number of iterations. 

Lastly, for this test problem, we vary the noise level from $0.1$ to $1\%$ and we summarize the results of the obtained RREs for a group of methods we consider in Table \ref{Table: tomo_stream_3prob_errors_noiselevles_differ}. Such results emphasize that the methods we propose are robust with respect to the noise levels we consider.

\begin{table}[h]
\centering
\fontsize{8pt}{12pt}\selectfont
   \centering 
   \begin{tabular}{|c|c|c|c|c|c|c|c|} 
   \hline
   \hline
   \textbf{$\sigma$} & HyBR 1st & MM-GKS 1st & HyBR all &  MM-GKS all & HyBR-rec & RMM-GKS all & s-RMM-GKS\\
     \hline
   0.1\% &   0.6283  & 0.5340 & 0.1777 & 0.0039 &0.3523& 0.0055 & 0.0623\\
   \hline
   0.5\% & 0.6391  &  0.5459 & 0.3618 &
 0.0333 & 0.3607&  0.0391 & 0.1156 \\
   \hline
   1\% & 0.6448 &   0.5819 &  0.4203 & 0.0743 & 0.4154 & 0.0860 & 0.1584\\
   \hline
   \hline
   \end{tabular}
   \caption{Test 1: RRE for noise levels $0.1\%$, $0.5\%$, and $1\%$. Algorithm setting: Maximum memory limit is $40$, Image: Computerized tomography. Size: Image of size $500 \times 500$, All methods are run for $200$ iterations.}
   \label{Table: tomo_stream_3prob_errors_noiselevles_differ}
\end{table}

\subsubsection{Test 2} 
In the second case we assume that the data are streamed randomly as they are not available all at once. For this problem set up, we use IRTools \cite{gazzola2018ir} to generate four problems where 45  angels are selected equally spaced from the following angle intervals: $0^{\circ}-44^{\circ}$, $45^{\circ}-89^{\circ}$, $90^{\circ}-134$, and $135^{\circ}-179$. The four individual linear system matrices generated this way we will refer to as $\bB_i$, $i=1,\ldots,4$, and the corresponding sinogram data vectors in each of the four instances are $\bb_i$, $i=1,\ldots,4$. We then have $\bB_1, \bB_2, \bB_3, \bB_4 \in \R^{63630\times 1000^2}$ and $\bb_1, \bb_2, \bb_3, \bb_4 \in \R^{63630\times 1}$. To each right-hand side we add $0.1\%$ white Gaussian noise.

\begin{figure}[htbp]
\centering
  \begin{tabular}{ccccc}
  $\bx_{\rm true}$ & $\bb_1$: $0^{\circ}-44^{\circ}$ & $\bb_2$: $45^{\circ}-89^{\circ}$ &  $\bb_2$: $90^{\circ}-134^{\circ}$ & $\bb_2$: $135^{\circ}-179^{\circ}$\\
  \includegraphics[height = 0.15\textwidth, width = .15\textwidth]{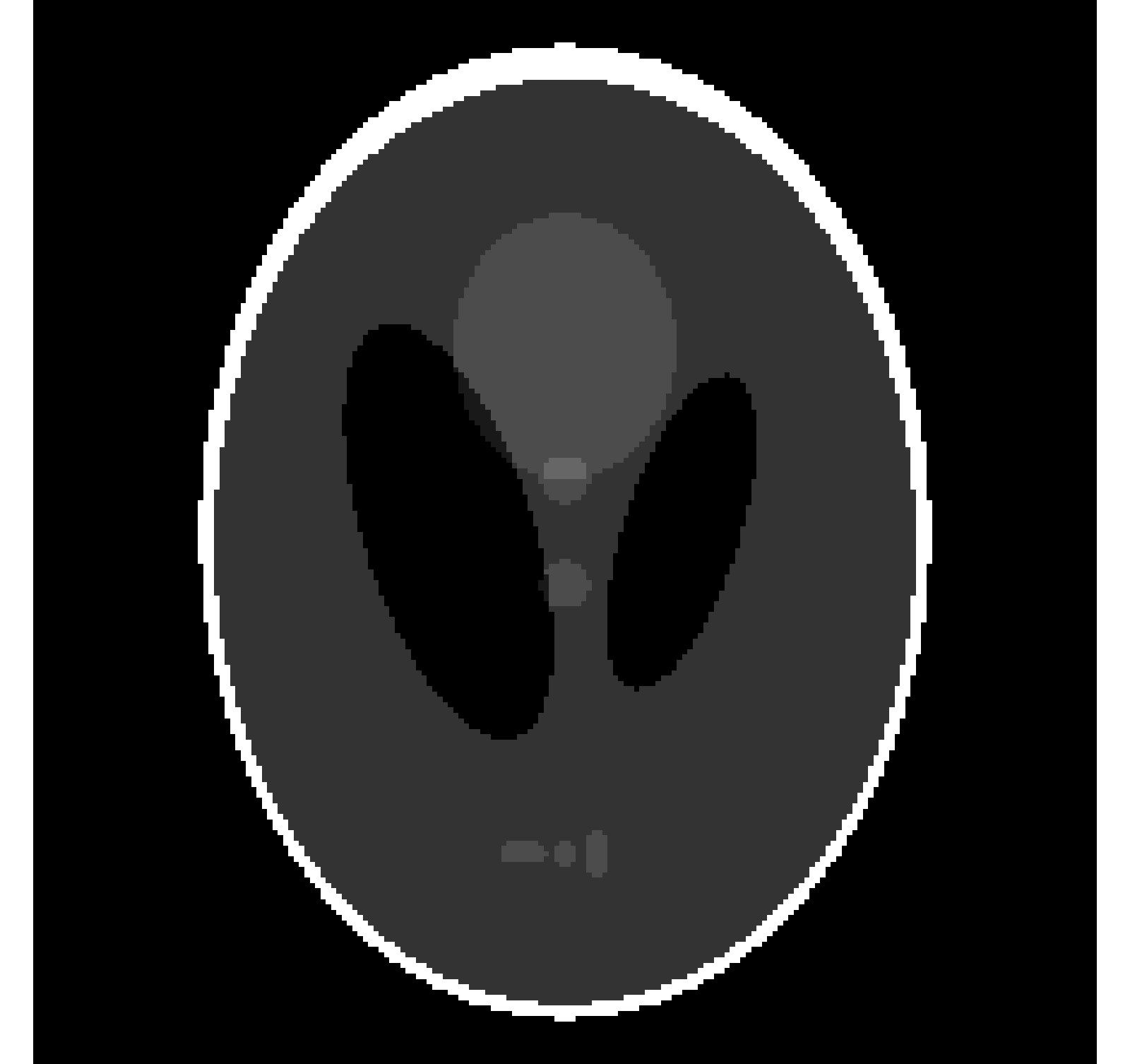} &
  \includegraphics[angle = 90,height = 0.15\textwidth, width = .15\textwidth]{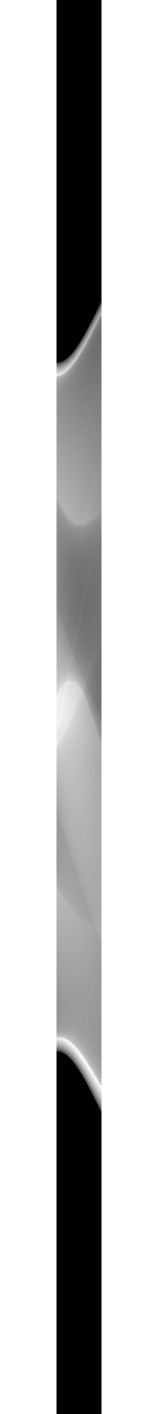} &
  \includegraphics[angle = 90,height = 0.15\textwidth, width = .15\textwidth]{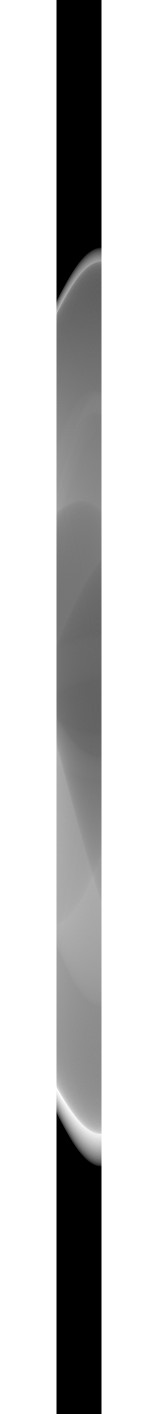} &
  \includegraphics[angle = 90,height = 0.15\textwidth, width = .15\textwidth]{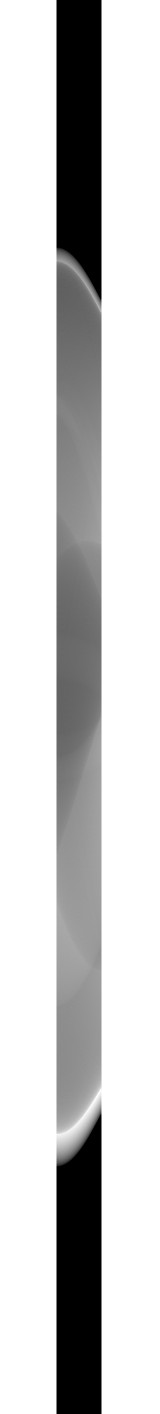} &
  \includegraphics[angle = 90,height = 0.15\textwidth, width = .15\textwidth]{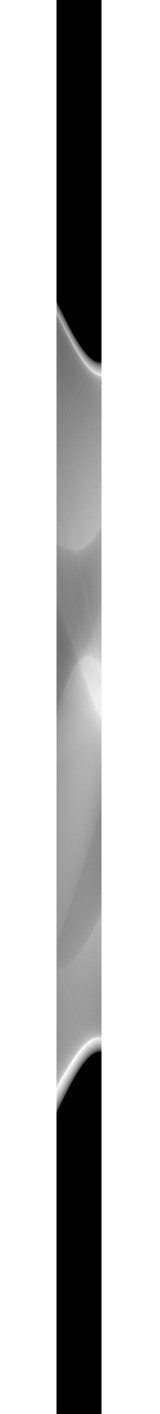}
  \\
  (a) & (b) & (c) & (d) & (e)
\end{tabular}
  \caption{Test 2: Streaming tomography example. The true image is provided in (a), along with three observed sinograms (b) $\bd_1$, (c) $\bd_2$, (d) $\bd_3$, and (e) $\bd_4$ corresponding to projections taken at intervals from $\bd_1$: $0^{\circ}-44^{\circ}$,  $\bd_2$: $45^{\circ}-89^{\circ}$, $\bd_3$: $90^{\circ}-134^{\circ}$, and $135^{\circ}-179^{\circ}$, respectively.}
  \label{fig: Streaming6prob_true}
\end{figure}

The true image is shown in Figure~\ref{fig: Streaming6prob_true}(a). We show the sinograms ($\bb_1, \bb_2, \bb_3, \bb_4$) in Figure~\ref{fig: Streaming6prob_true}(b), (c), (d), and (e).
The entire system that we want to solve is
$\bA \bx = \bd,$ where $\bA$ and $\bd$ are obtained by stacking the $\bB_i$
and $\bb_i$, respectively.

\begin{figure}
\centering
\includegraphics[width=0.95\textwidth, height = 0.3\textheight]{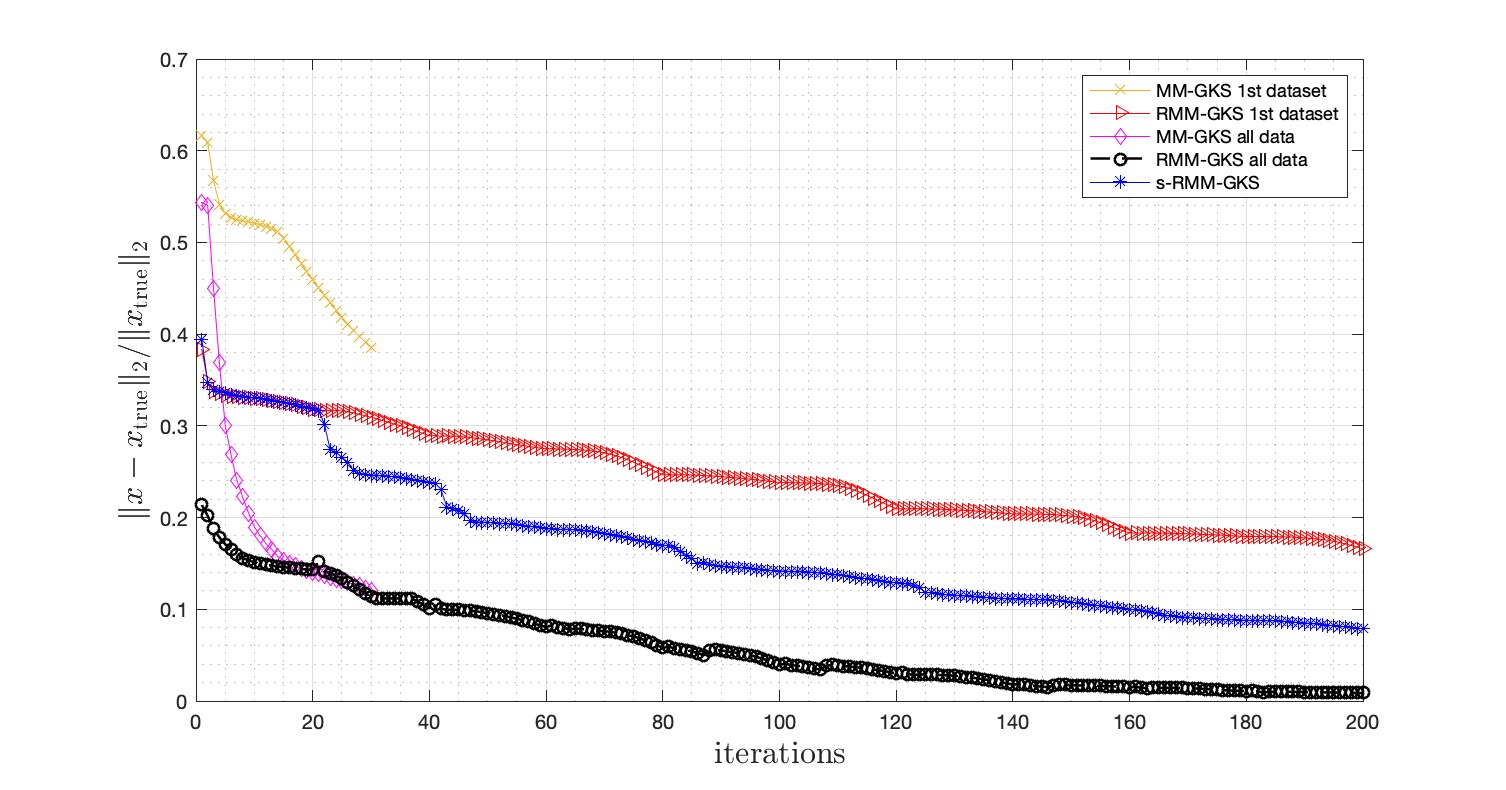}
\caption{Test 2: Streaming tomography example with 6 problems. Relative reconstruction norms.}
\label{Fig: Error_example2}
\end{figure}

\begin{figure}[ht!]
\centering
\begin{tabular}{ccccc}
 MM-GKS 1st &  RMM-GKS 1st & MM-GKS all & RMM-GKS all & s-RMM-GKS\\
\end{tabular}
  \begin{tabular}{ccccc}
  \includegraphics[height = 0.13\textwidth, width = .15\textwidth]{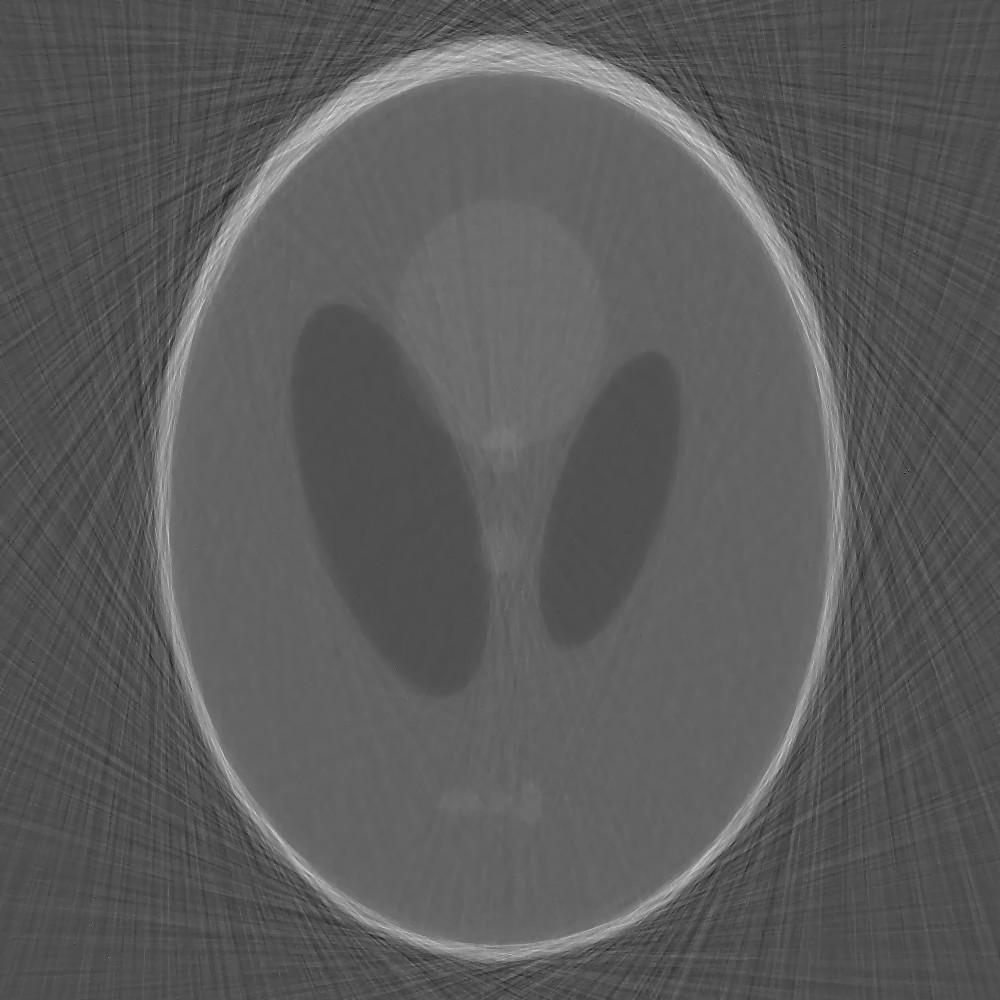} &
  \includegraphics[height = 0.13\textwidth, width = .15\textwidth]{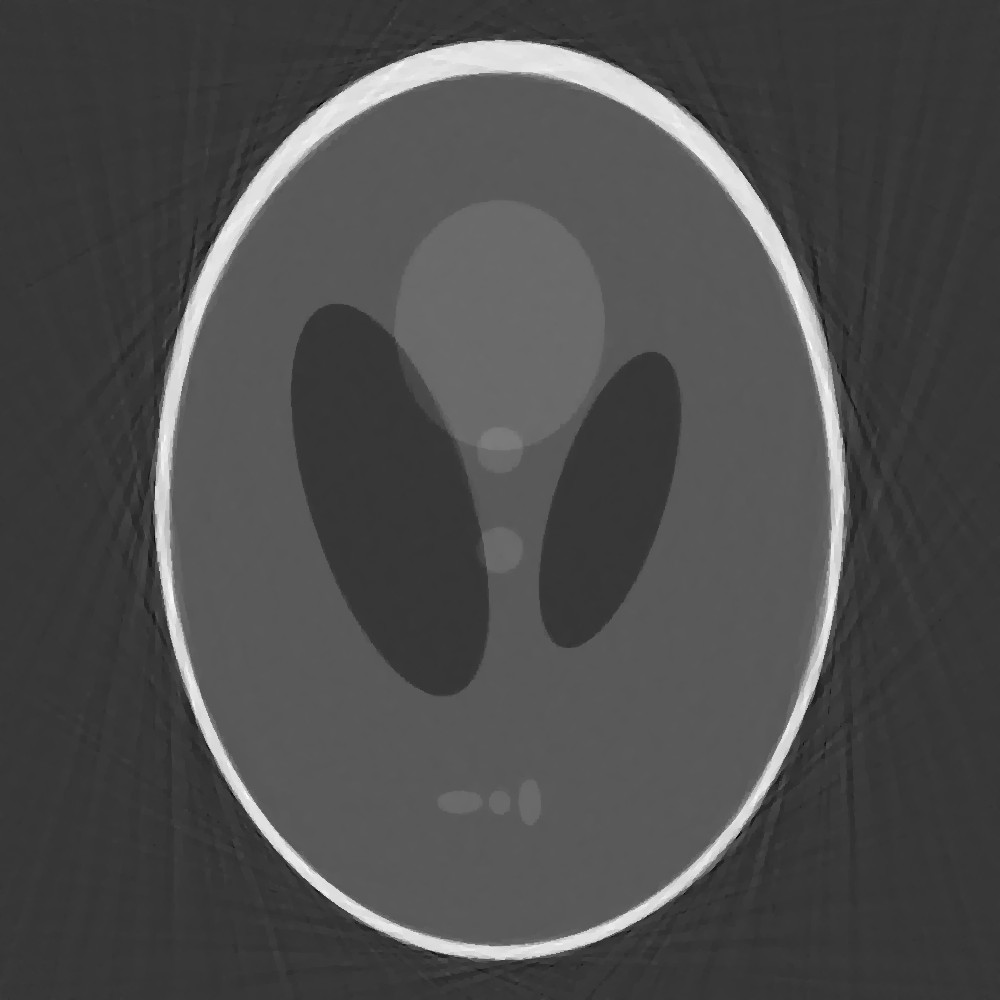} &
  \includegraphics[height = 0.13\textwidth, width = .15\textwidth]{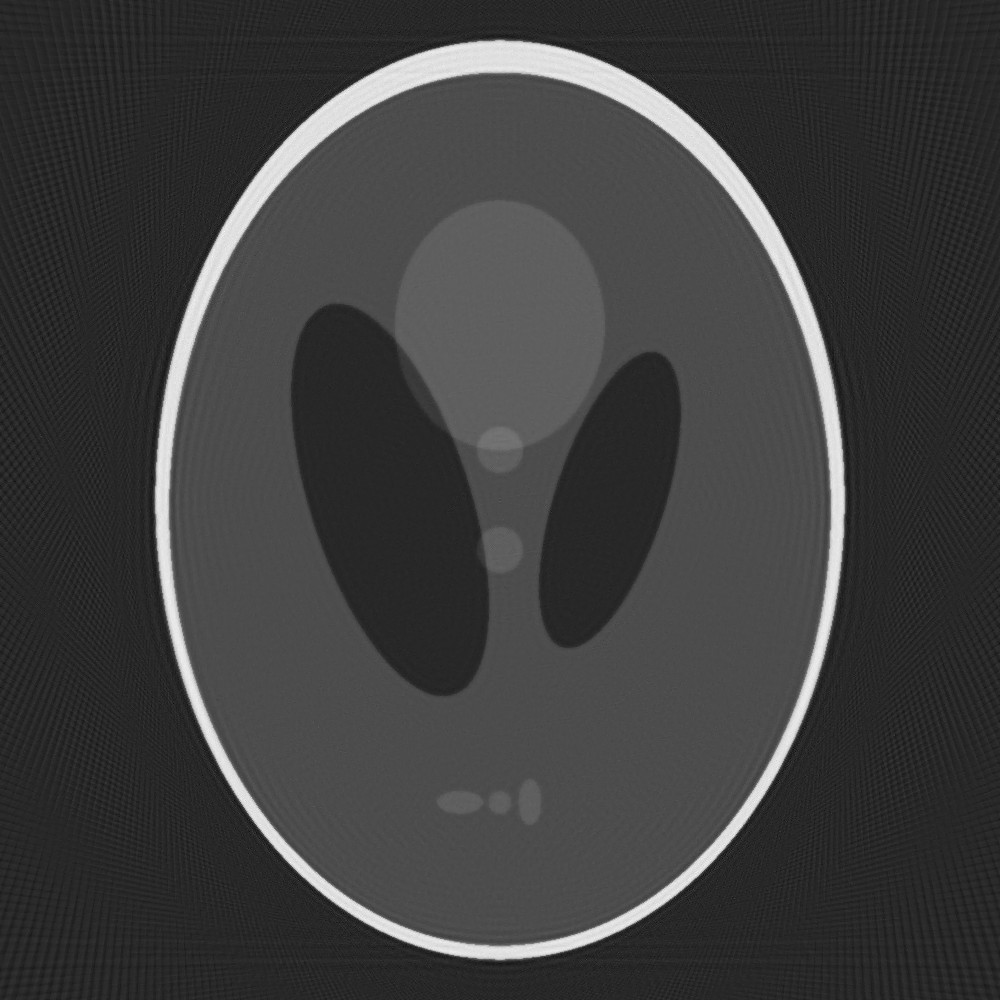}& 
    \includegraphics[height = 0.13\textwidth, width = .15\textwidth]{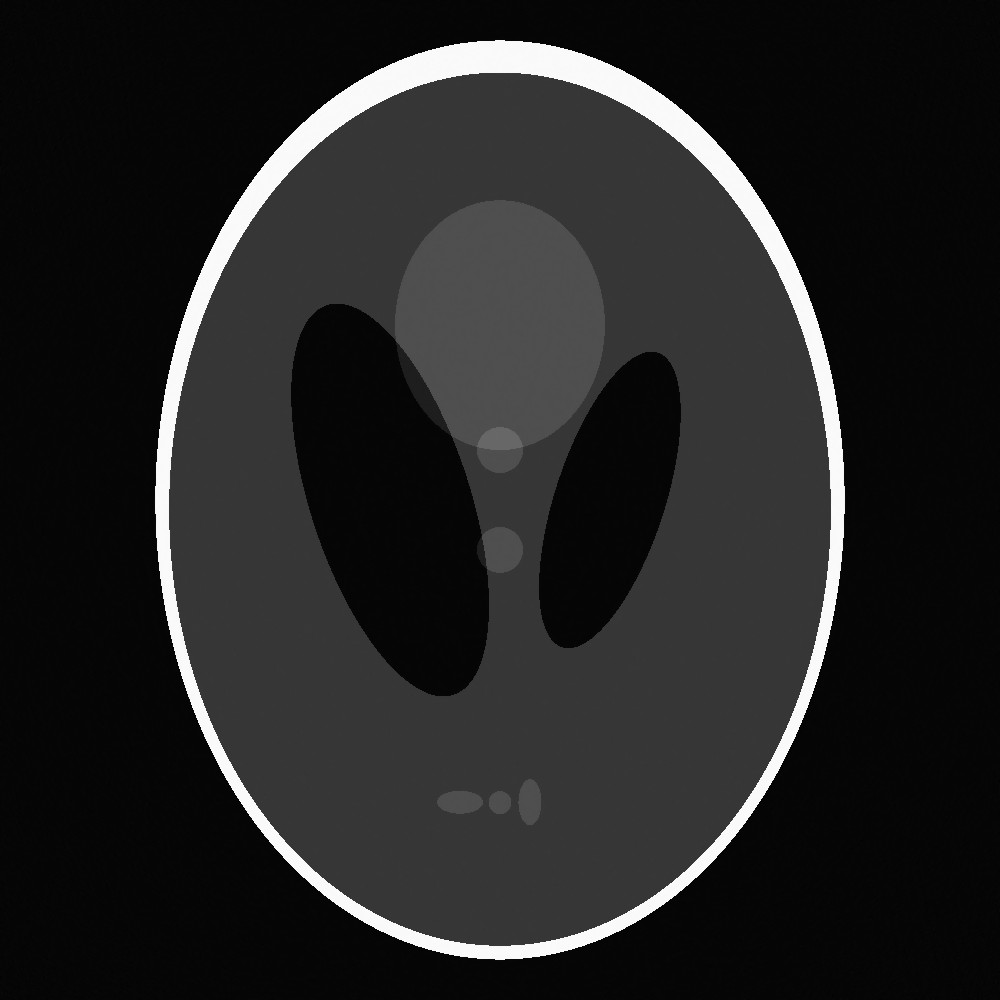} 
    & 
    \includegraphics[height = 0.13\textwidth, width = .15\textwidth]{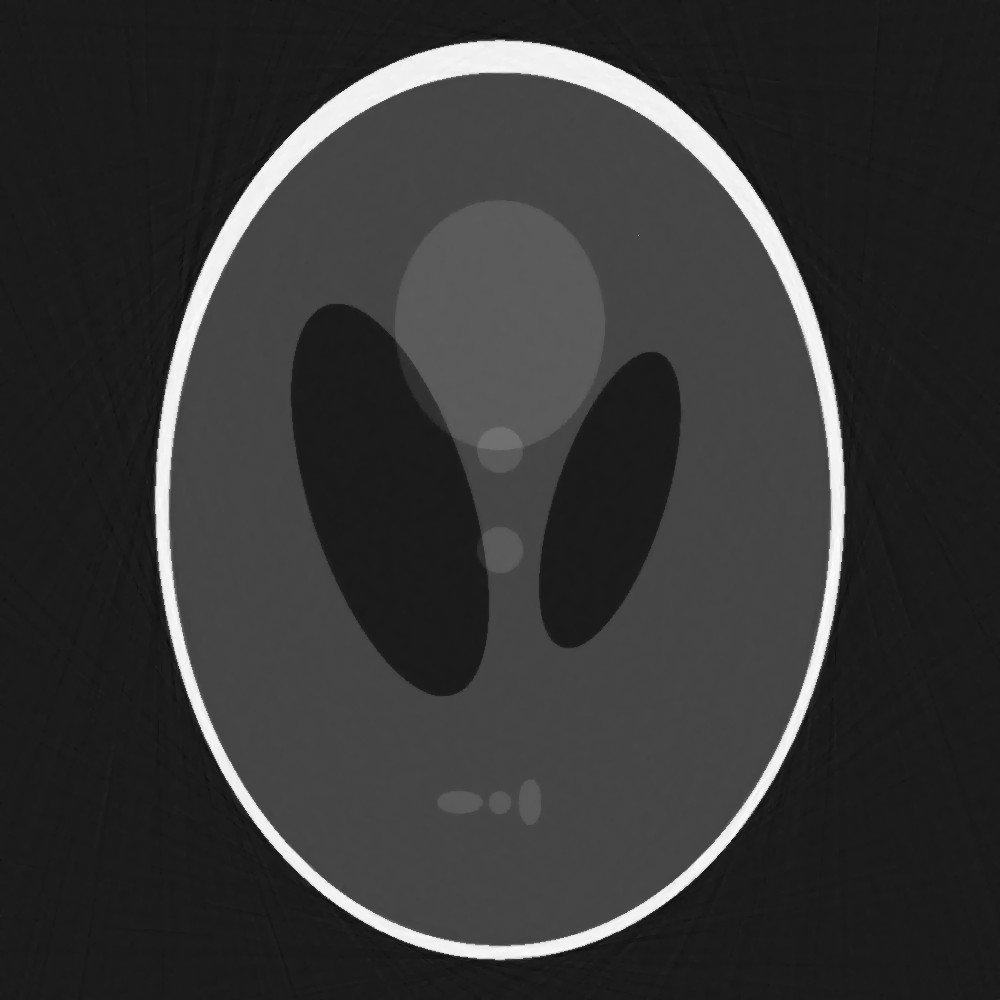} 
\end{tabular}
  \begin{tabular}{ccccc}
  \includegraphics[height = 0.13\textwidth, width = .15\textwidth]{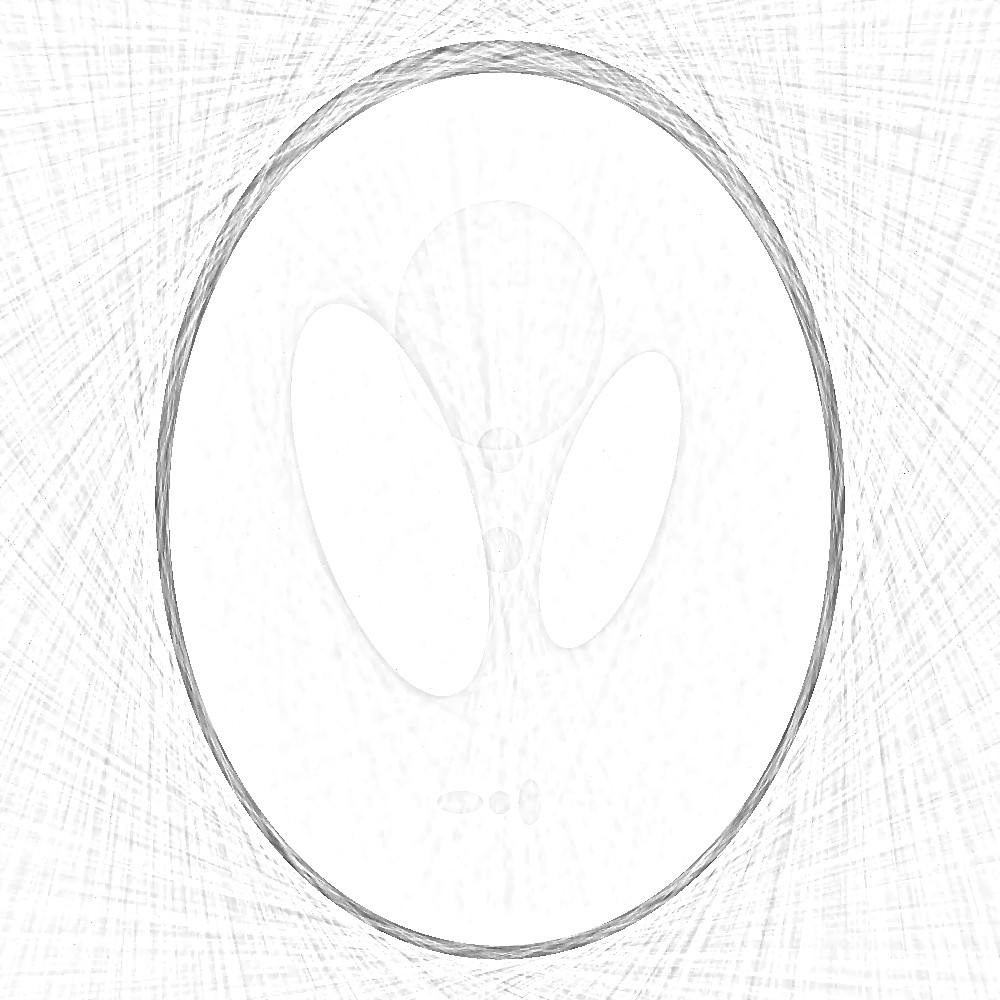} &
  \includegraphics[height = 0.13\textwidth, width = .15\textwidth]{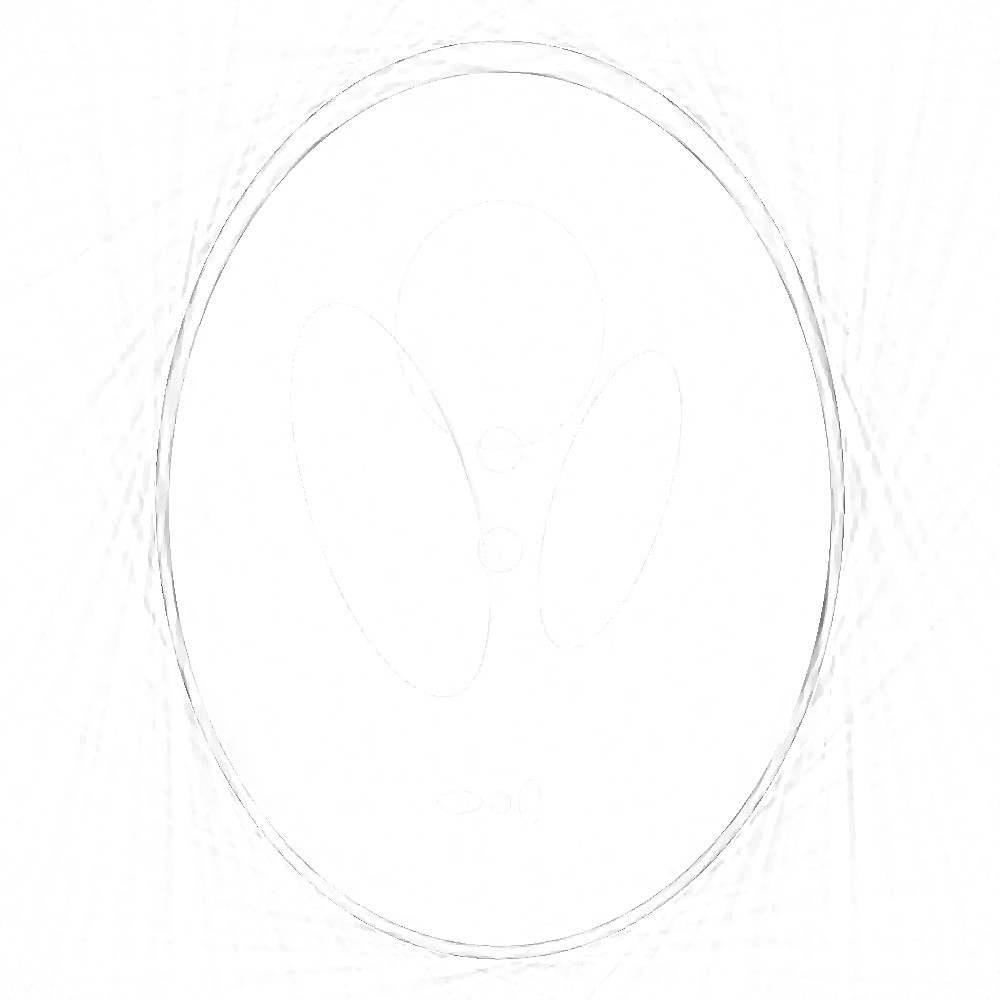} &
  \includegraphics[height = 0.13\textwidth, width = .15\textwidth]{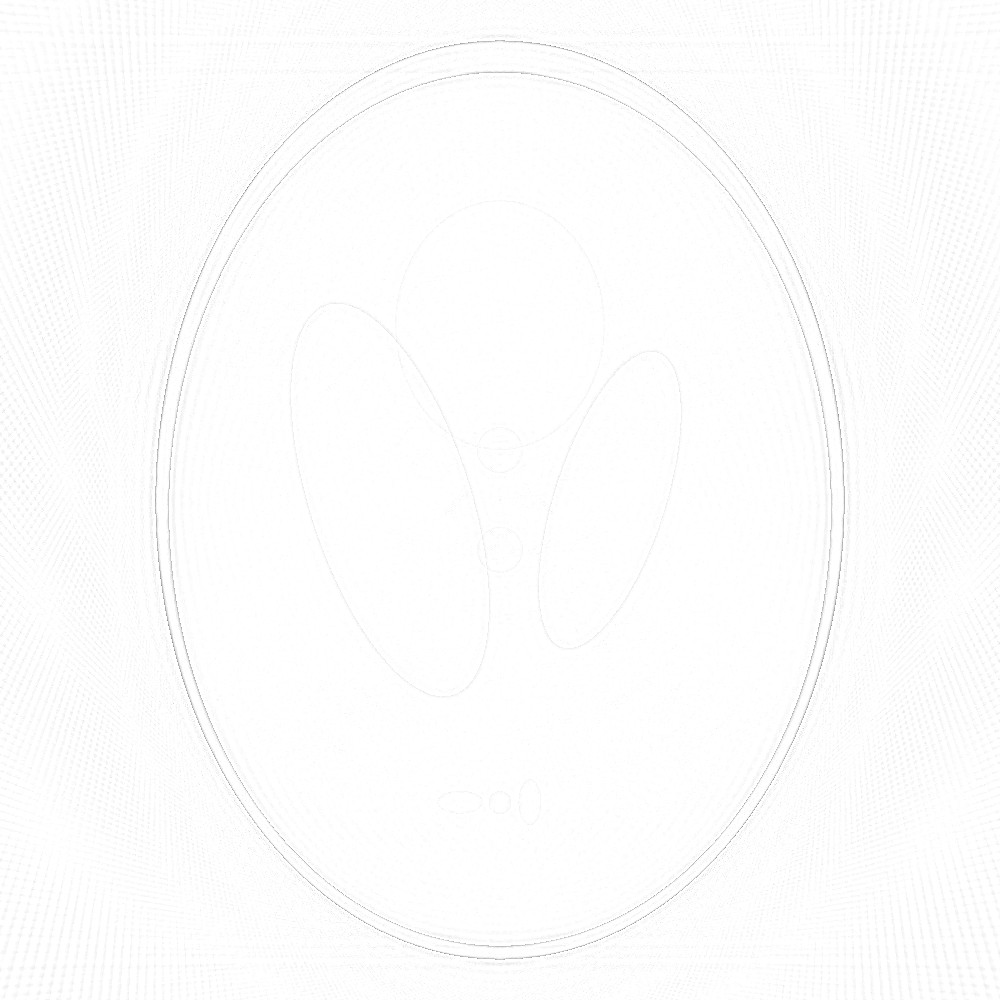}& 
    \includegraphics[height = 0.13\textwidth, width = .15\textwidth]{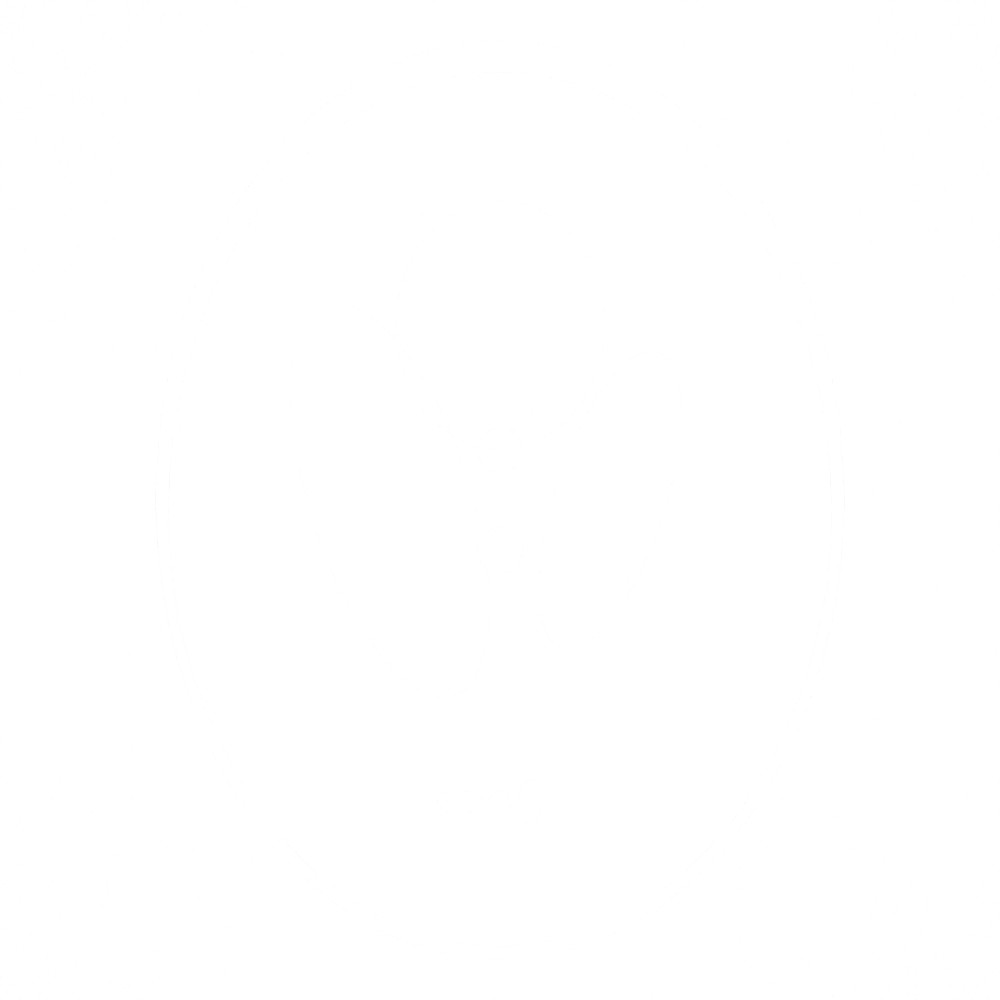} 
    & 
    \includegraphics[height = 0.13\textwidth, width = .15\textwidth]{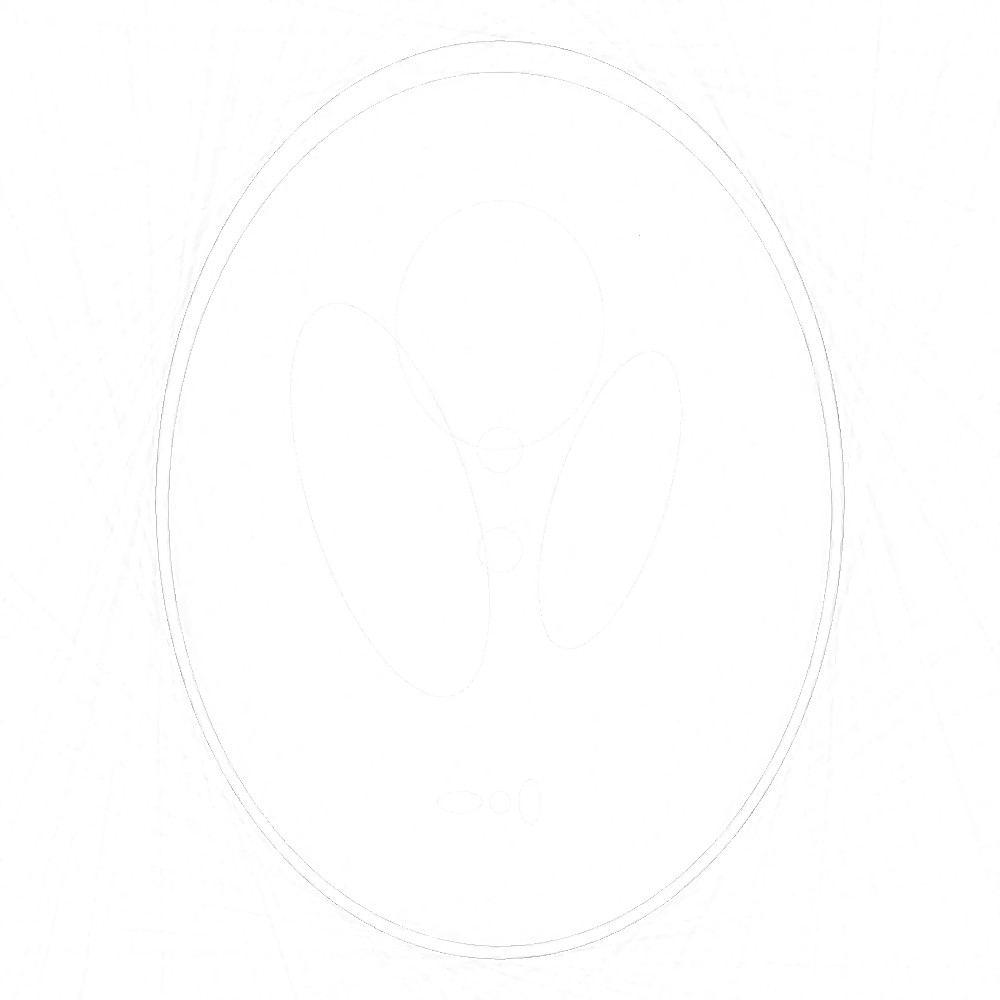} 
\end{tabular}
  \caption{Test 2: Streaming tomography example. Reconstructions for MM-GKS and RMM-GKS on the limited data from the first dataset and on all the data, from left to right, respectively. The last image shows the  reconstruction using s-RMM-GKS.}
  \label{fig: Streaming6prob_reconstructions}
\end{figure}

Assume that during the data acquisition process the data is streamed in $6$ blocks, i.e., for each block we collect data through $30$ random angles (not necessary equally spaced, but each block has the same size). From the large system $\bA \bx = \bd$, we obtain $6$ sub-systems to solve whose matrices we denote by $\bA_1, \bA_2, \bA_3, \bA_4, \bA_5, \bA_6 \in \R^{42420 \times 1000^2}$ and with corresponding data subvectors $\bd_1, \bd_2, \bd_3, \bd_4, \bd_5, \bd_6 \in \R^{42420 \times 1}$. To test and compare the performance of RMM-GKS on streaming data, we consider the following scenarios.
\begin{enumerate}
    \item Run MM-GKS on the first subproblem (MM-GKS 1st). Also run MM-GKS on all the data, i.e., solve \eqref{eq:rtomoproblemall} with $n_t=6$ (MM-GKS all) for 30 iterations, where 30 corresponds to the maximum number of vectors allowed in the subspace.
    \item Run RMM-GKS on the first subproblem (RMM-GKS 1st). Also run RMM-GKS on the full data problem (RMM-GKS all) for 200 iterations. 
    \item Run s-RMM-GKS method first for a total of 360 iterations (here we run a fixed number of iterations for each subproblem). Then, we consider the case of stopping the iterations based on some stopping criteria. In particular we stop if either the maximum number of iterations (100 for each subproblem) or a desired tolerance is achieved (we set $tol_1 = 10^{-3}$). 
\end{enumerate}
For this example we focus our comparison on MM-GKS under a limited memory assumption. Specifically, assume that the memory capacity is 30 basis vectors. To mimic this, we run MM-GKS type methods and we report their result for only 30 iterations. RMM-GKS and s-RMM-GKS keep the memory bounded by 30 basis vectors in the solution subspace, too, but the enlarging and compressing technique allows us to run for a larger number of iterations. The RREs for the methods we consider are shown in Figure \ref{Fig: Error_example2}. We observe that RMM-GKS applied on all the data outperforms all the methods considered. Moreover, we are able to achieve such results with only the memory capacity vectors stored. The reconstructed images with all the methods along with the error images in the inverted colormap are displayed in Figure~\ref{fig: Streaming6prob_reconstructions}\footnote{The last two images of the second row appear almost white due to the small error they represent. They correspond to the dotted-black and starred blue lines on Figure~\ref{Fig: Error_example2}}. Further, we show the reconstructed images at each stage of streaming along with the error images in the inverted colormap. An obvious increase in the reconstructed quality is observed until convergence. 

\begin{figure}[ht!]
\centering

  \begin{tabular}{cccccc}
  \includegraphics[height = 0.13\textwidth, width = .13\textwidth]{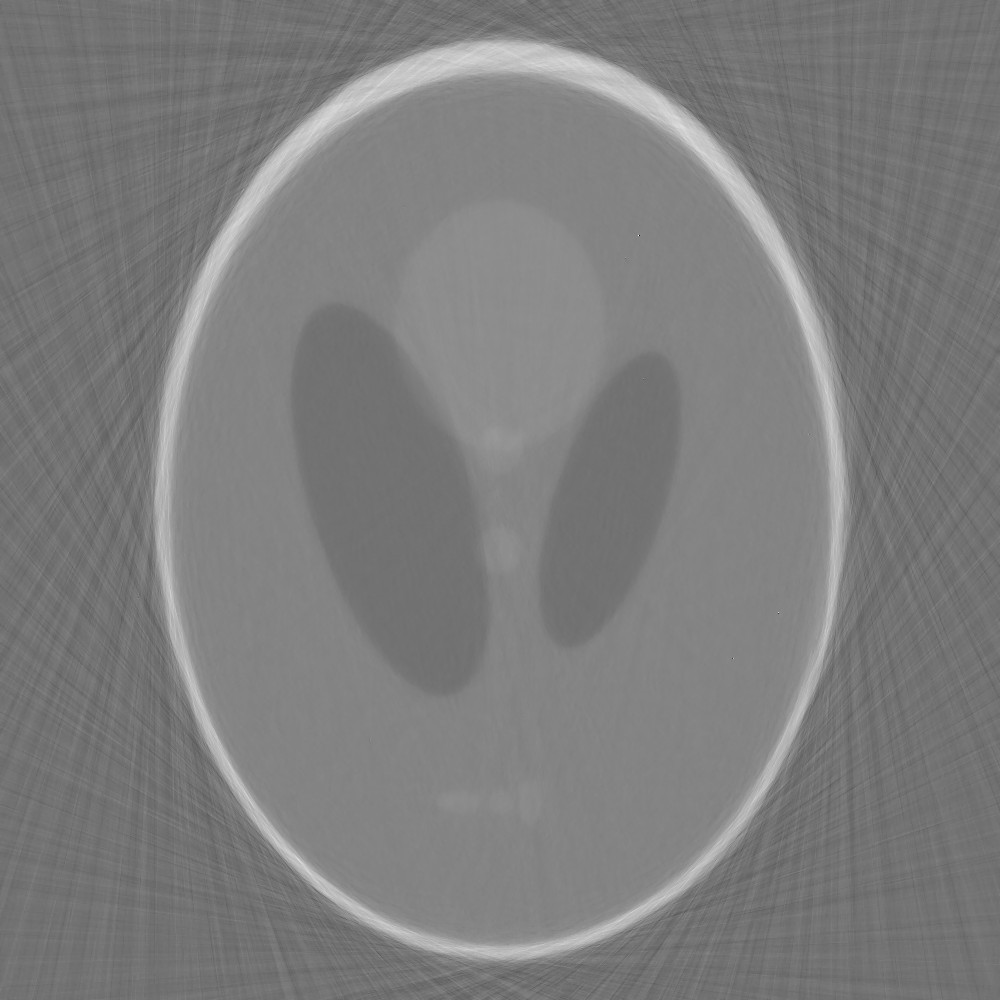} &
  \includegraphics[height = 0.13\textwidth, width = .13\textwidth]{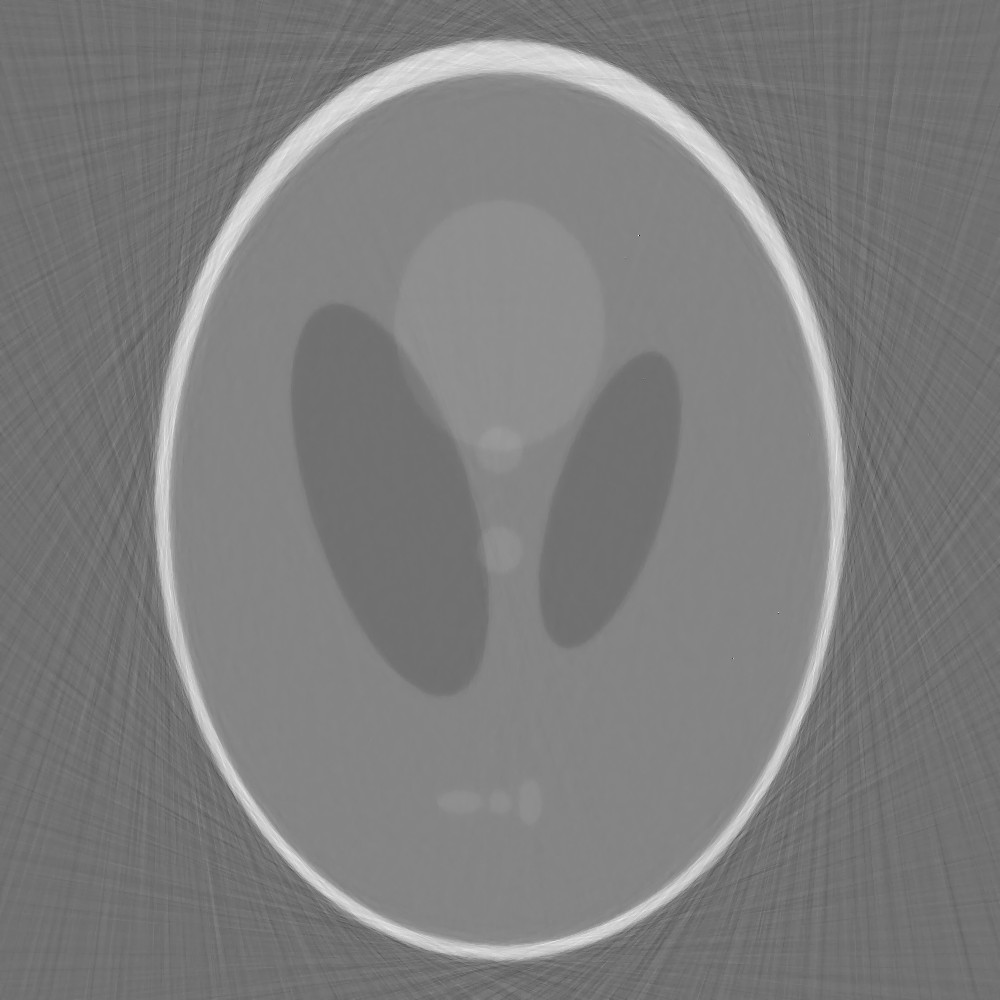} &
  \includegraphics[height = 0.13\textwidth, width = .13\textwidth]{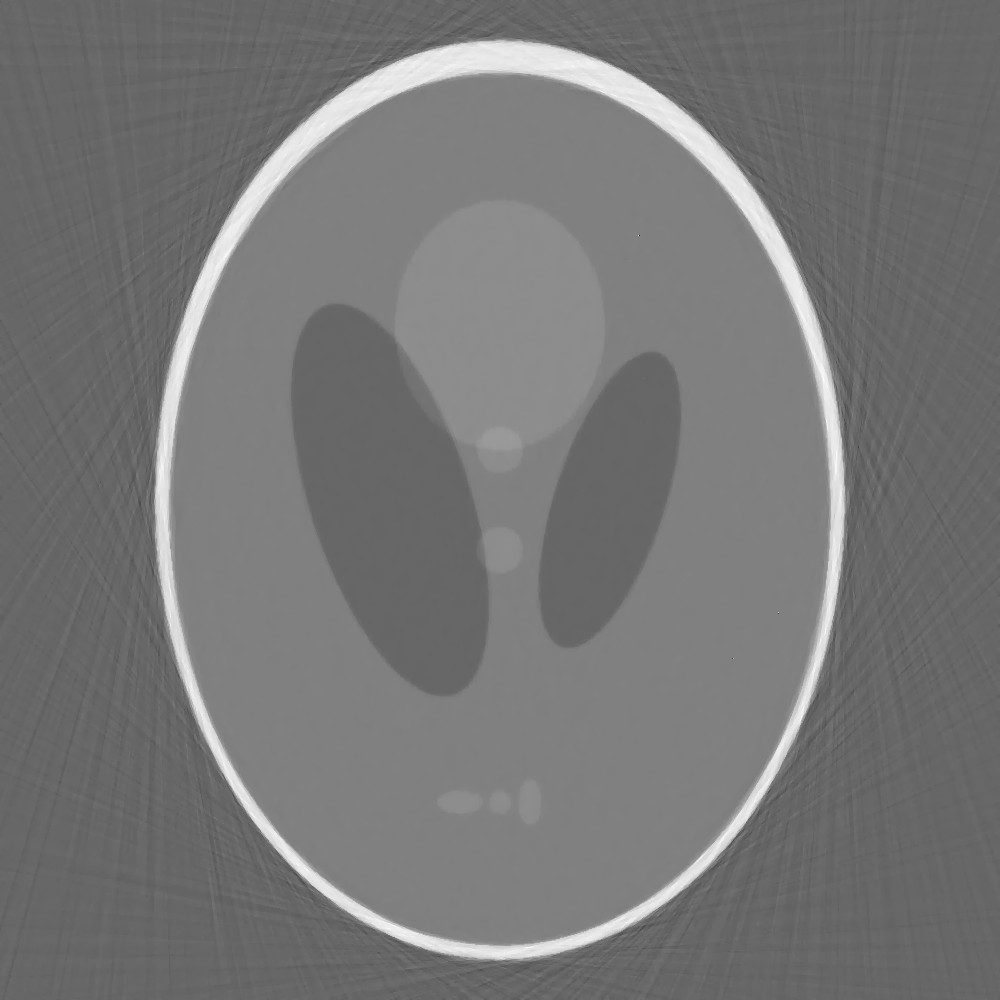}& 
    \includegraphics[height = 0.13\textwidth, width = .13\textwidth]{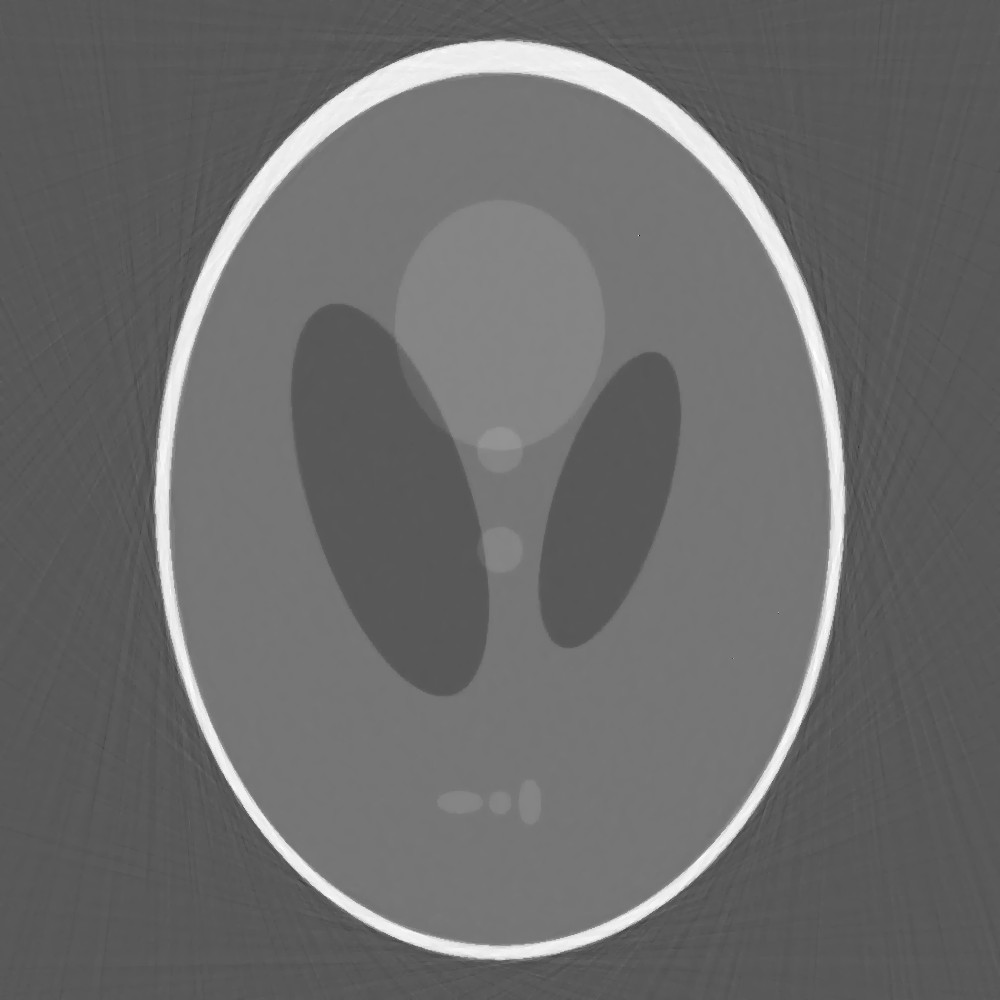} 
    & 
    \includegraphics[height = 0.13\textwidth, width = .13\textwidth]{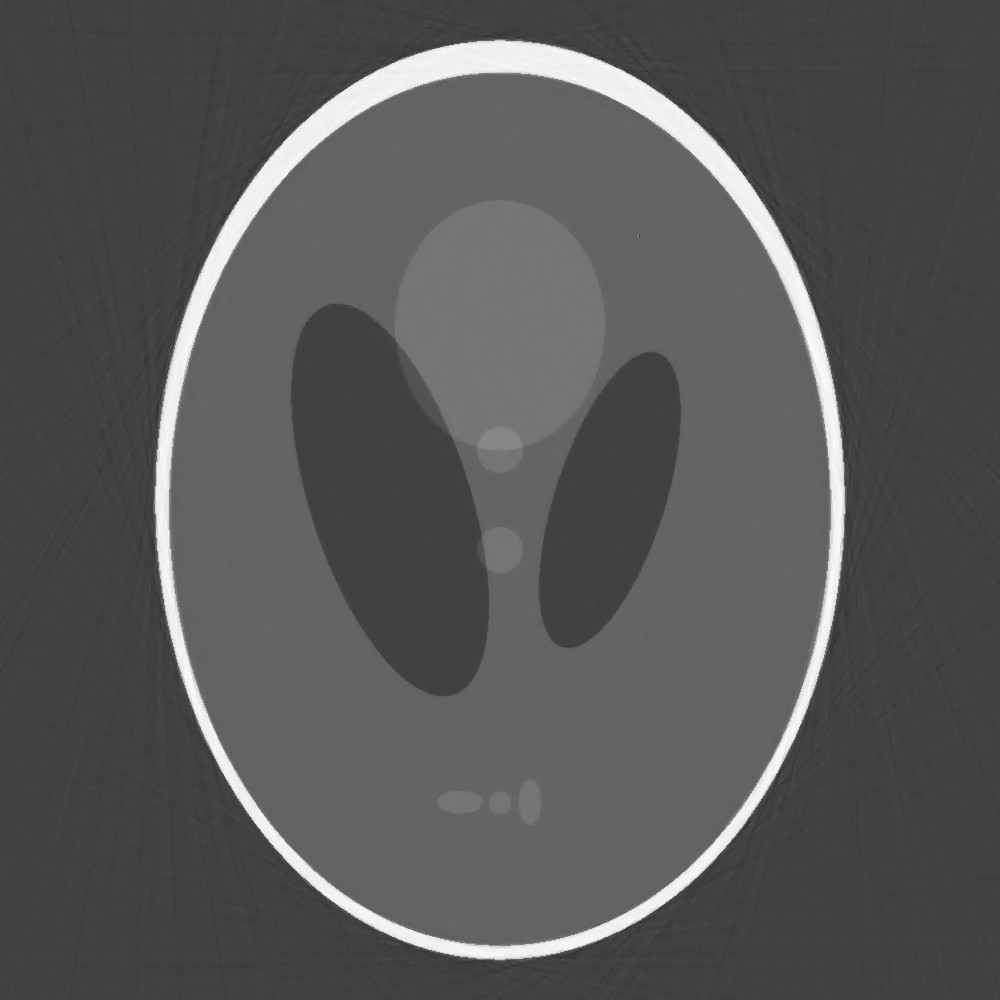} 
        & 
    \includegraphics[height = 0.13\textwidth, width = .13\textwidth]{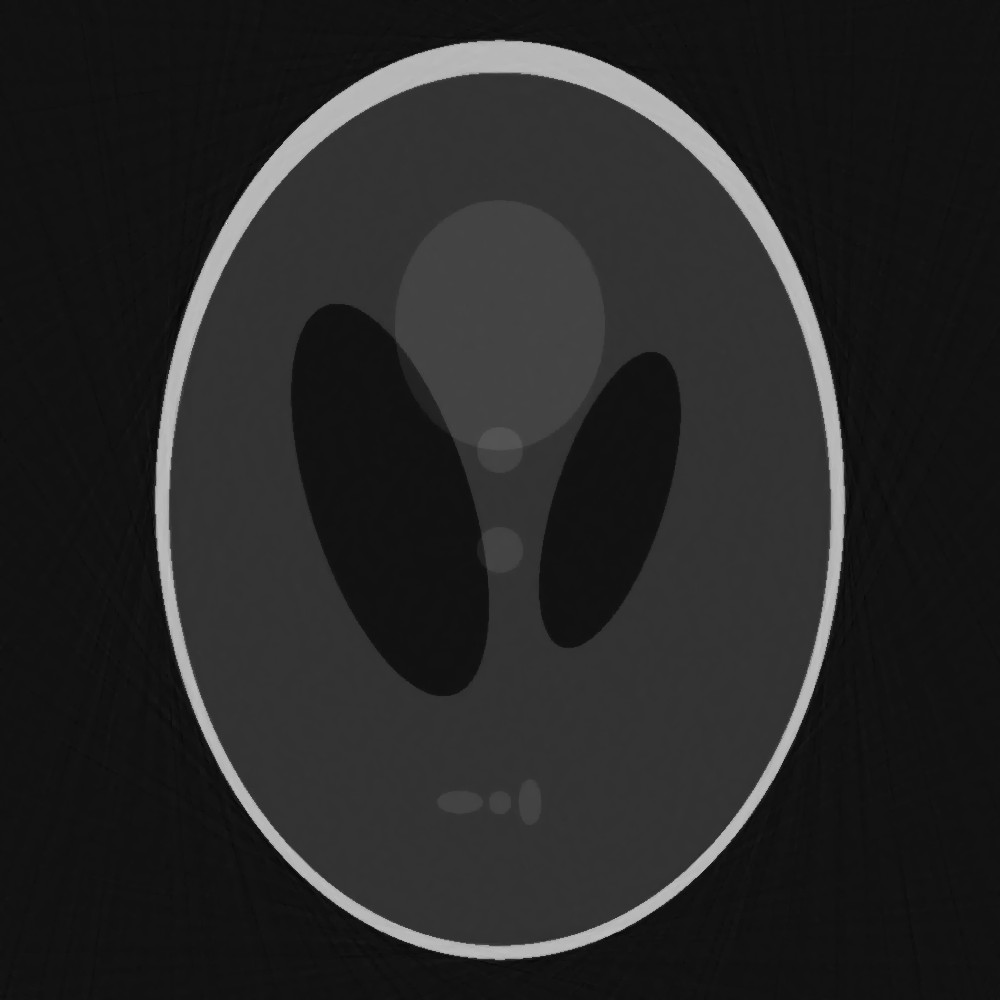} 
\end{tabular}
  \begin{tabular}{cccccc}
  \includegraphics[height = 0.13\textwidth, width = .13\textwidth]{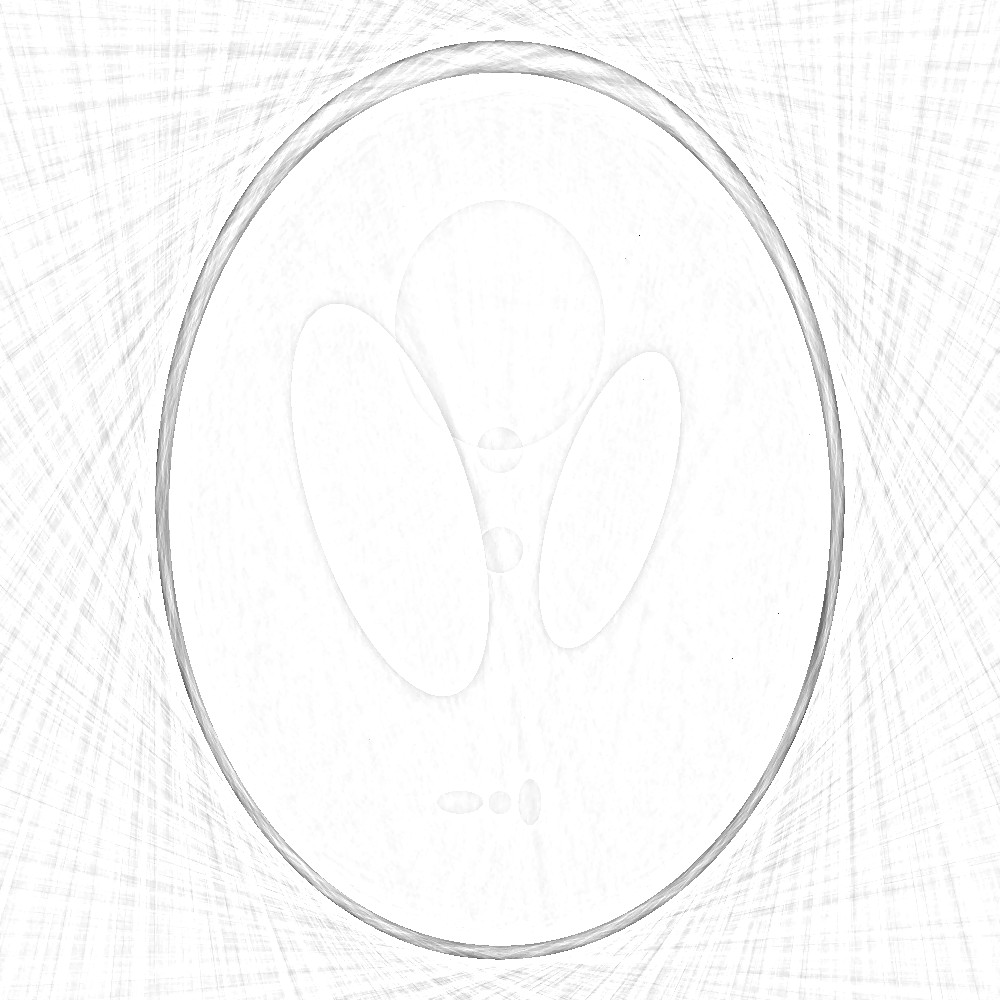} &
  \includegraphics[height = 0.13\textwidth, width = .13\textwidth]{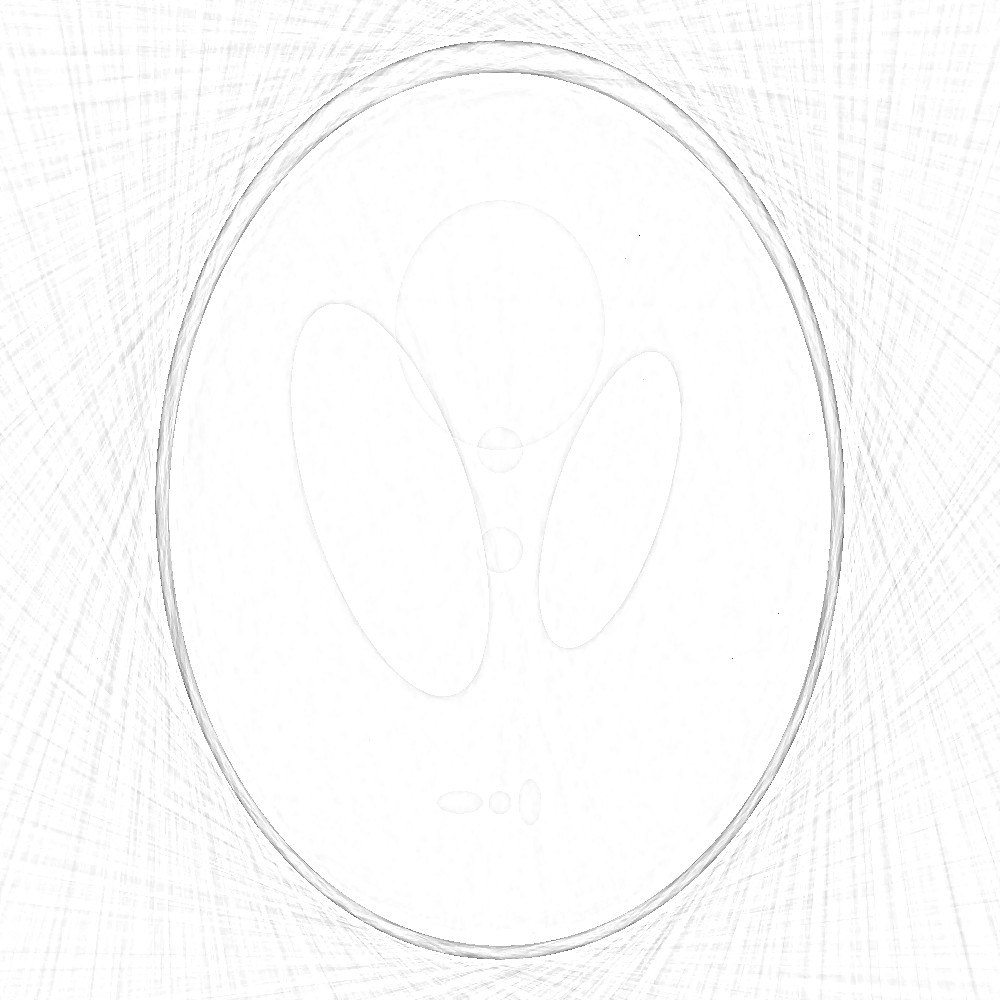} &
  \includegraphics[height = 0.13\textwidth, width = .13\textwidth]{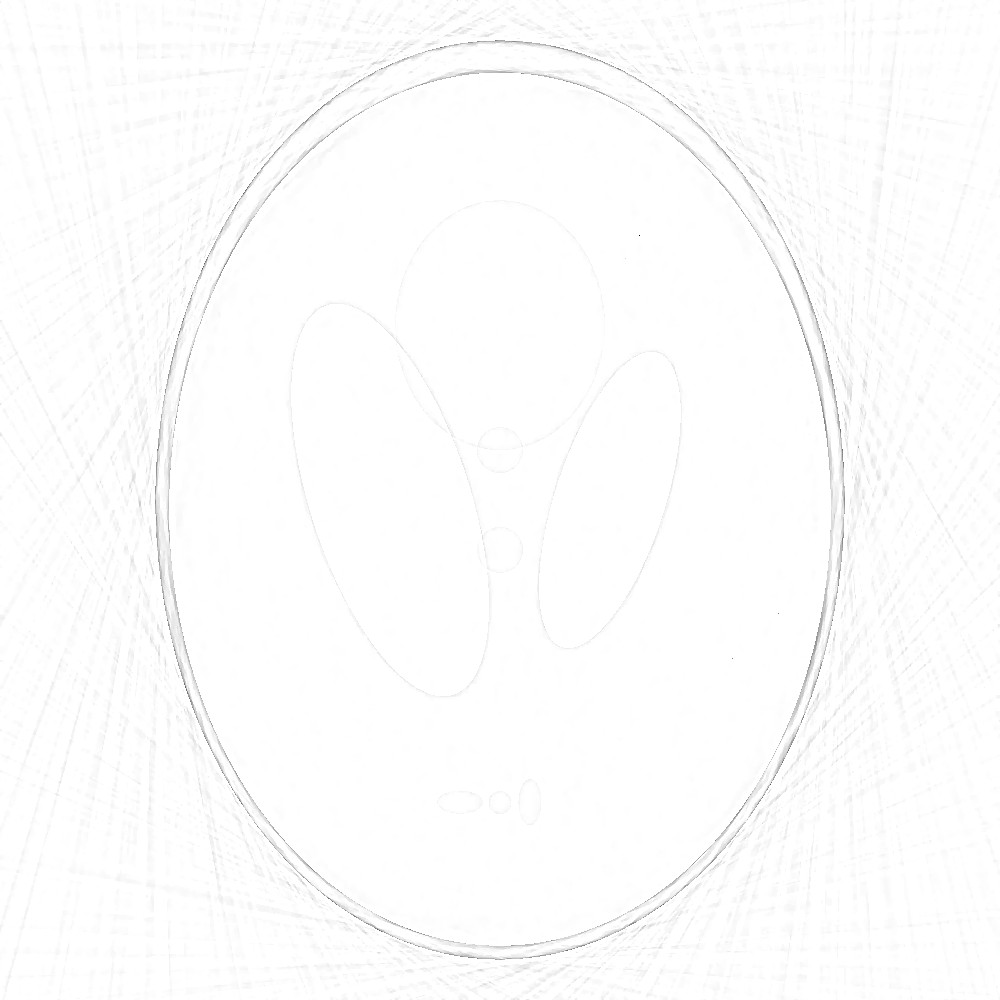}& 
    \includegraphics[height = 0.13\textwidth, width = .13\textwidth]{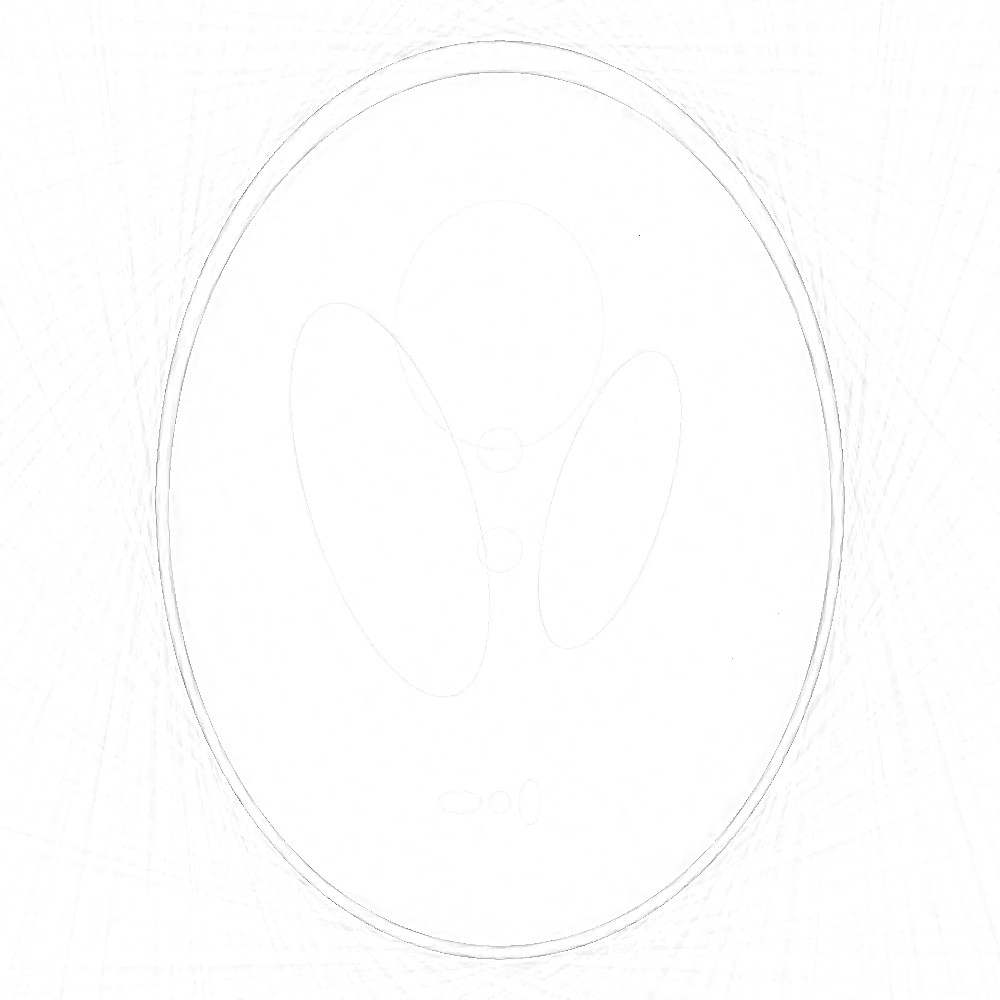} 
    & 
    \includegraphics[height = 0.13\textwidth, width = .13\textwidth]{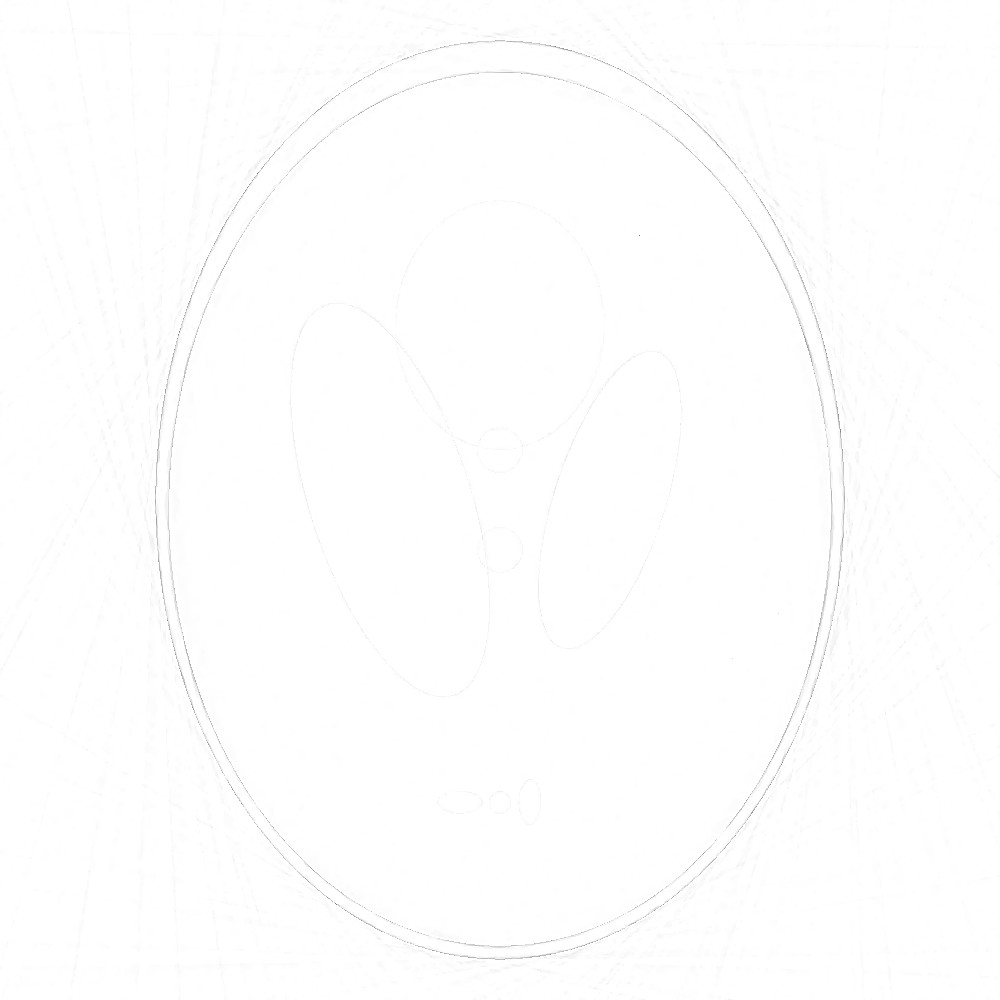} 
        & 
    \includegraphics[height = 0.13\textwidth, width = .13\textwidth]{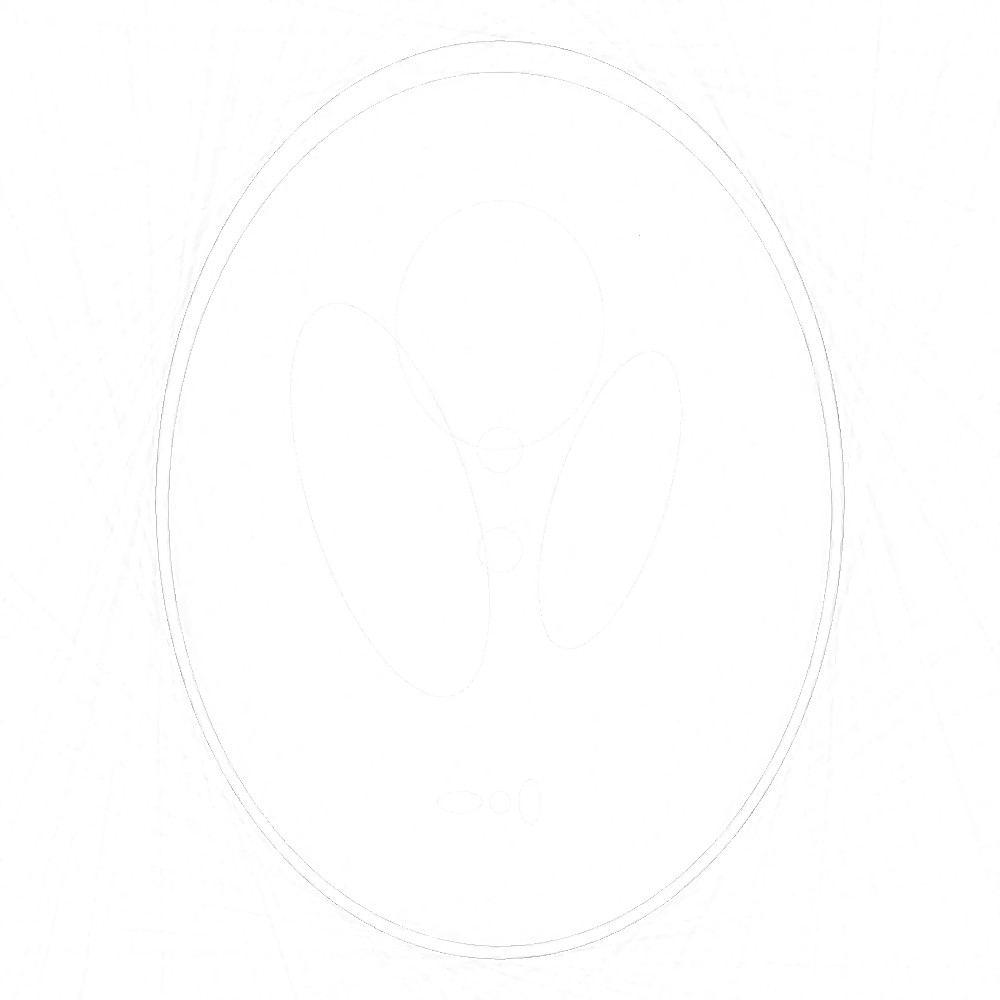} 
\end{tabular}
\begin{tabular}{cccccc}
$31.84\% \qquad$ & $\rm 23.78\% \qquad$& $\rm 16.97\% \quad$& \qquad $\rm  12.87\%$&$ \qquad \rm 9.96\%$&$ \qquad \rm 7.88\%$\\
\end{tabular}
  \caption{Test 2: Streaming tomography example. Reconstructions for s-RMM-GKS where each problem is ran until the stopping criteria is satisfied.}
  \label{fig: Streaming6prob_reconstructions_sRMM-GKS}
\end{figure}

\subsection{Dynamic photoacoustic tomography (PAT)}
An emerging hybrid imaging modality that shows great potential for pre-clinical research is known as photoacoustic tomography. PAT combines the rich contrast of optical imaging with the high resolution of ultrasound imaging aiming to produce high resolution images with lower cost and less side effects than other imaging modalities. 
In this example we consider a discrete PAT problem where we let $\bx_{i}\in \R^{n_x\times n_y}$ be the discretized desired solution at time instance $i$ at locations of transducers $\bd_i$, $i = 1, 2, \dots, n_t$. This yields a dynamic model. It is not the first time that PAT is considered in the dynamic framework. For instance, it has been already considered in \cite{chung2017motion,chung2018efficient,lucka2018enhancing}. At each location we assume there are $r$ radii such that the forward operator $\bA_{i} \in \R^{r \times n\cdot n_t}$, where $n = n_x \cdot n_y$. Hence, the measurements are obtained from spherical projections as
\begin{equation}\label{eq: lin_eq_spherical}
\bd_{i} = \bA_{i} \bx_{i} + \be_i,
\end{equation}
with $\be_i$ being Gaussian noise. The goal is to estimate approximate images $\bx_{i}$ when are give the observations $\bd_i$.

\begin{figure}[ht!]
\centering
\begin{tabular}{cccccc}
     $t =1$ & \qquad  $t =10$ & \qquad  $t =20$ &   \qquad  $t =30$ & \qquad $t =40$ &  \qquad $t =50$  \\
\end{tabular}
\begin{tabular}{cccccc}
\includegraphics[height = 0.13\textwidth, width = .13\textwidth]{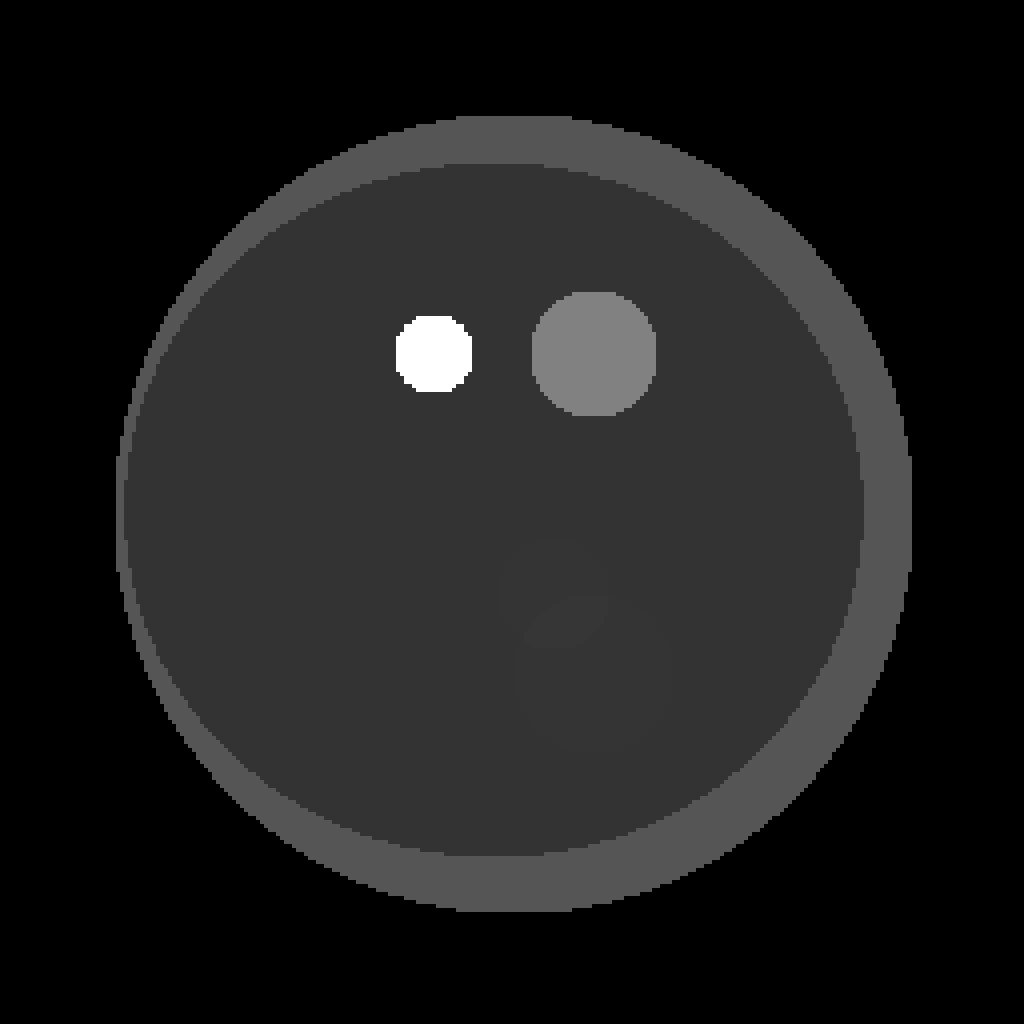} &
\includegraphics[height = 0.13\textwidth, width = .13\textwidth]{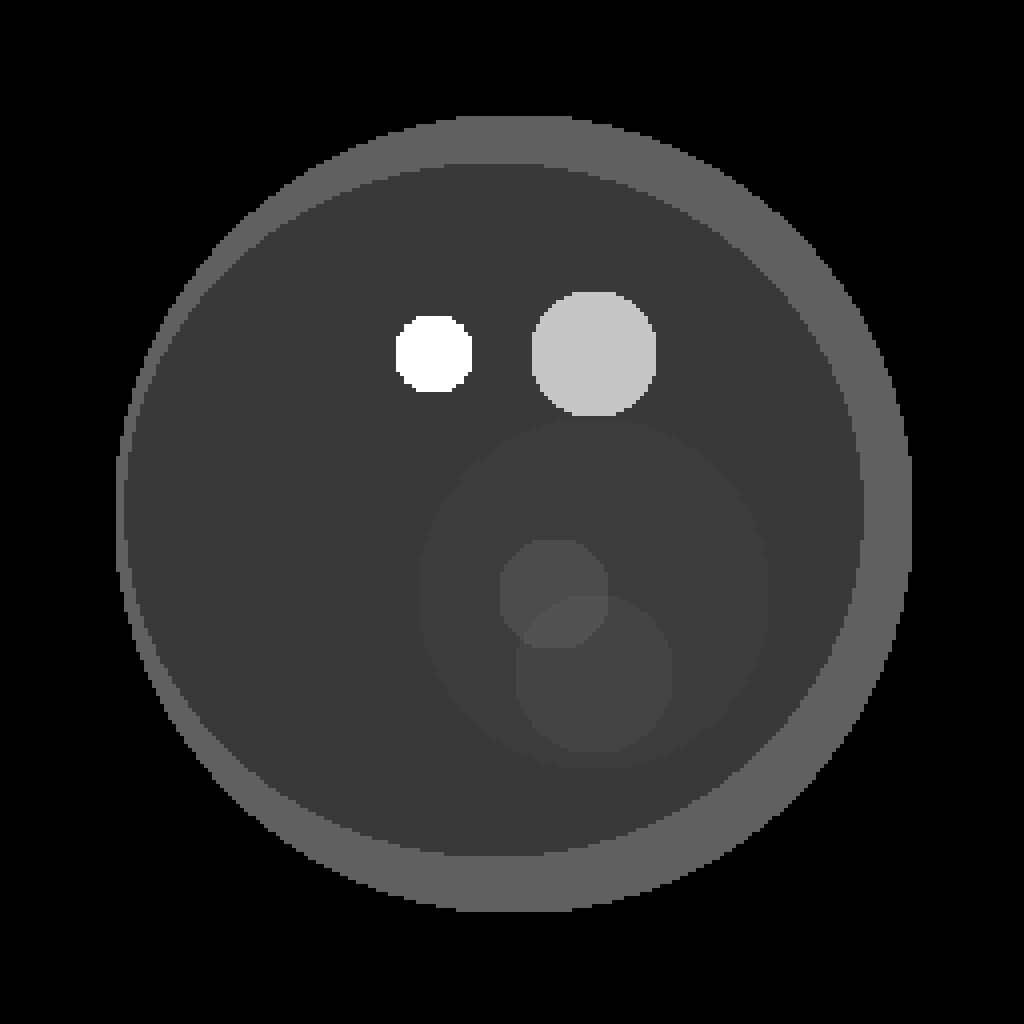} &
\includegraphics[height = 0.13\textwidth, width = .13\textwidth]{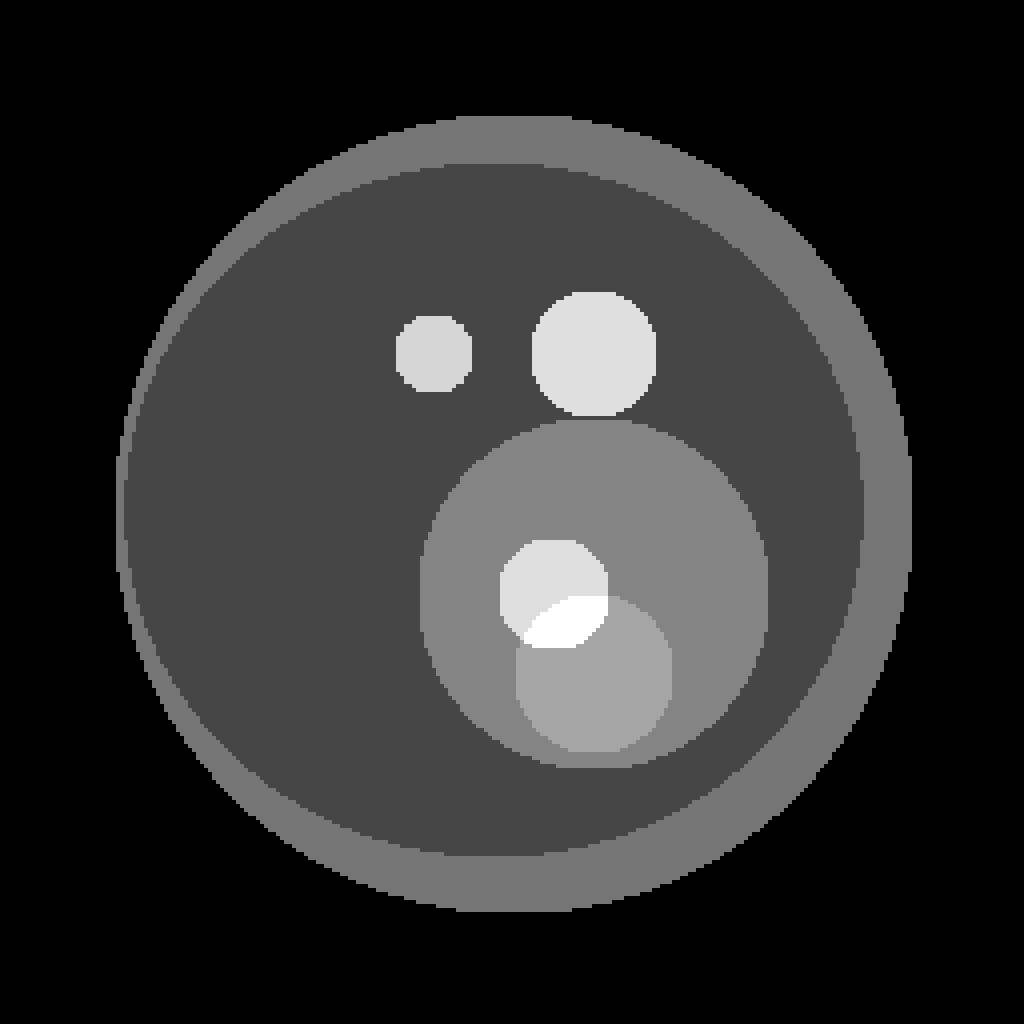}& \includegraphics[height = 0.13\textwidth, width = .13\textwidth]{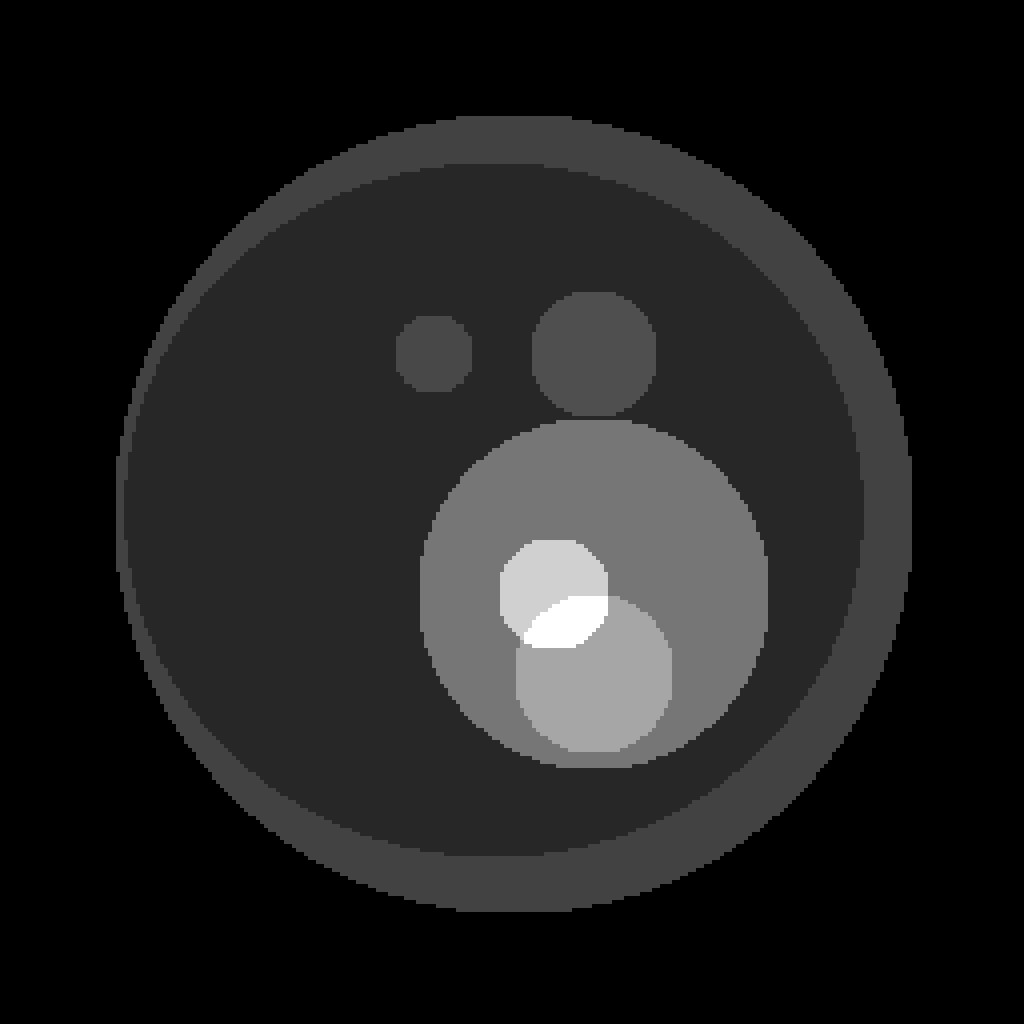} &
\includegraphics[height = 0.13\textwidth, width = .13\textwidth]{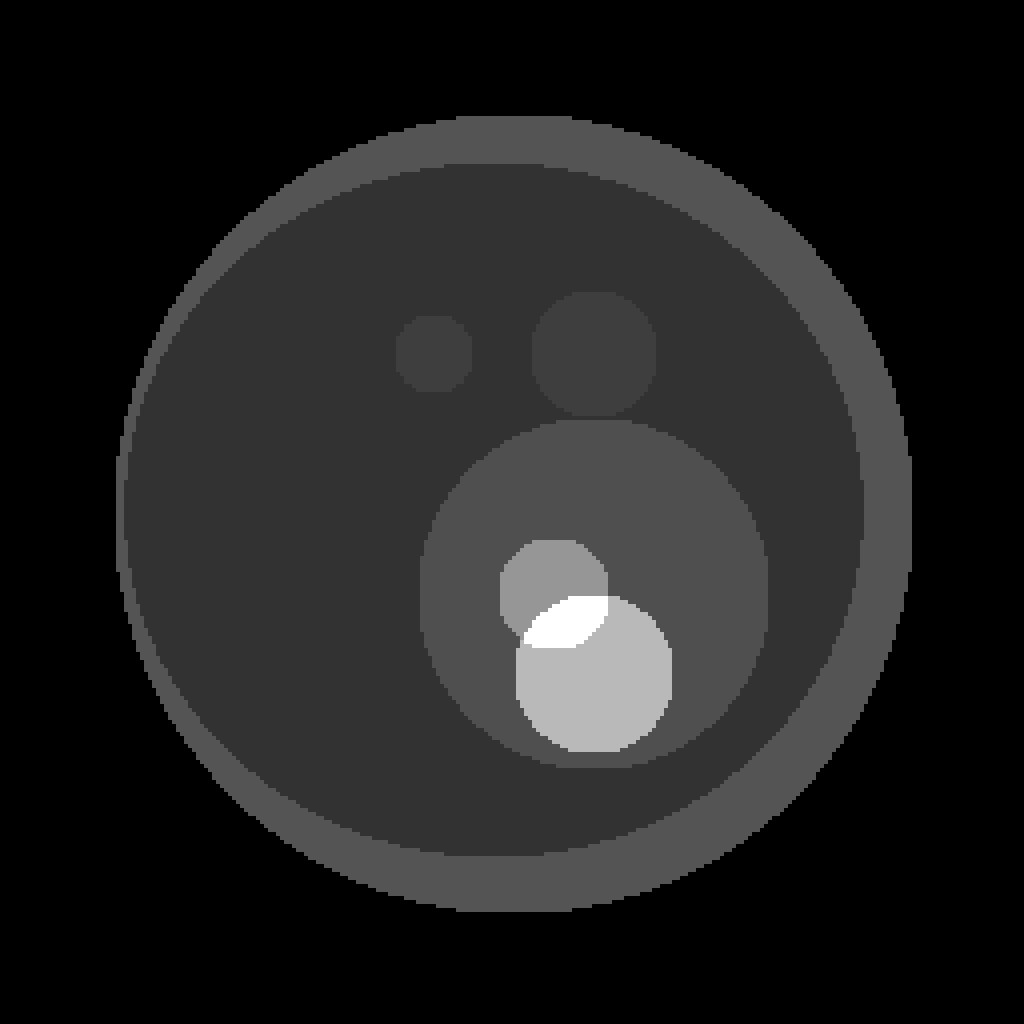}& 
\includegraphics[height = 0.13\textwidth, width = .13\textwidth]{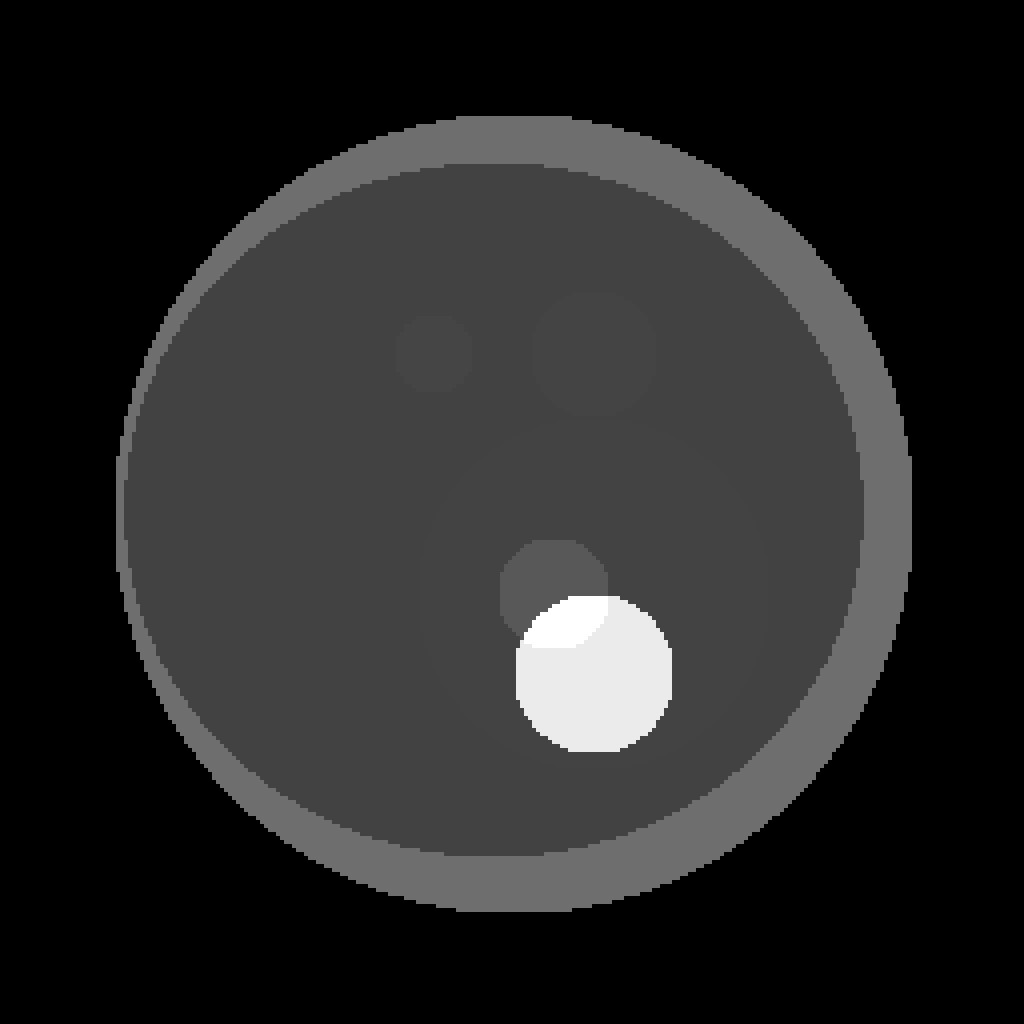}
\end{tabular}
\begin{tabular}{cccccc}
\includegraphics[height = 0.13\textwidth, width = .13\textwidth]{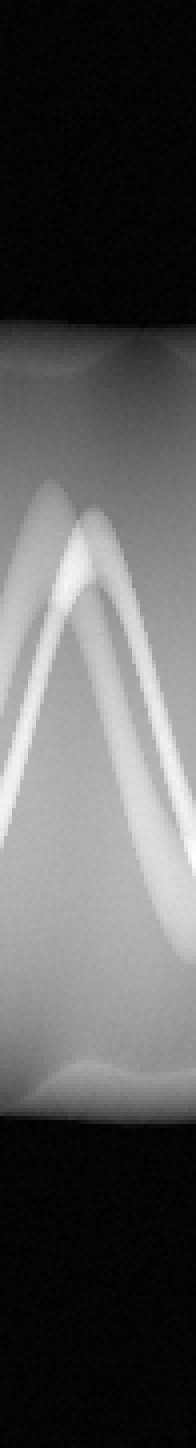} &
\includegraphics[height = 0.13\textwidth, width = .13\textwidth]{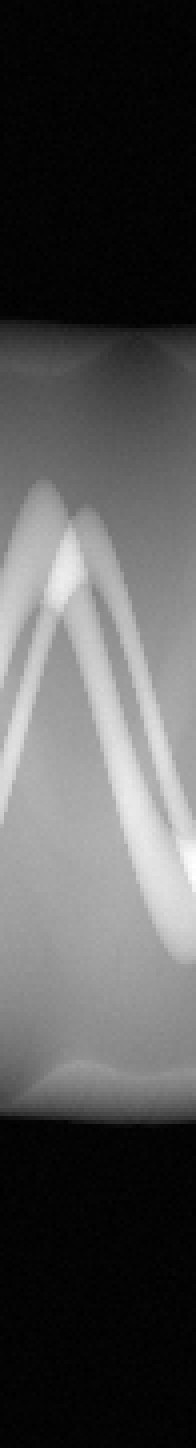} &   
\includegraphics[height = 0.13\textwidth, width = .13\textwidth]{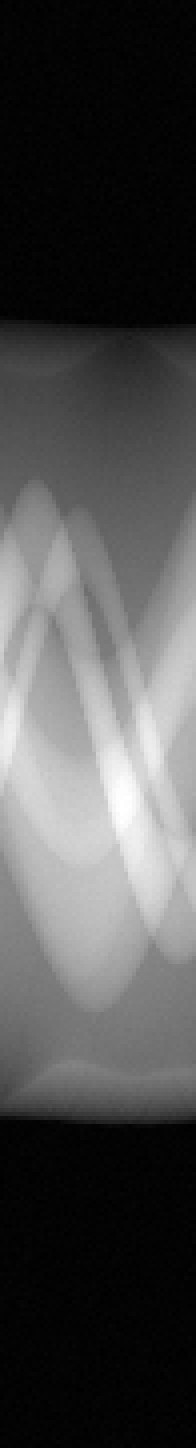} &
\includegraphics[height = 0.13\textwidth, width = .13\textwidth]{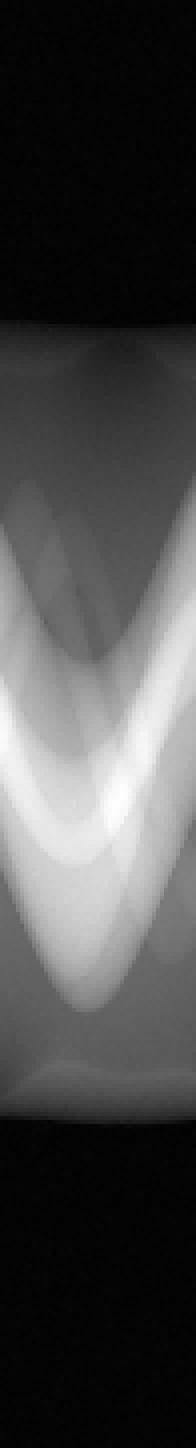} &
\includegraphics[height = 0.13\textwidth, width = .13\textwidth]{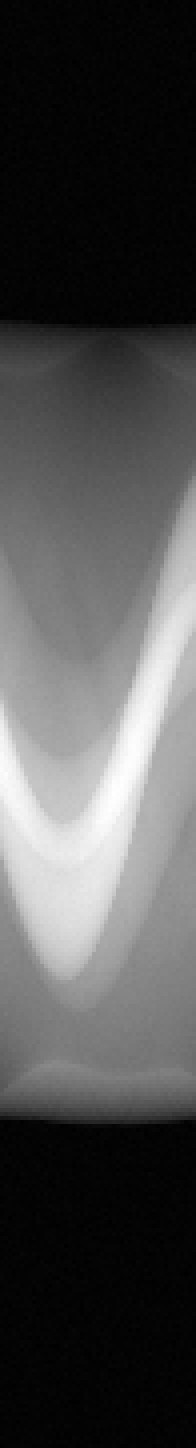} &
\includegraphics[height = 0.13\textwidth, width = .13\textwidth]{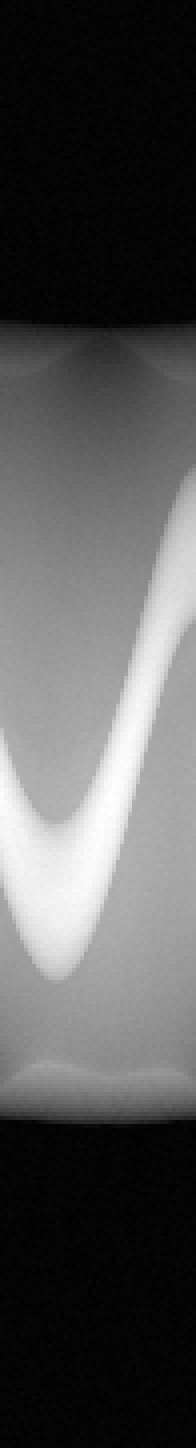} 
\end{tabular}
\caption{PAT test problem. First row: true images of size $256\times 256$, Second row: sinograms of size $362 \cdot 45$ at times steps $1, 10, 20, 30, 40, 50$, from left to right respectively.}
\label{Fig: PATtrueandsino}
\end{figure}

\begin{figure}[h!]
	\centering
	\begin{minipage}{0.15\textwidth}
		\centering
		\includegraphics[width=\textwidth]{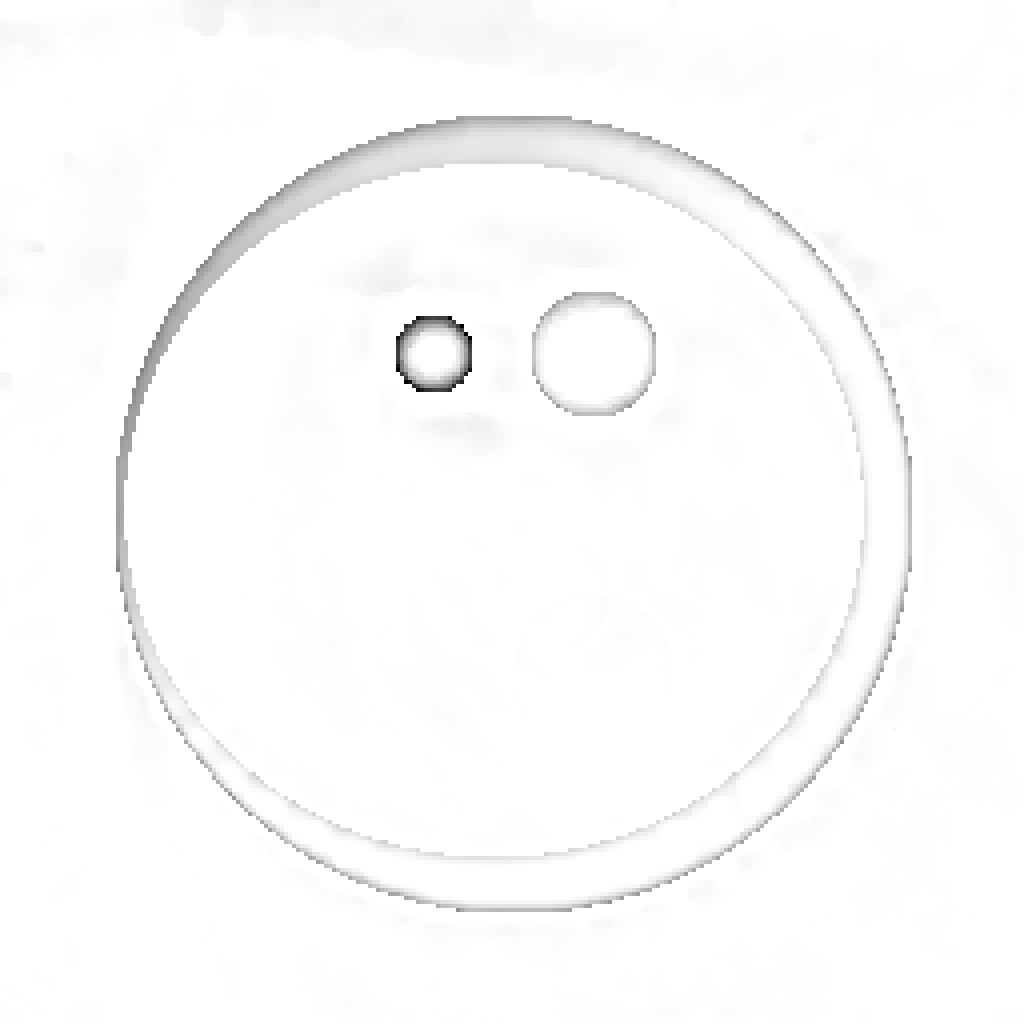}\\
	\end{minipage}
	\begin{minipage}{0.15\textwidth}
		\centering
		\includegraphics[width=\textwidth]{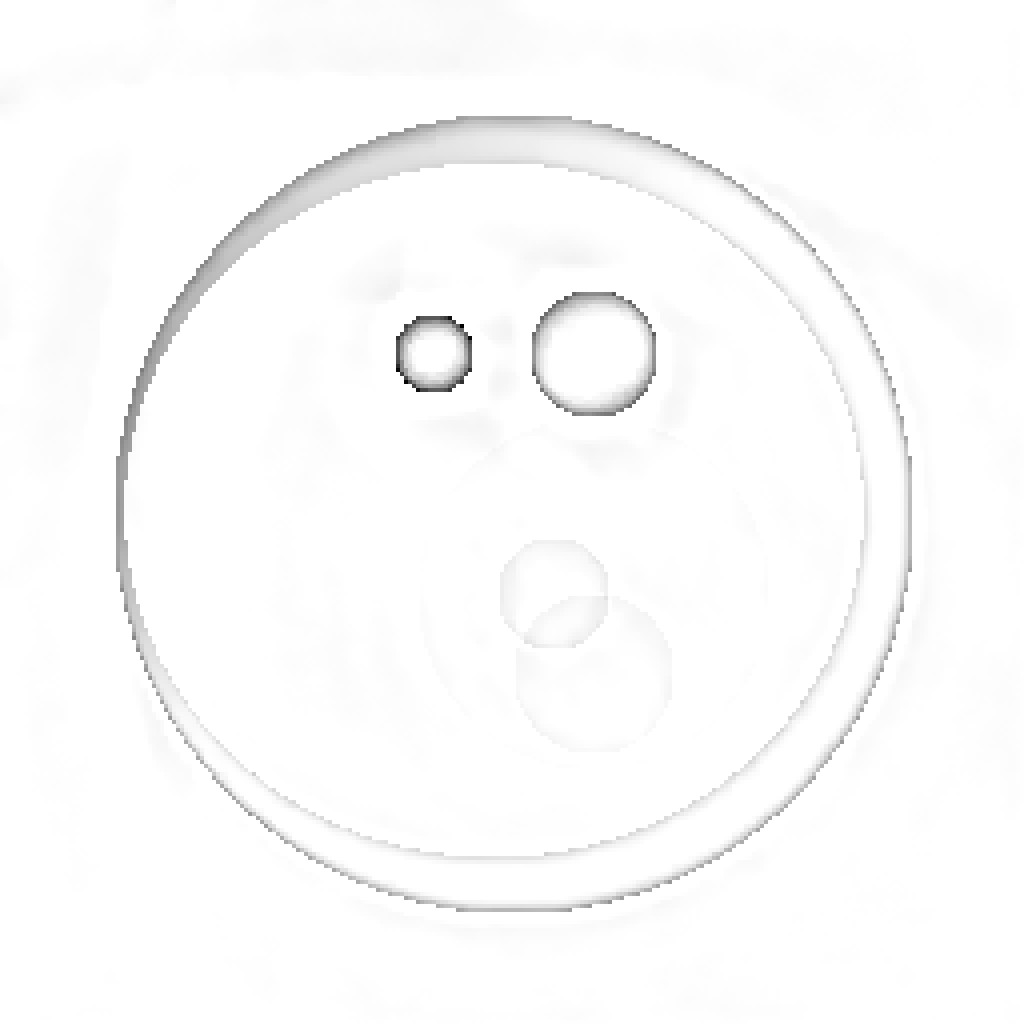}\\
	\end{minipage}
	\begin{minipage}{0.15\textwidth}
		\centering
		\includegraphics[width=\textwidth]{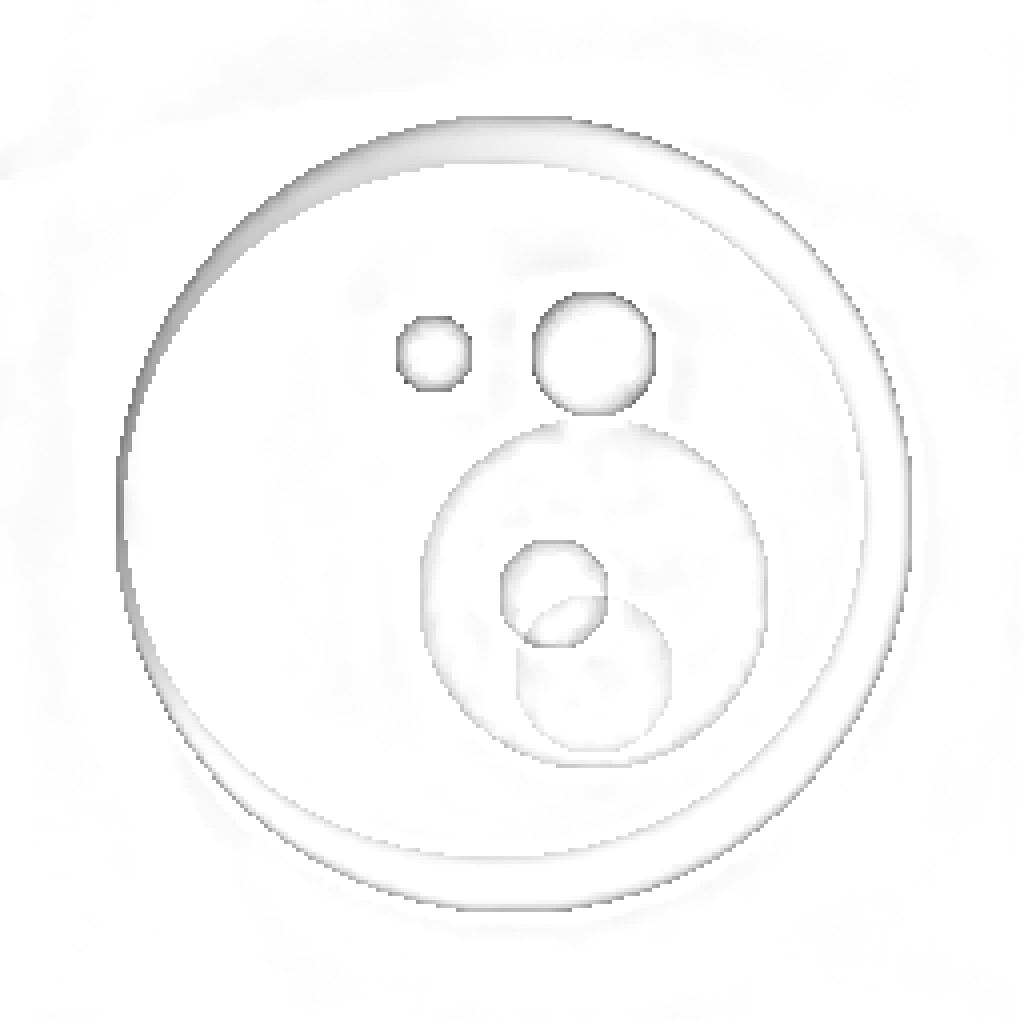}\\
	\end{minipage}
		\begin{minipage}{0.15\textwidth}
		\centering
		\includegraphics[width=\textwidth]{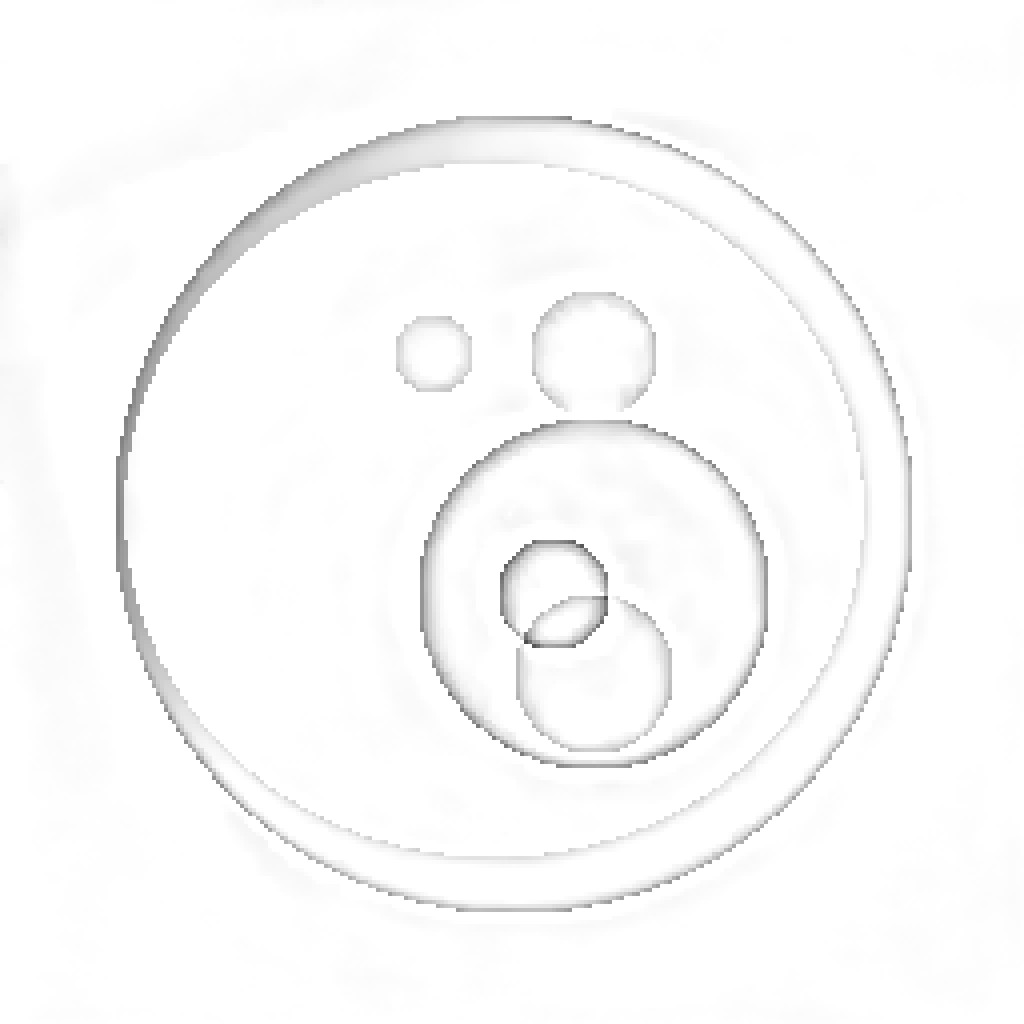}\\
	\end{minipage}
 \begin{minipage}{0.15\textwidth}
		\centering
		\includegraphics[width=\textwidth]{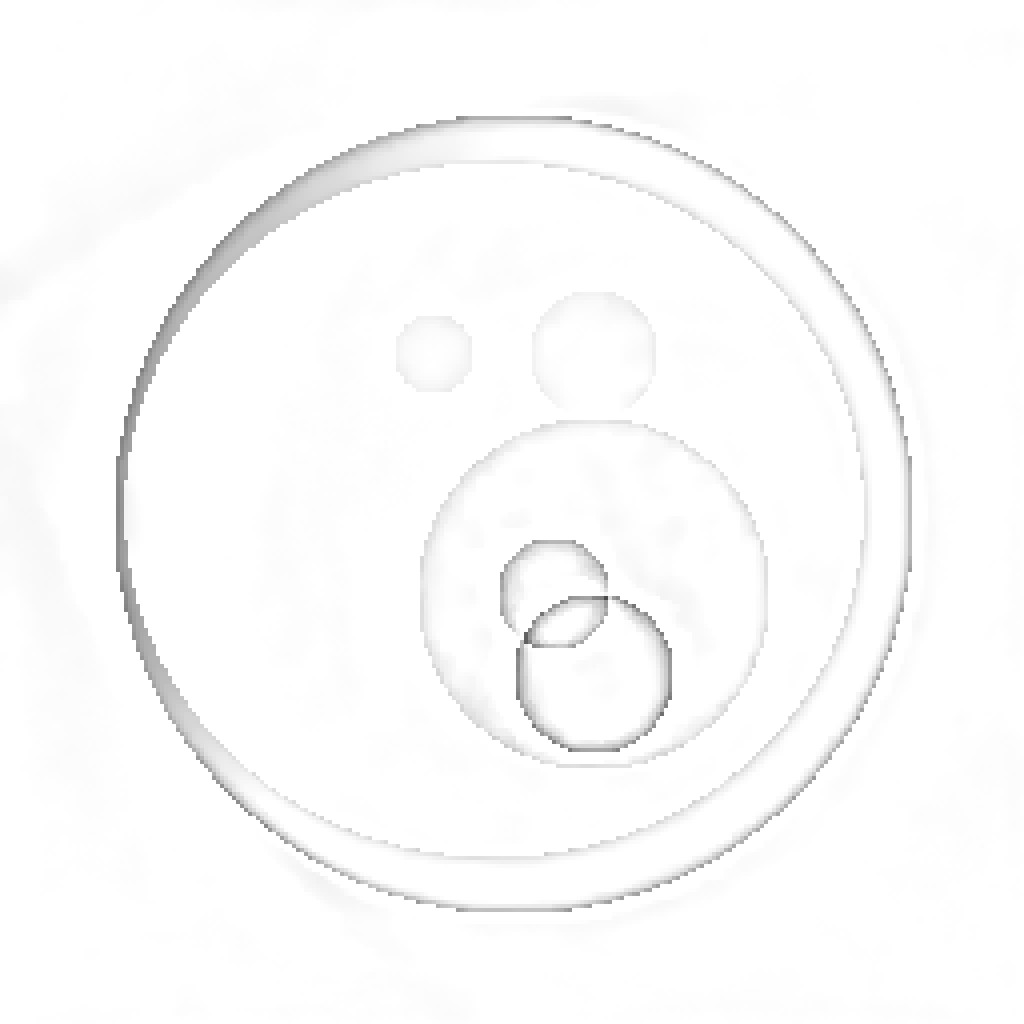}\\
	\end{minipage}
	\begin{minipage}{0.15\textwidth}
		\centering
		\includegraphics[width=\textwidth]{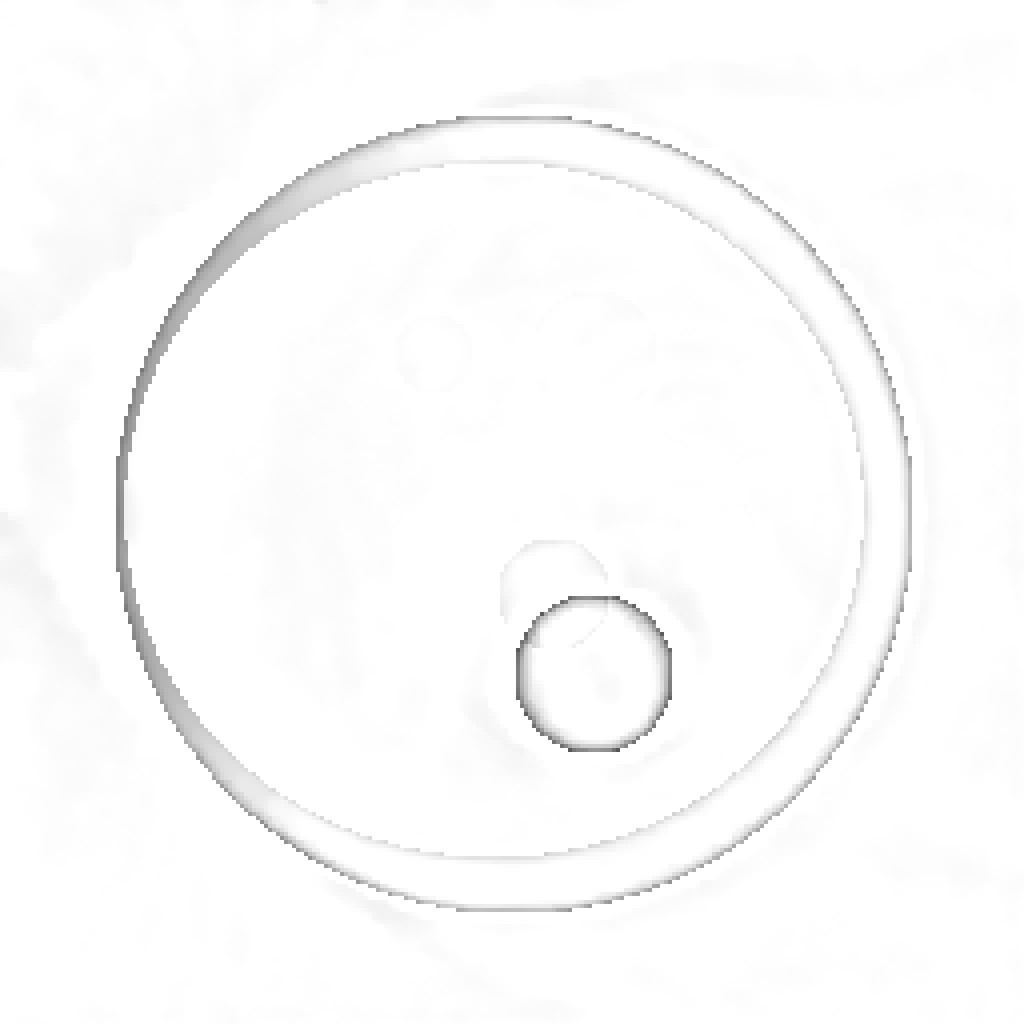}\\
	\end{minipage}
 
	\begin{minipage}{0.15\textwidth}
		\centering
		\includegraphics[width=\textwidth]{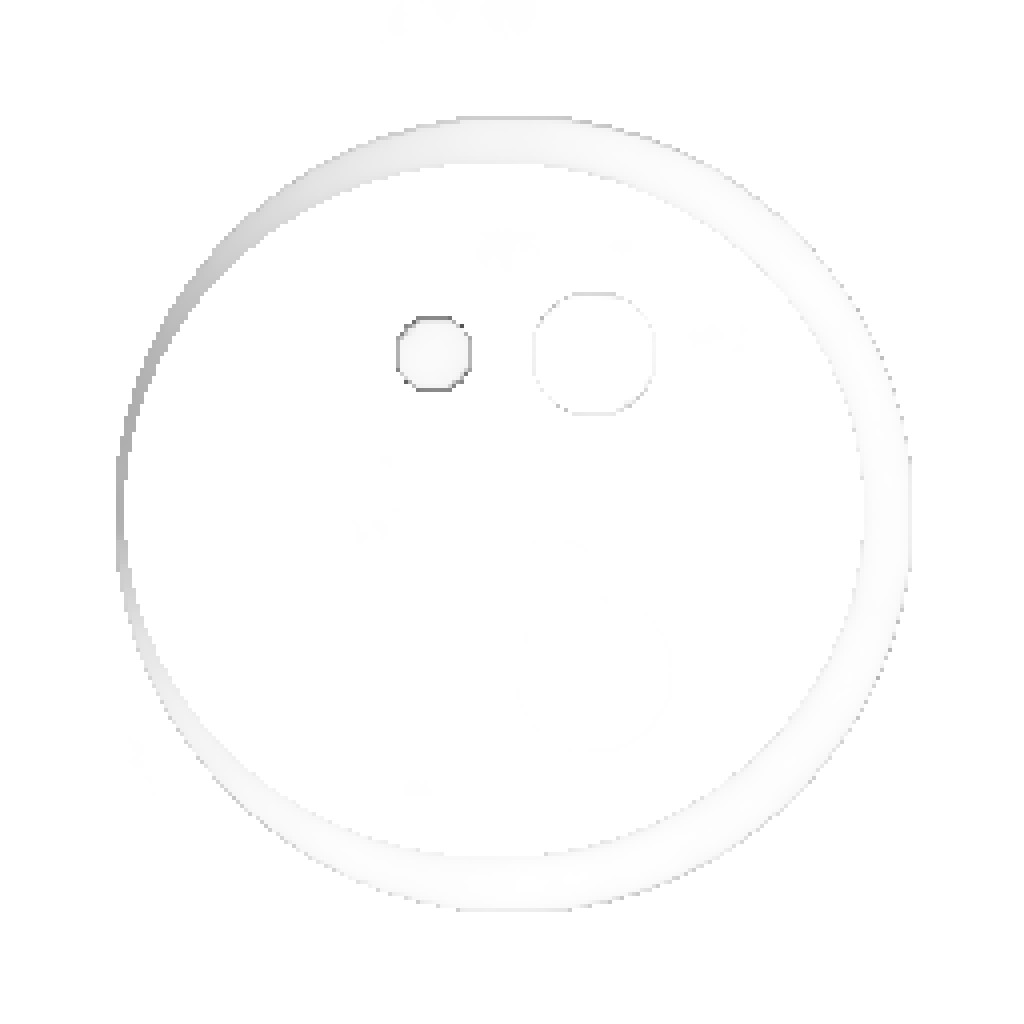}\\
	\end{minipage}
	\begin{minipage}{0.15\textwidth}
		\centering
		\includegraphics[width=\textwidth]{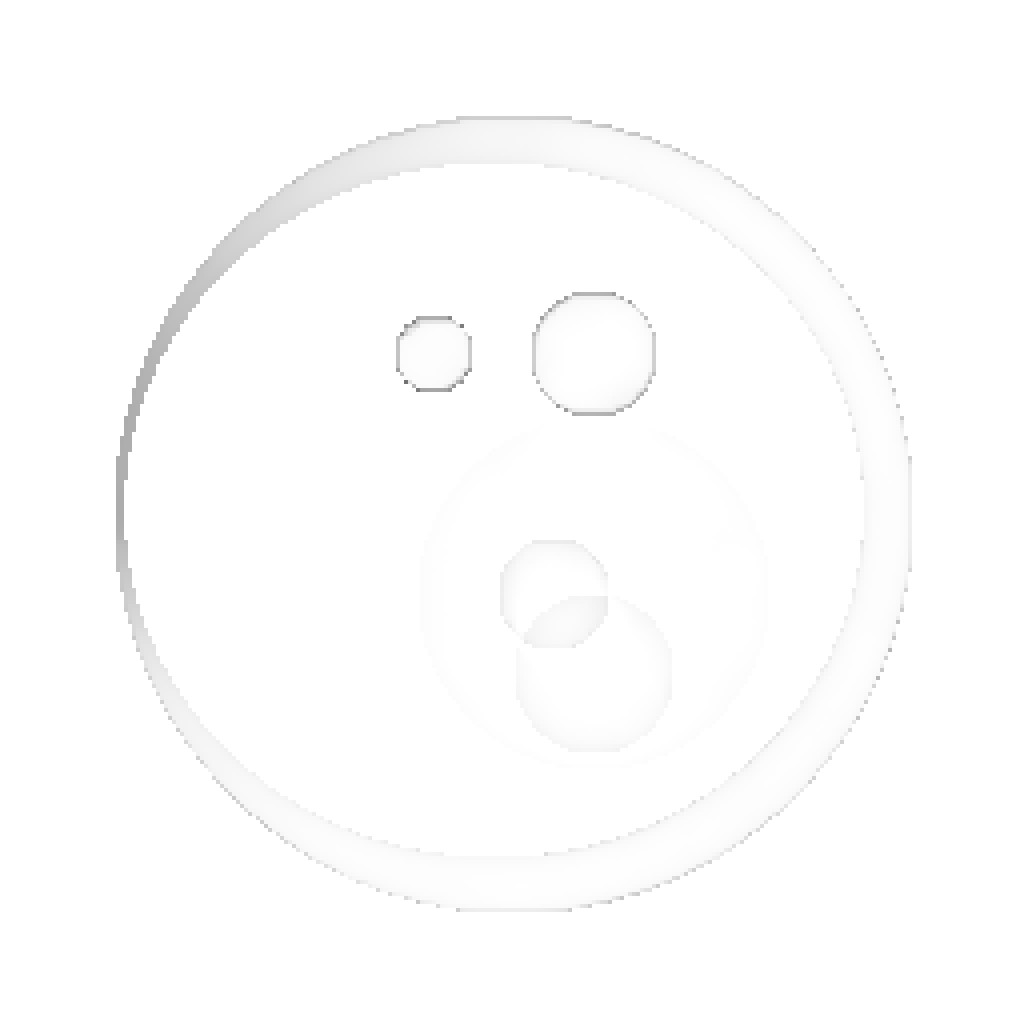}\\
	\end{minipage}
	\begin{minipage}{0.15\textwidth}
		\centering
		\includegraphics[width=\textwidth]{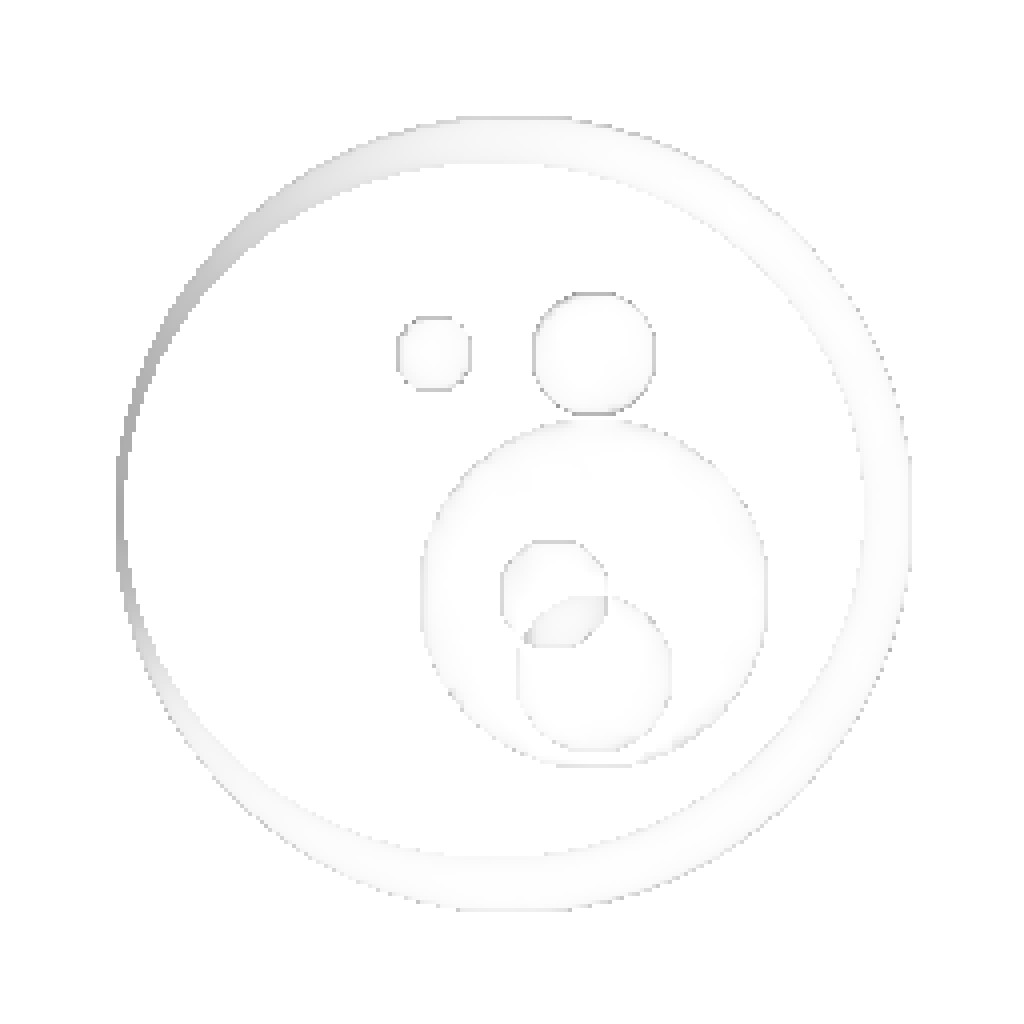}\\
	\end{minipage}
		\begin{minipage}{0.15\textwidth}
		\centering
		\includegraphics[width=\textwidth]{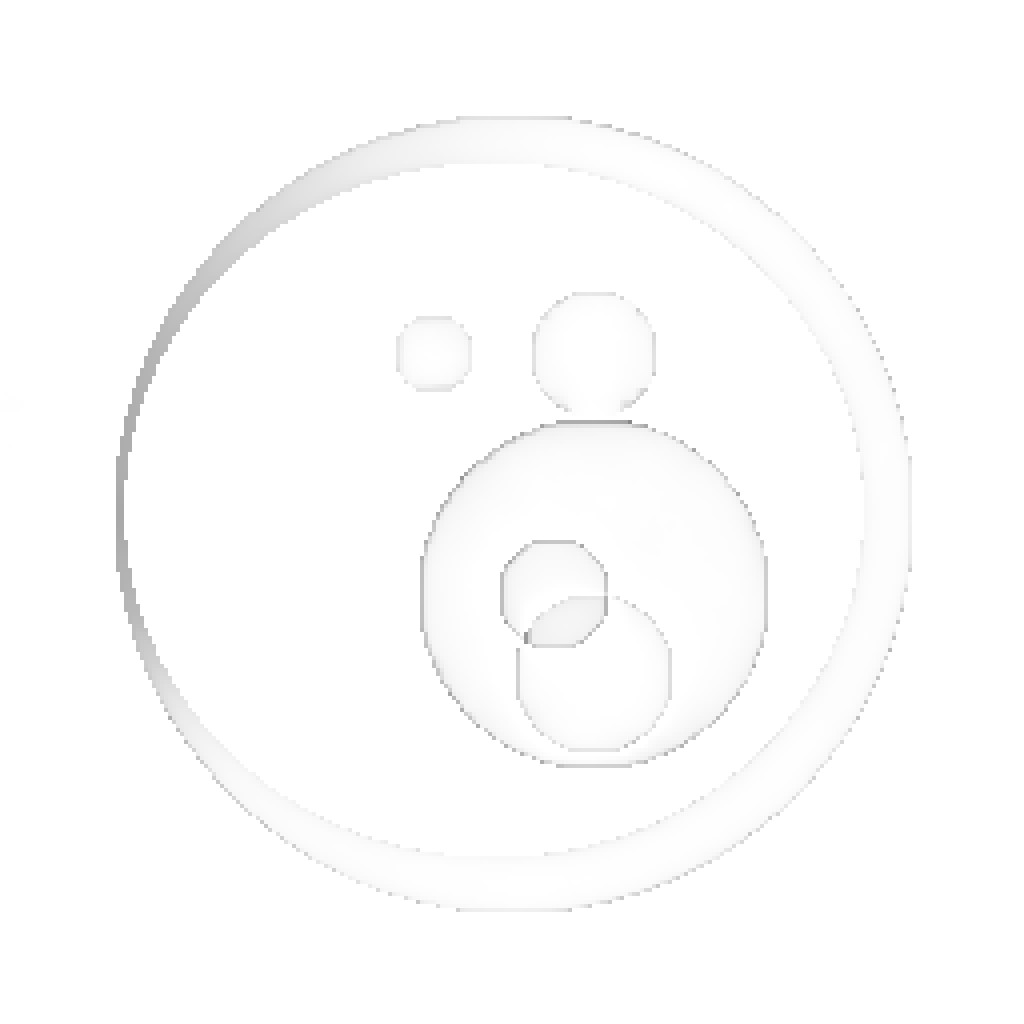}\\
	\end{minipage}
		\begin{minipage}{0.15\textwidth}
		\centering
		\includegraphics[width=\textwidth]{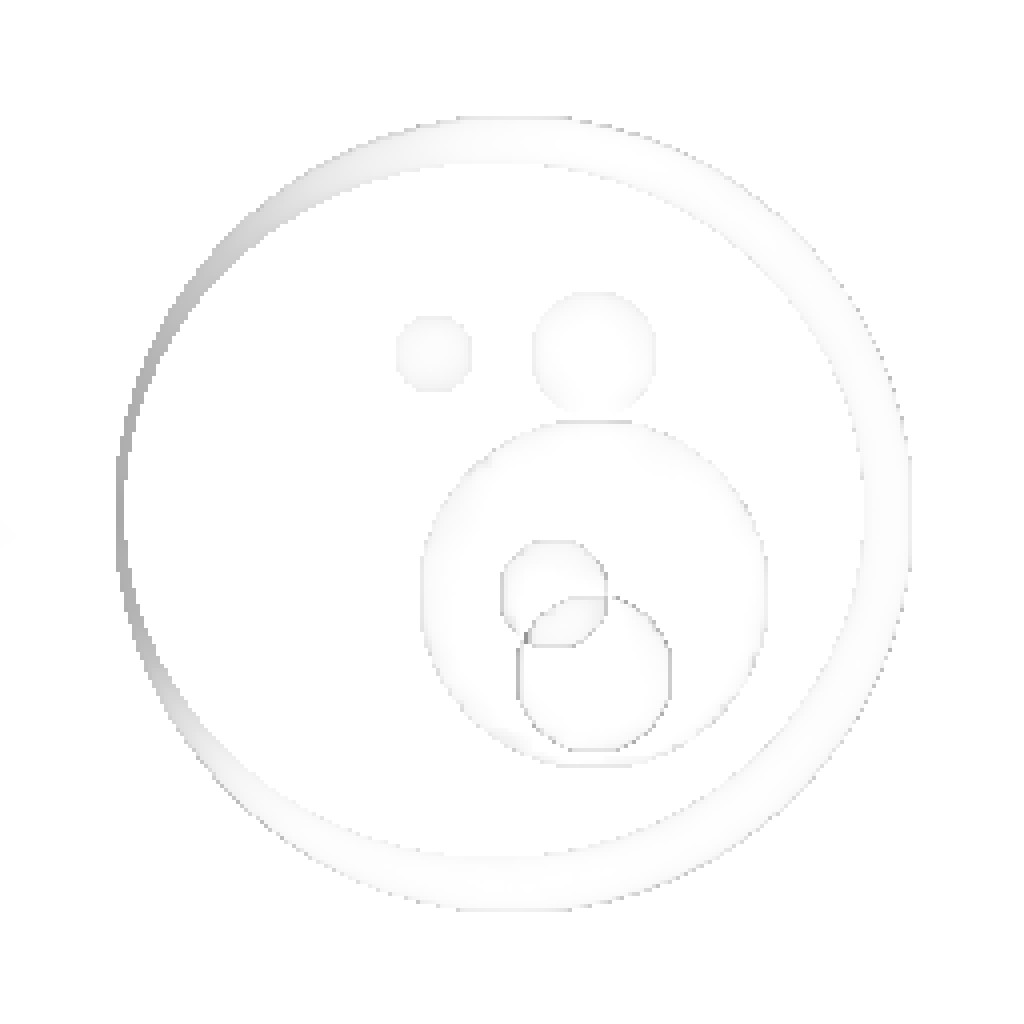}\\
	\end{minipage}
	\begin{minipage}{0.15\textwidth}
		\centering
		\includegraphics[width=\textwidth]{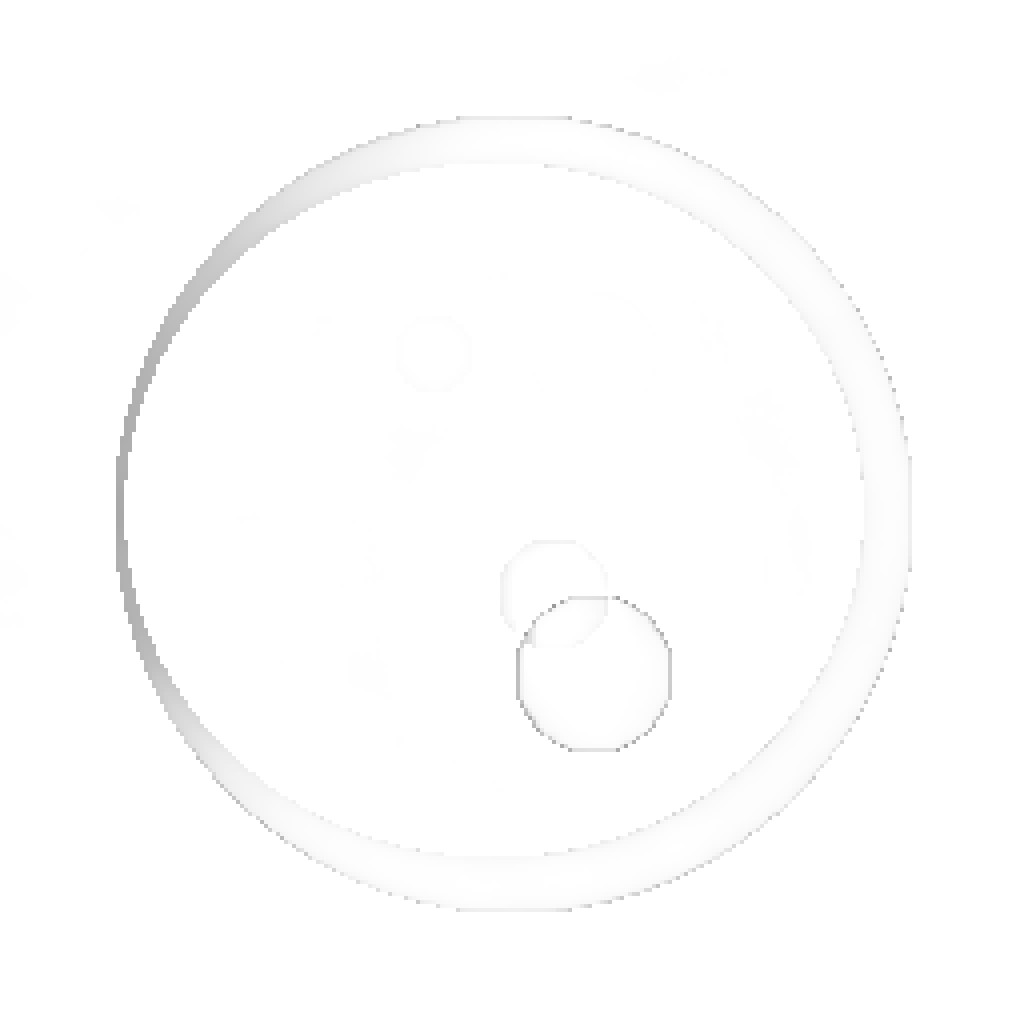}\\
	\end{minipage}

 	\begin{minipage}{0.15\textwidth}
		\centering
		\includegraphics[width=\textwidth]{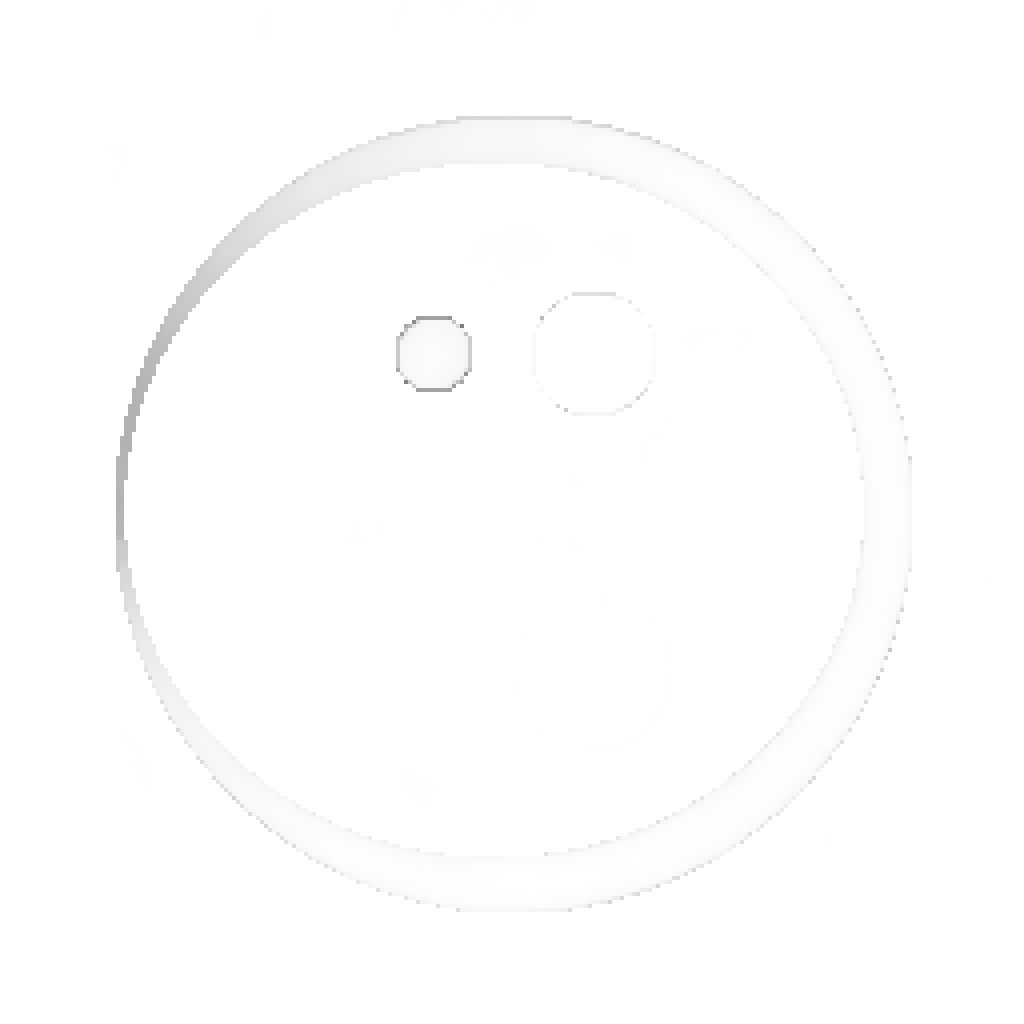}\\
	\end{minipage}
	\begin{minipage}{0.15\textwidth}
		\centering
		\includegraphics[width=\textwidth]{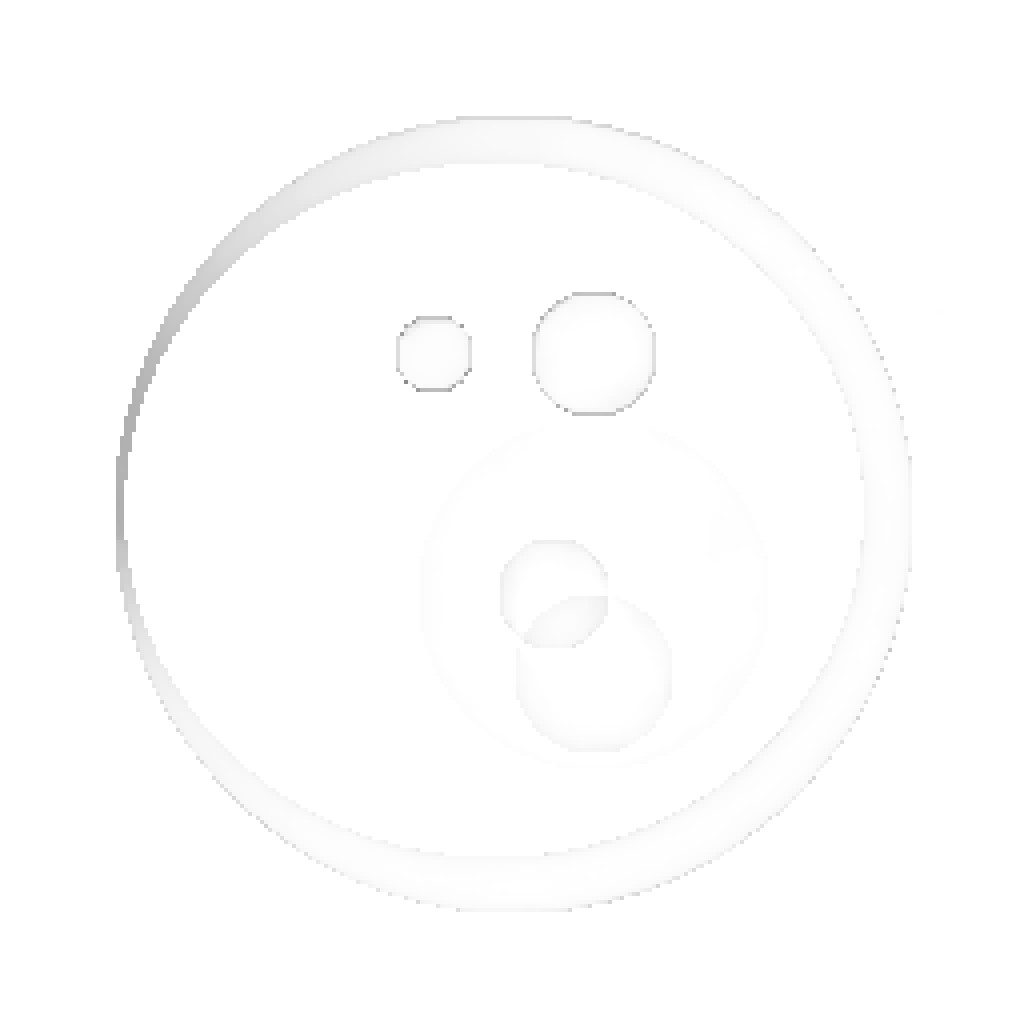}\\
	\end{minipage}
	\begin{minipage}{0.15\textwidth}
		\centering
		\includegraphics[width=\textwidth]{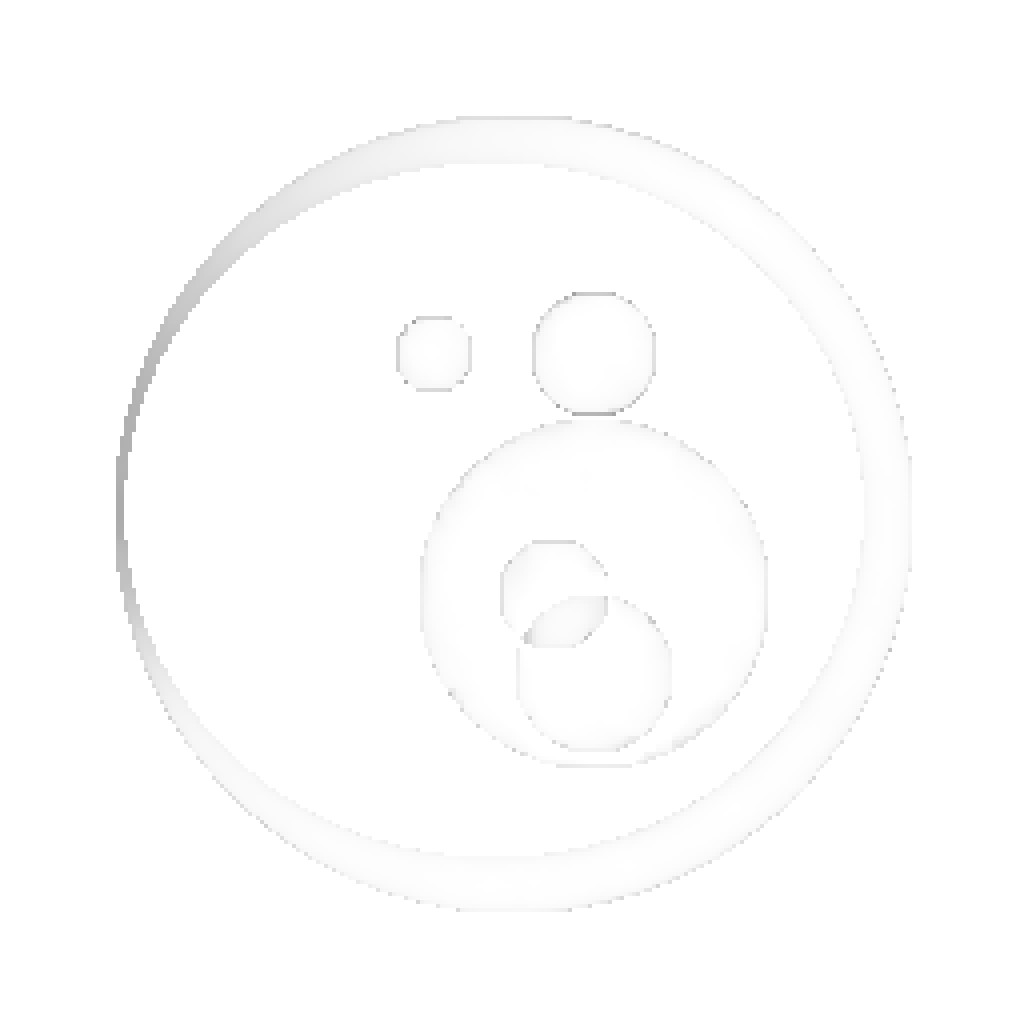}\\
	\end{minipage}
		\begin{minipage}{0.15\textwidth}
		\centering
		\includegraphics[width=\textwidth]{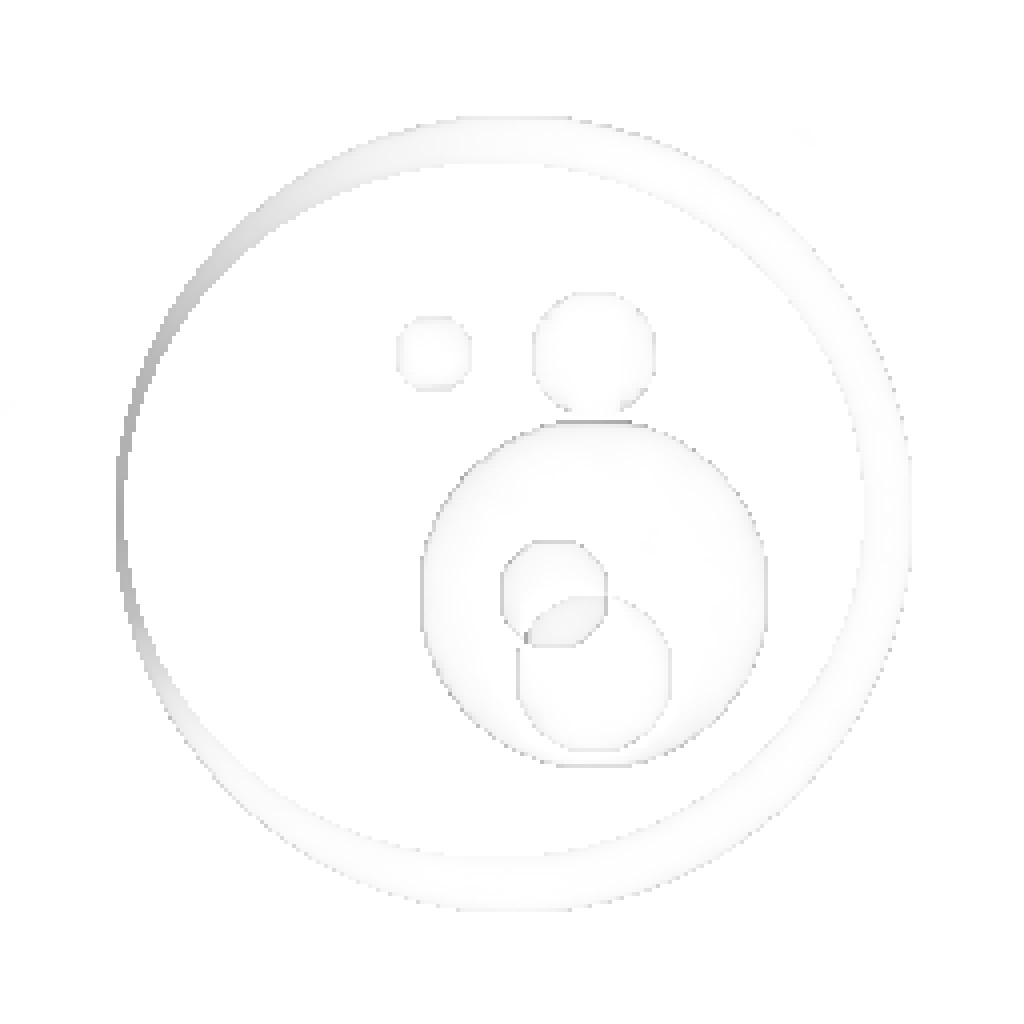}\\
	\end{minipage}
		\begin{minipage}{0.15\textwidth}
		\centering
		\includegraphics[width=\textwidth]{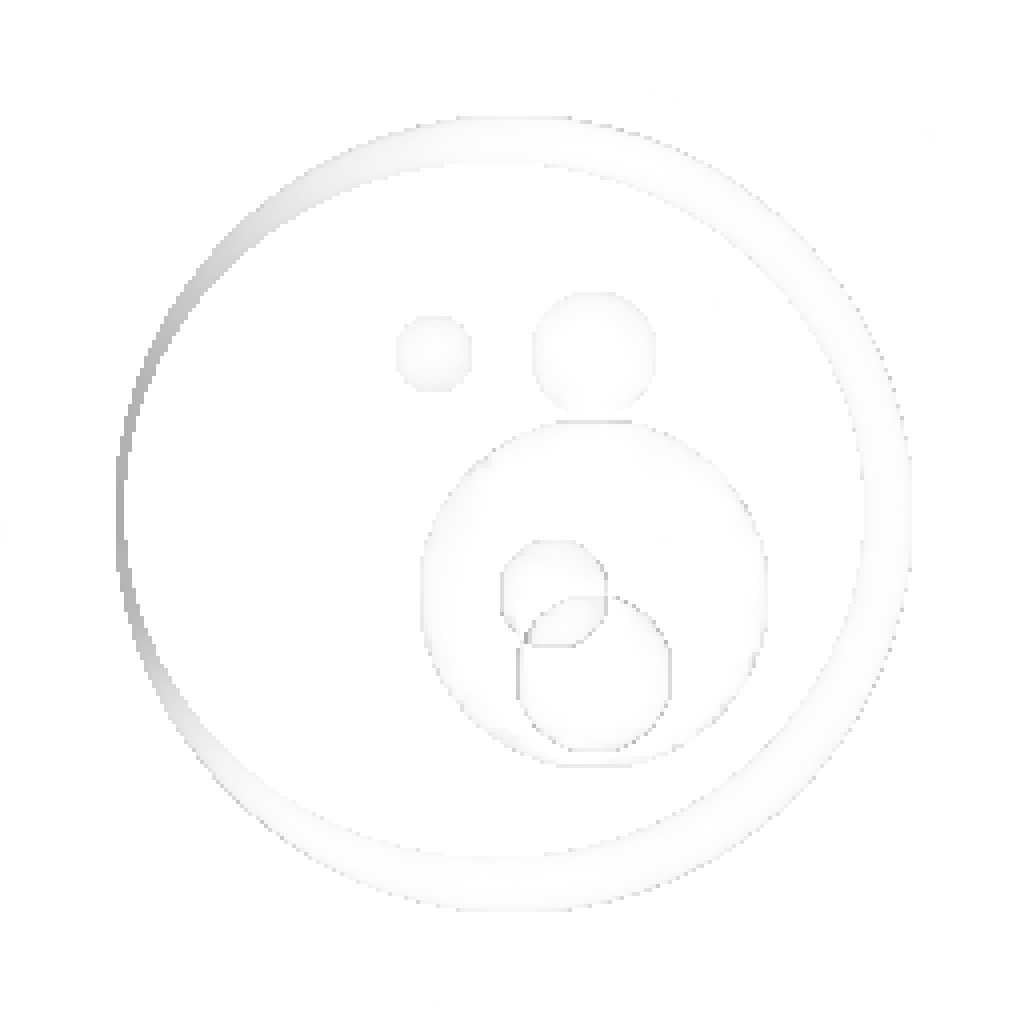}\\
	\end{minipage}
	\begin{minipage}{0.15\textwidth}
		\centering
		\includegraphics[width=\textwidth]{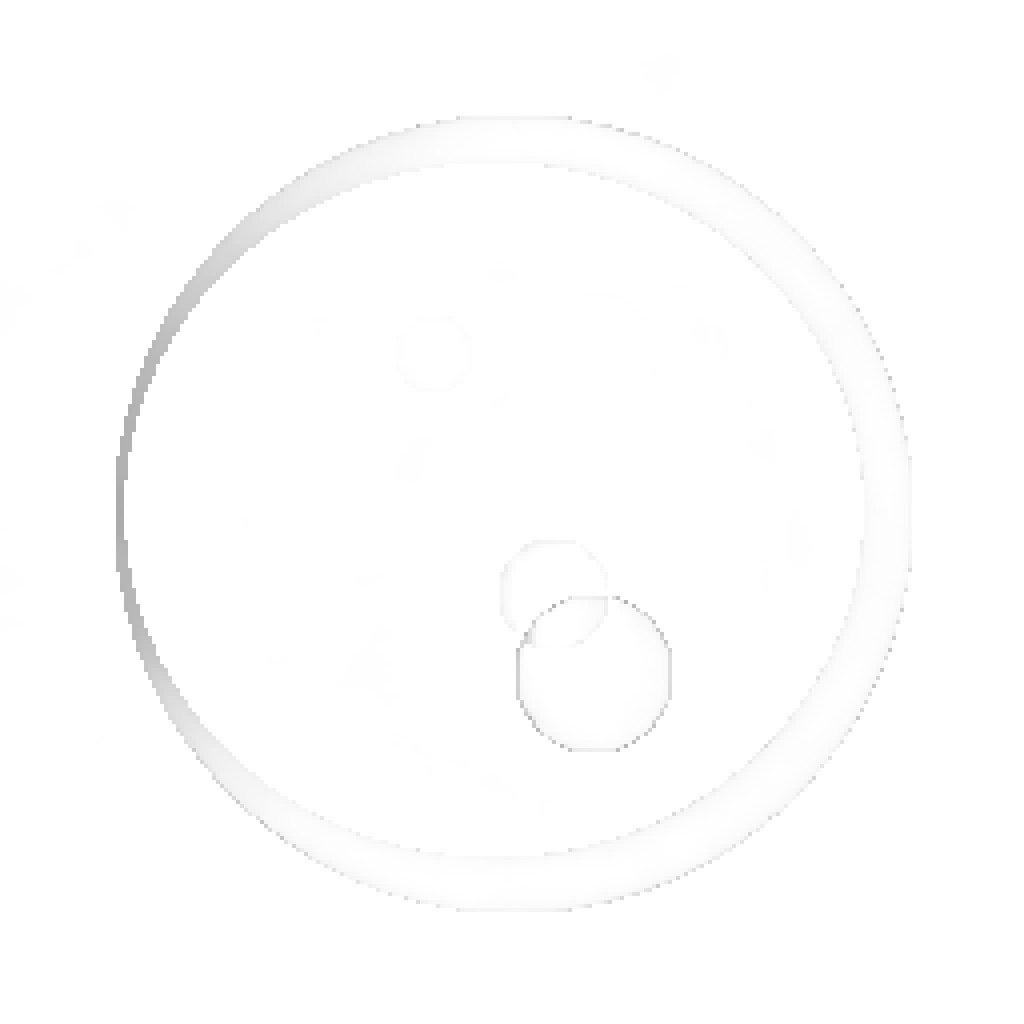}\\
	\end{minipage}

	\caption{Dynamic PAT test problem:  $t=1,10,20,30, 40, 50$ from left to right, respectively.}
	\label{fig: PATRec}
\end{figure}

\begin{figure}[ht!]
\centering
\includegraphics[width=0.99\textwidth]{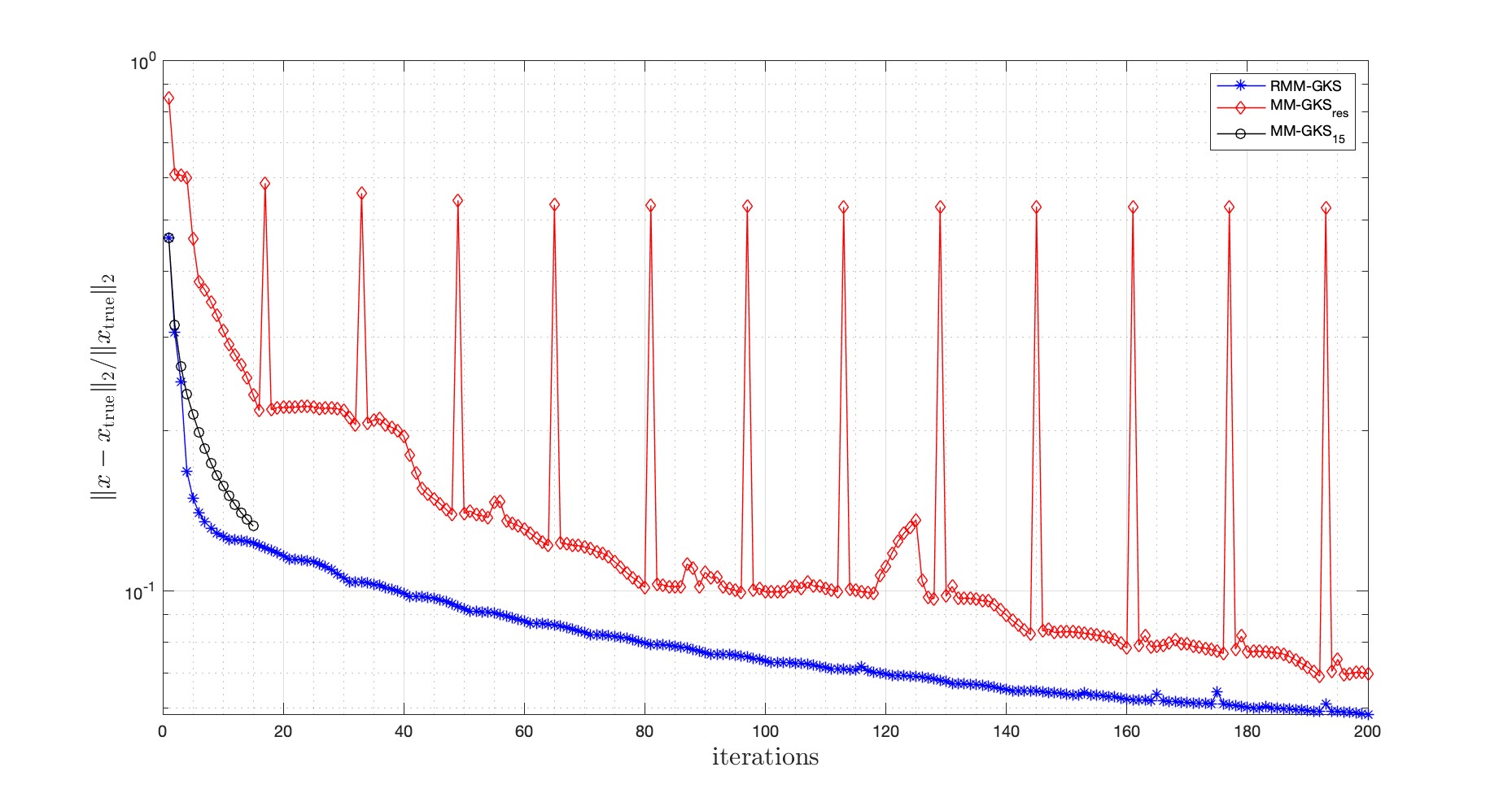}
\caption{PAT test problem: RRE for MM-GKS (200 iterations), MM-GKS (15 iterations), MM-GKS$_{\rm res}$ (200 iterations), RMM-GKS (200 iterations).}
\label{Fig: Error_telescope}
\end{figure}

\paragraph{Comparison of RMM-GKS with MM-GKS and MM-GKS$_{\rm res}$}
In this example we consider a phantom with $50$ time frames. Each image is of size $256\times 256$ and it represents the superposition of six circles. A sample of the true phantoms at time steps $i = 1, 10, 20, 30, 40, 50$ is given in the first row of Figure \ref{Fig: PATtrueandsino}. Measurements were taken from 49 distinct equidistant angles between 1 to 290 at 5 degree intervals and each projections consists of 362 radii. The measurements, also known as sinograms, $\bd_i \in \R^{362 \cdot 49}$ contain identically distributed and independent Gaussian noise with noise level $1\%$ (i.e. $\frac{\|\epsilon\|_2}{\|\bA_i(\bx_{\rm true})_i\|_2} = 0.01$). A sample of sinograms at time steps $t = 1,10, 20, 30, 40, 50$ is shown in the second row of Figure \ref{Fig: PATtrueandsino}. The total measured sinogram has 886900 pixels, i.e., $\bd\in \R^{886900}$. Note that the problem is severely underdetermined. Given the observed sinogram, the goal is to reconstruct the approximate solution that has a total of 3276800 pixels throughout 50 time instances. Let the problem of interest be formulated as follows \begin{equation}\label{eq: blockF}
\underbrace{
\begin{bmatrix} \bA_{(1)} 
\\ &  \ddots\\
& & \bA_{(n_t)} 
\end{bmatrix}}_{\bA}
\underbrace{\begin{bmatrix}
    \bx_1\\
    \vdots\\
    \bx_{n_t}
\end{bmatrix}}_{\bx}
+
\underbrace{\begin{bmatrix}
    \be_1\\
    \vdots\\
    \be_{n_t}
\end{bmatrix}}_{\be}
=
\underbrace{
\begin{bmatrix}
    \bd_1\\
    \vdots\\
    \bd_{n_t}
\end{bmatrix}}
_{\bd}.
\end{equation}
We consider solving the large-scale minimization problem 
\begin{equation}
    \min_{\bx \in \R^{n_x\cdot n_y \cdot n_t}} \|\bA\bx - \bd\|_2^2 + \lambda\|\Psi\bx\|_1,
\text{with} \quad \Psi  = \begin{bmatrix} \bI_{n_t}\otimes\bI_{n_y}\otimes\bL_{x} \\  
\bI_{n_t}\otimes \bL_{y}\otimes\bI_{n_x} \\
\bL_{t}\otimes\bI_{n_y}\otimes\bI_{n_x} 
\end{bmatrix}.
\end{equation}
The above is known as a dynamic inverse problem. Other ways to define $\Psi$ for dynamic edge-preserving inverse problems can be found in \cite{pasha2021efficient} and on the Beysian framework \cite{lan2023spatiotemporal}.
The matrix $\bL_d$ represents the discretizations of the first derivative operators in the $d$-direction, with $d = x$ (vertical direction), $d=y$ (horizontal direction), and $d=t$ (time direction) resulting in $\Psi \in \R^{9764864  \times  3276800
}$. In this example we compare our proposed approach RMM-GKS where the compression used is based on TSVD with the original MM-GKS \cite{huang2017majorization} method MM-GKS$_{\rm res}$. We let the stopping criteria for compression to be the maximum number of basis vectors, which we set to 15. For all the methods we set the maximum number of iterations to $200$. Note that for a fair comparison we consider one iteration every time that a new basis vector is added in the solution subspace which is a different iteration count for RMM-GKS from what is shown in Algorithm \eqref{Alg: RMMGKS}. Requiring that MM-GKS performs $200$ iterations means that we have to store and compute the QR factorizations for $\bA\bV_{200}$ and $\Psi\bV_{200}$ and the reorthogonalize step involves matrix vector multiplication with $200$ vectors of size  $886900$ and $9764864$ ($\bA\bV_{200}$) and $\Psi\bV_{200}$, respectively. To emphasize the need for recycling, we assume that the memory limit is $15$ solution vectors and we show the reconstruction of MM-GKS at $15$ iterations. We show the reconstructed images on time steps $t = 1,10, 20, 30, 40, 50$ in Figure \eqref{fig: PATRec} for MM-GKS at 15 iterations (first row),  MM-GKS$_{\rm res}$(second row), and RMM-GKS (third row). The RRE history for 200 iterations is shown in Figure~\ref{Fig: Error_telescope}. Results of varying the noise level from $0.1-5\%$ and computing the RRE, SSIM, and HP for all three methods are presented in Table~\ref{Table: RRE_telescope}.

\begin{table}[]
\centering
\fontsize{9pt}{10pt}\selectfont
\begin{tabular}{|l|lll|lll|lll|}
\hline
\hline
\multirow{2}{*}{\begin{tabular}[c]{@{}l@{}}Noise\\ level\end{tabular}} & \multicolumn{3}{c|}{MM-GKS$_{15}$ }                               & \multicolumn{3}{c|}{MM-GKS$_{\rm res}$ }                              & \multicolumn{3}{c|}{RMM-GKS}                                      \\ \cline{2-10} 
                                                                       & \multicolumn{1}{l|}{RRE}   & \multicolumn{1}{l|}{SSIM}  & HP & \multicolumn{1}{l|}{RRE}   & \multicolumn{1}{l|}{SSIM}  & HP & \multicolumn{1}{l|}{RRE}   & \multicolumn{1}{l|}{SSIM}  & HP \\ \hline
0.001                                                                  & \multicolumn{1}{l|}{0.126} & \multicolumn{1}{l|}{0.928} & 0.939   & \multicolumn{1}{l|}{0.119} & \multicolumn{1}{l|}{0.932} & 0.956   & \multicolumn{1}{l|}{0.057} & \multicolumn{1}{l|}{0.982} & 0.989   \\
0.005                                                                  & \multicolumn{1}{l|}{0.128} & \multicolumn{1}{l|}{0.925} & 0.937   & \multicolumn{1}{l|}{0.102} & \multicolumn{1}{l|}{0.941} & 0.967   & \multicolumn{1}{l|}{0.056} & \multicolumn{1}{l|}{0.981} & 0.989   \\
0.01                                                                   & \multicolumn{1}{l|}{0.132} & \multicolumn{1}{l|}{0.915} & 0.931   & \multicolumn{1}{l|}{0.069} & 0.975  & 0.983 & \multicolumn{1}{l|}{0.058} & \multicolumn{1}{l|}{0.981} & 0.989 \\
0.05                                                                   & \multicolumn{1}{l|}{0.166} & \multicolumn{1}{l|}{0.843} & 0.884   & \multicolumn{1}{l|}{0.084} & \multicolumn{1}{l|}{0.963} & 0.978   & \multicolumn{1}{l|}{0.079} & \multicolumn{1}{l|}{0.963} & 0.981   \\ \hline 
\hline
\end{tabular}
\caption{Algorithm setting: Maximum memory limit is 15, Image: dynamic PAT. Size: $50$ images of size $256 \times 256$, MM-GKS$_{\rm res}$ and RMM-GKS are run for $200$ iterations.}
   \label{table3}
\end{table}

\section{Conclusions and outlook}\label{sec: conclusion}
In this paper we proposed a recycling technique for MM-GKS that efficiently combines enlarging of the solution subspace with compression to overcome huge memory requirements. 
The approach we propose iteratively refines the solution subspace by keeping the maximum memory requirement fixed. A streaming variation of the method is described for scenarios when 
large-scale data exceeds memory limitations or when all the data to be processed is not available at once. Such techniques allow us to improve the computed solution through edge preserving properties, reduce memory requirements, and automatically select a regularization parameter at each iteration at a low computational cost. 
Numerical examples arising from a wide range of applications such as computerized tomography, image deblurring, and time-dependent inverse problems are used to illustrate the effectiveness of the described approach. 

While the numerical results demonstrate that our proposed recycling techniques for MM-GKS work well, in future work, we intend to analyze the effect of compression on the rate of convergence and establish strong convergence results for the recycling version of the algorithm. 

\section*{Acknowledgments}
MP gratefully acknowledges support from the NSF under award No.\ 2202846.
MP would like to further acknowledge partial support from the NSF-AWM Mentoring Travel
and the Isaac Newton Institute (INI) for Mathematical Sciences, Cambridge, for hospitality during the programme ``Rich and Nonlinear Tomography - a multidisciplinary approach'' where partial work on this manuscript was undertaken.
This material is based upon work supported by the National Science Foundation under Award No. 2208470.
MK's work is partially supported by NSF HDR grant CCF-1934553.  MK would like to acknowledge the Turner-Kirk Charitable Trust for support provided by a Kirk Distinguished Visiting Fellowship to attend the aforementioned INI programme where partial work on this manuscript was undertaken.
\appendix
\section{Dynamic computerized X-ray tomography}\label{app:mm}
\begin{figure}[h!]
\begin{centering}
    	\begin{minipage}{0.18\textwidth}
		\includegraphics[angle = -90, width=\textwidth]{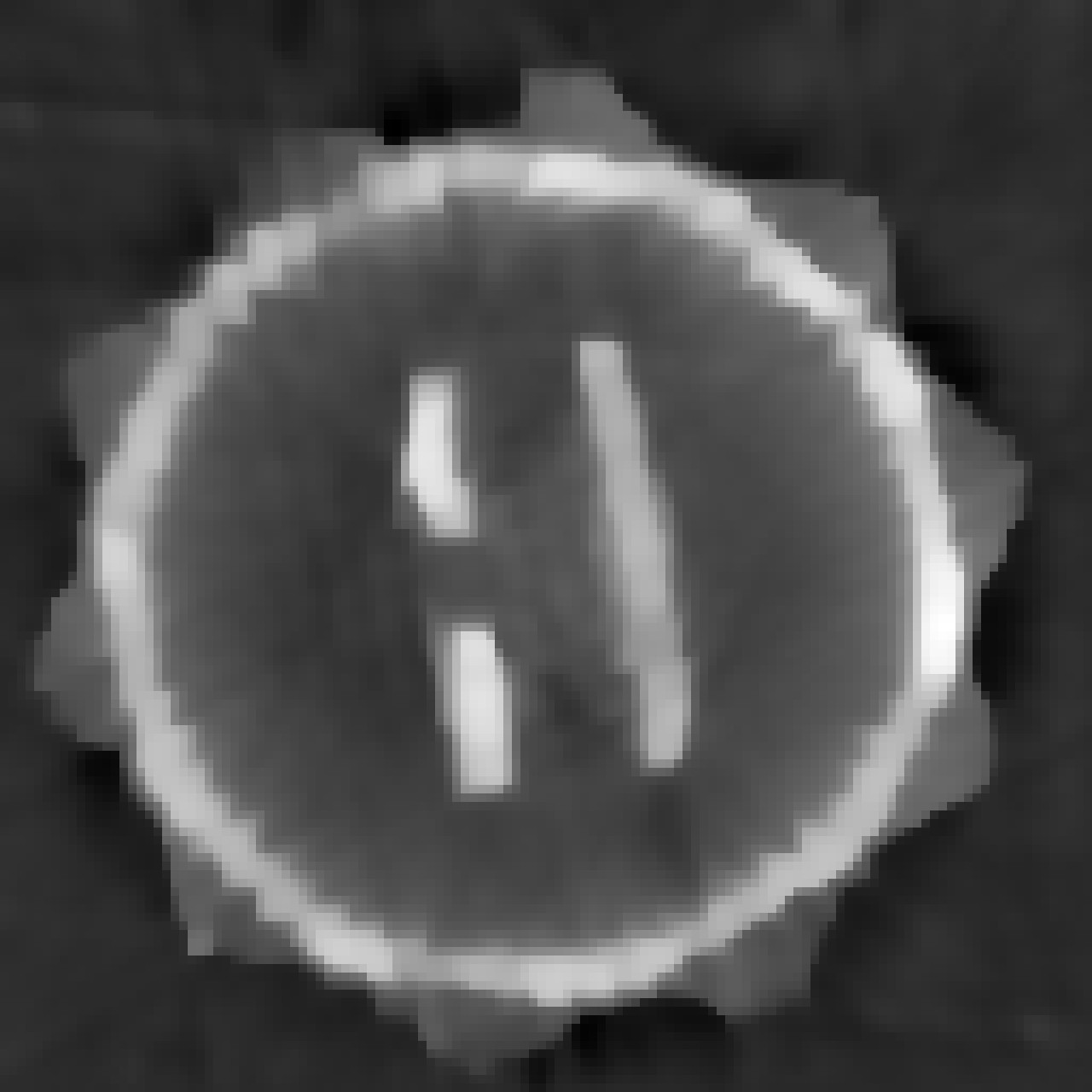}
	\end{minipage}
	\begin{minipage}{0.18\textwidth}
		\includegraphics[angle = -90, width=\textwidth]{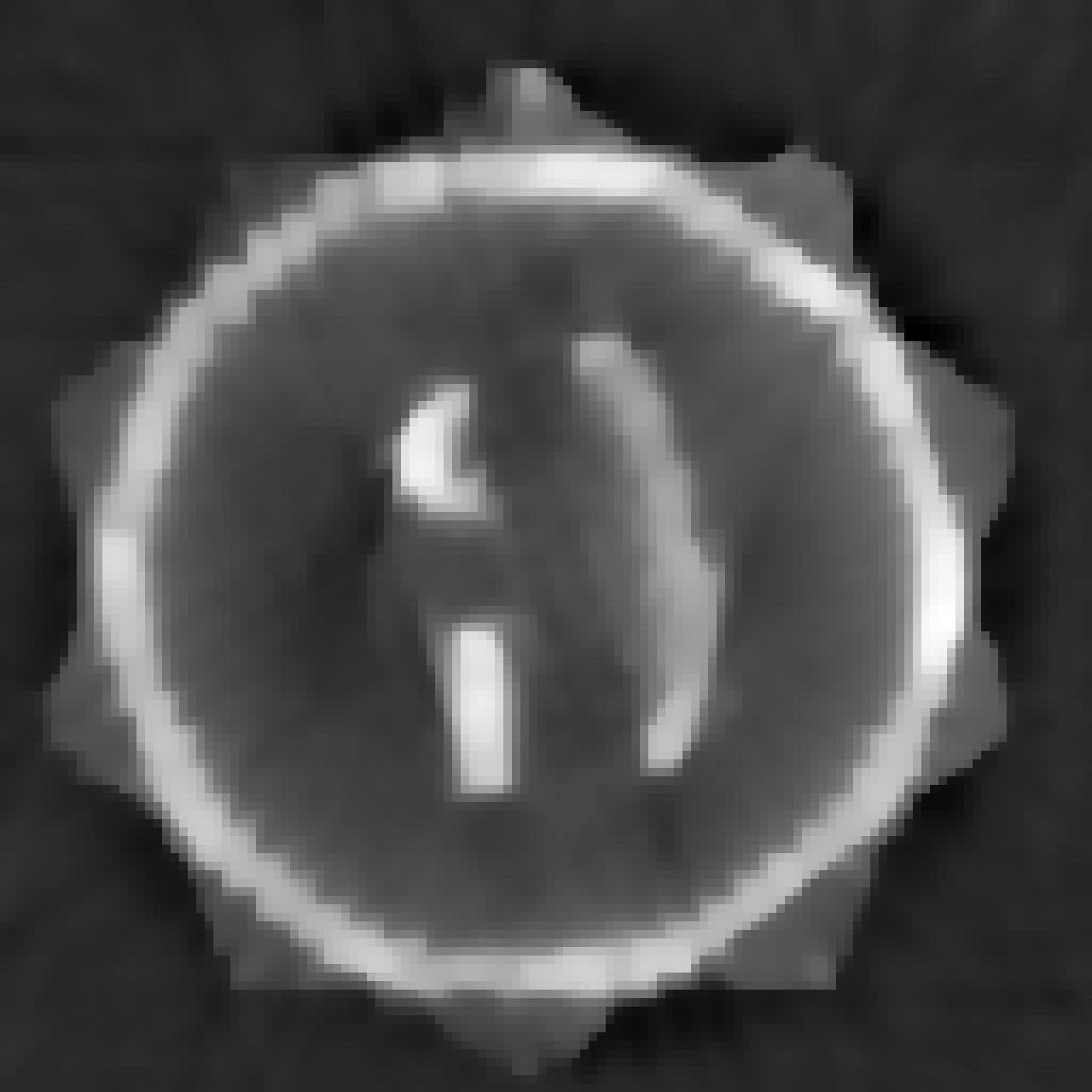}
	\end{minipage}
	\begin{minipage}{0.18\textwidth}
		\includegraphics[angle = -90, width=\textwidth]{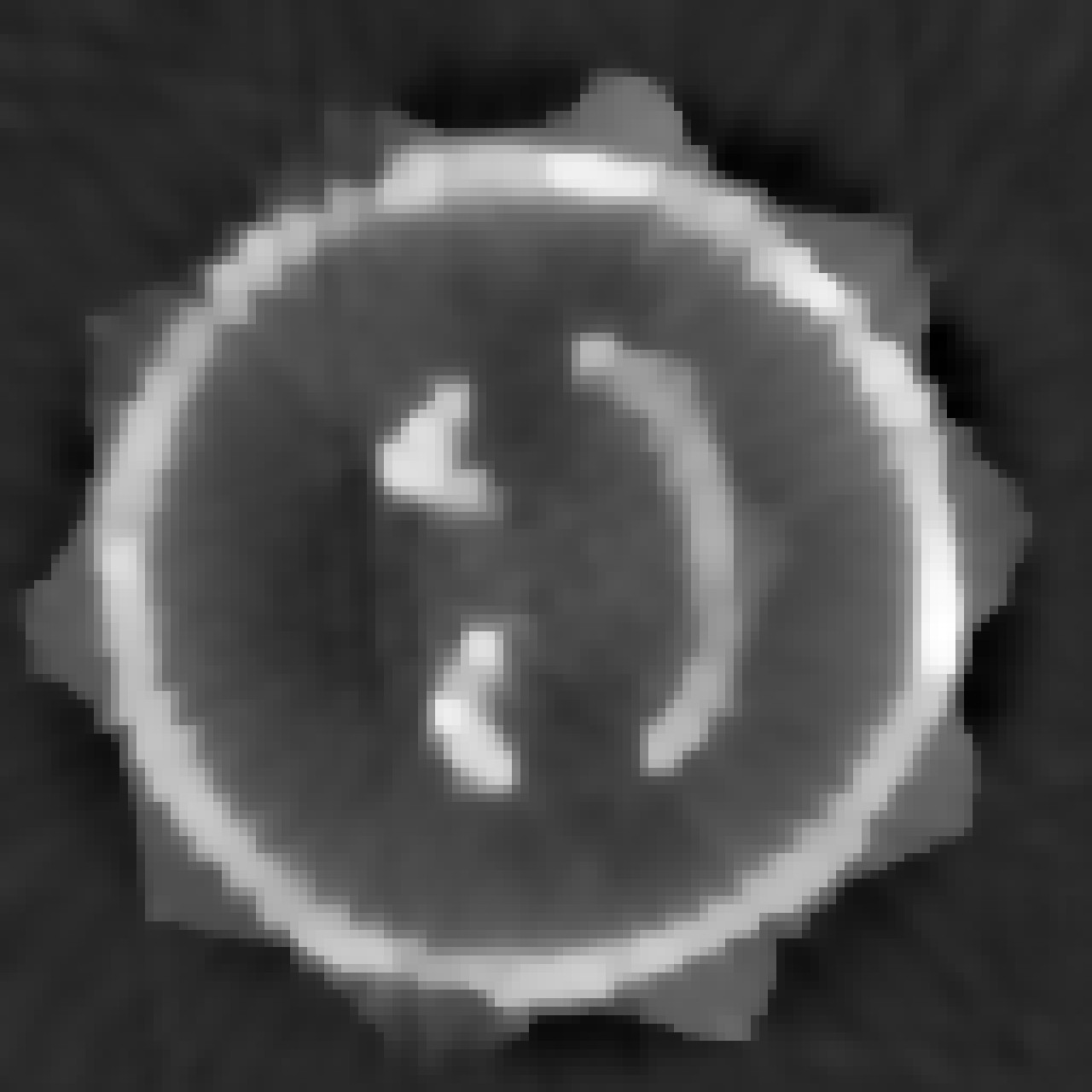}
	\end{minipage}
	\begin{minipage}{0.18\textwidth}
		\includegraphics[angle = -90, width=\textwidth]{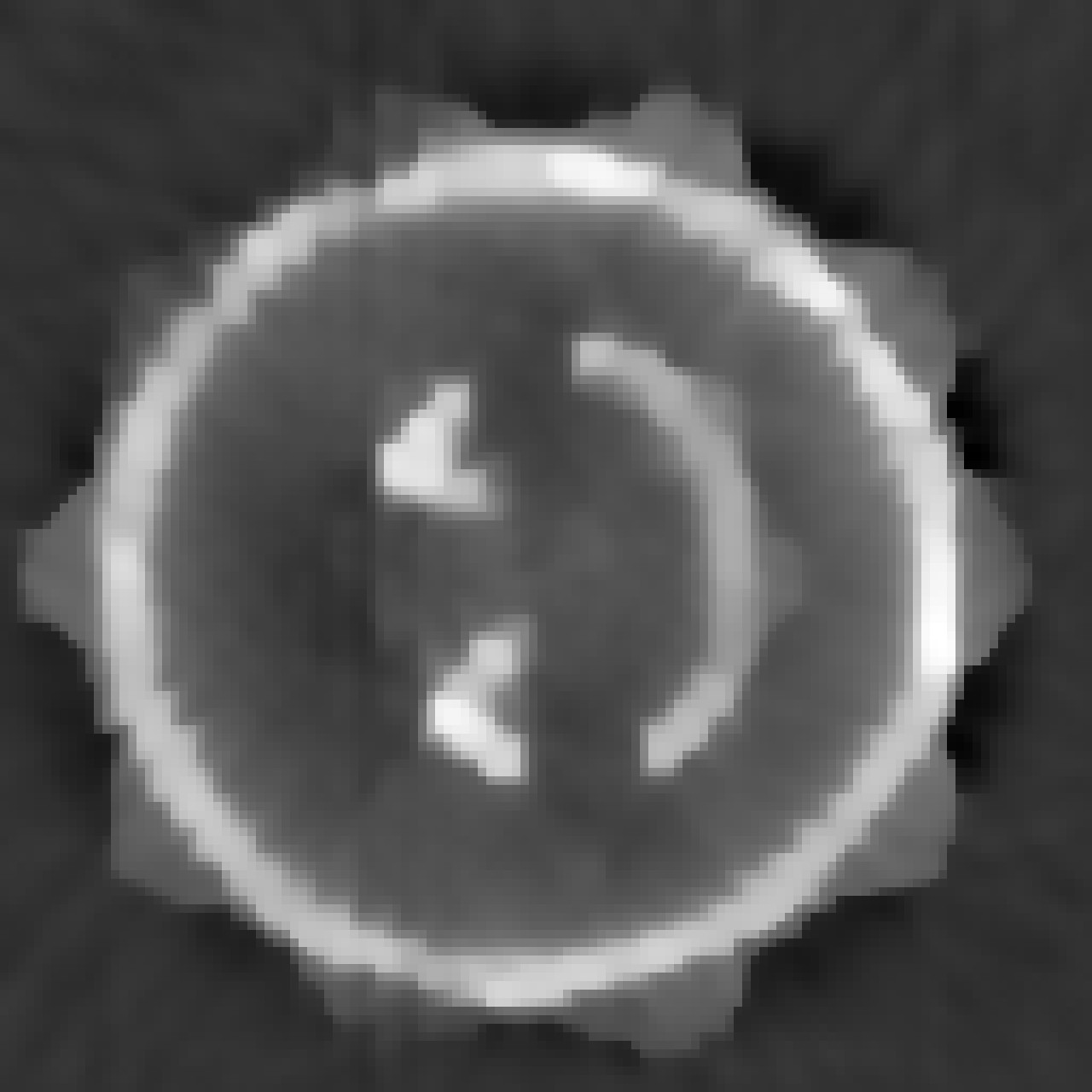}
	\end{minipage}
 	\begin{minipage}{0.18\textwidth}
		\includegraphics[angle = -90, width=\textwidth]{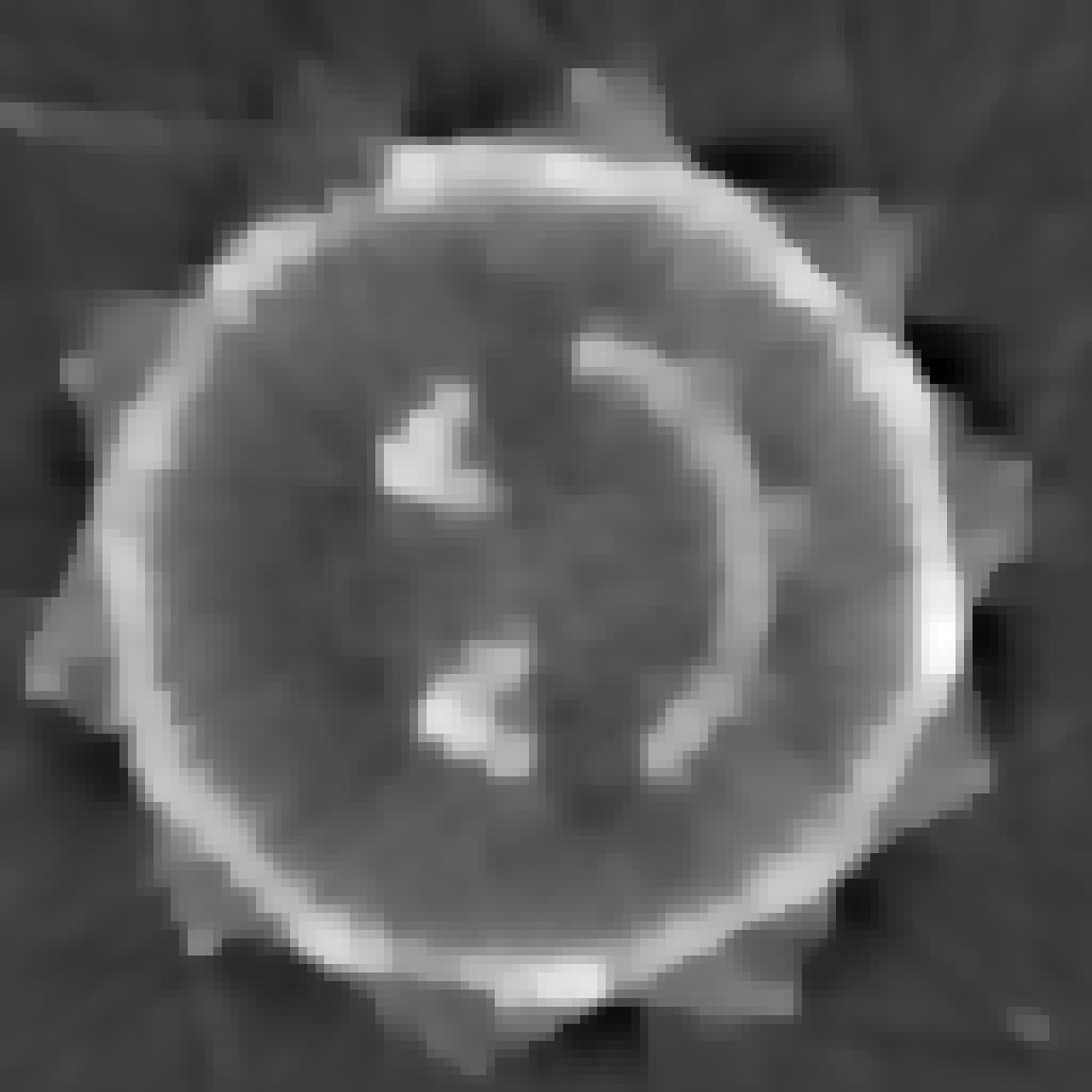}
	\end{minipage}\\
	\begin{minipage}{0.18\textwidth}
		\includegraphics[width=\textwidth,angle =-90]{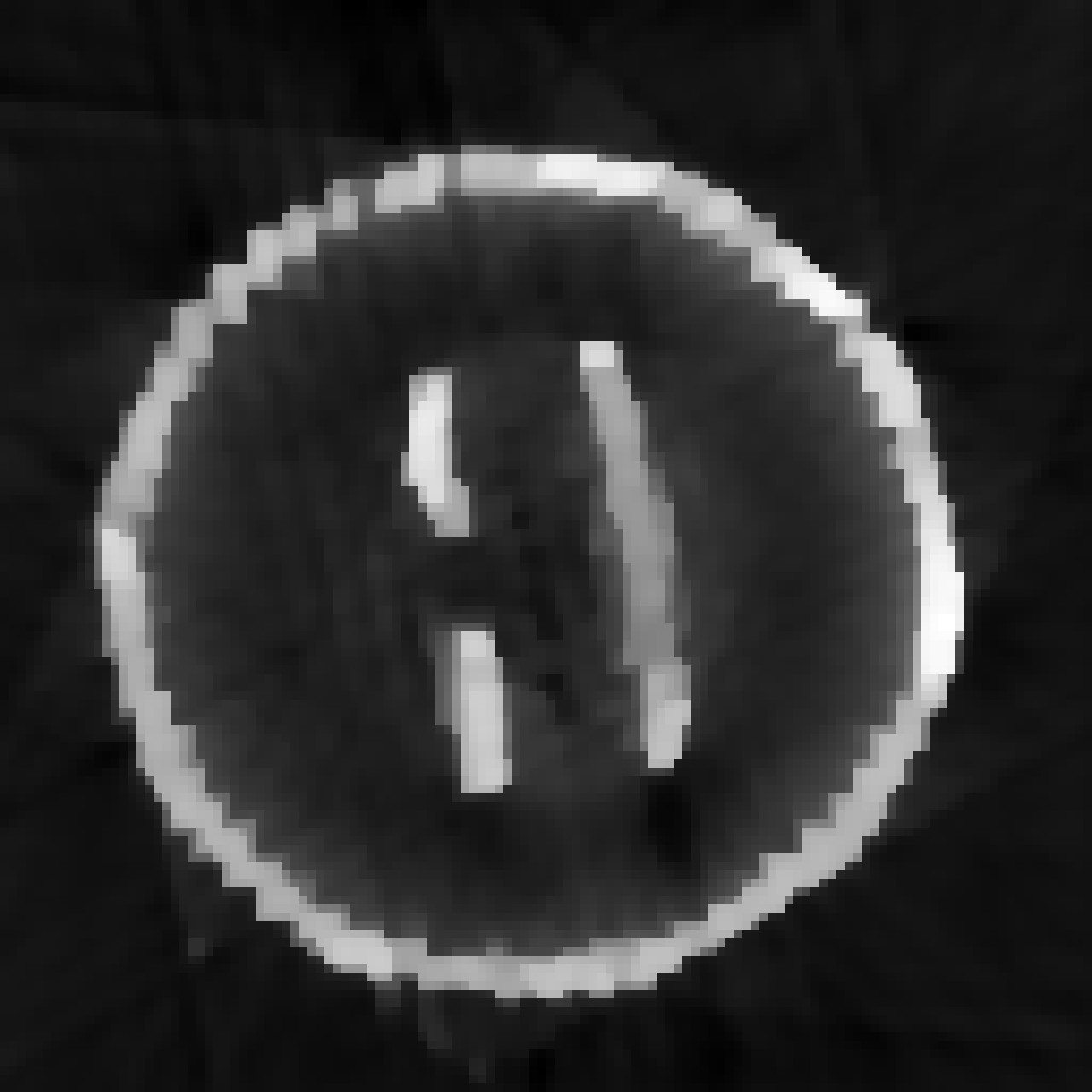}
	\end{minipage}
	\begin{minipage}{0.18\textwidth}
		\includegraphics[width=\textwidth,angle =-90]{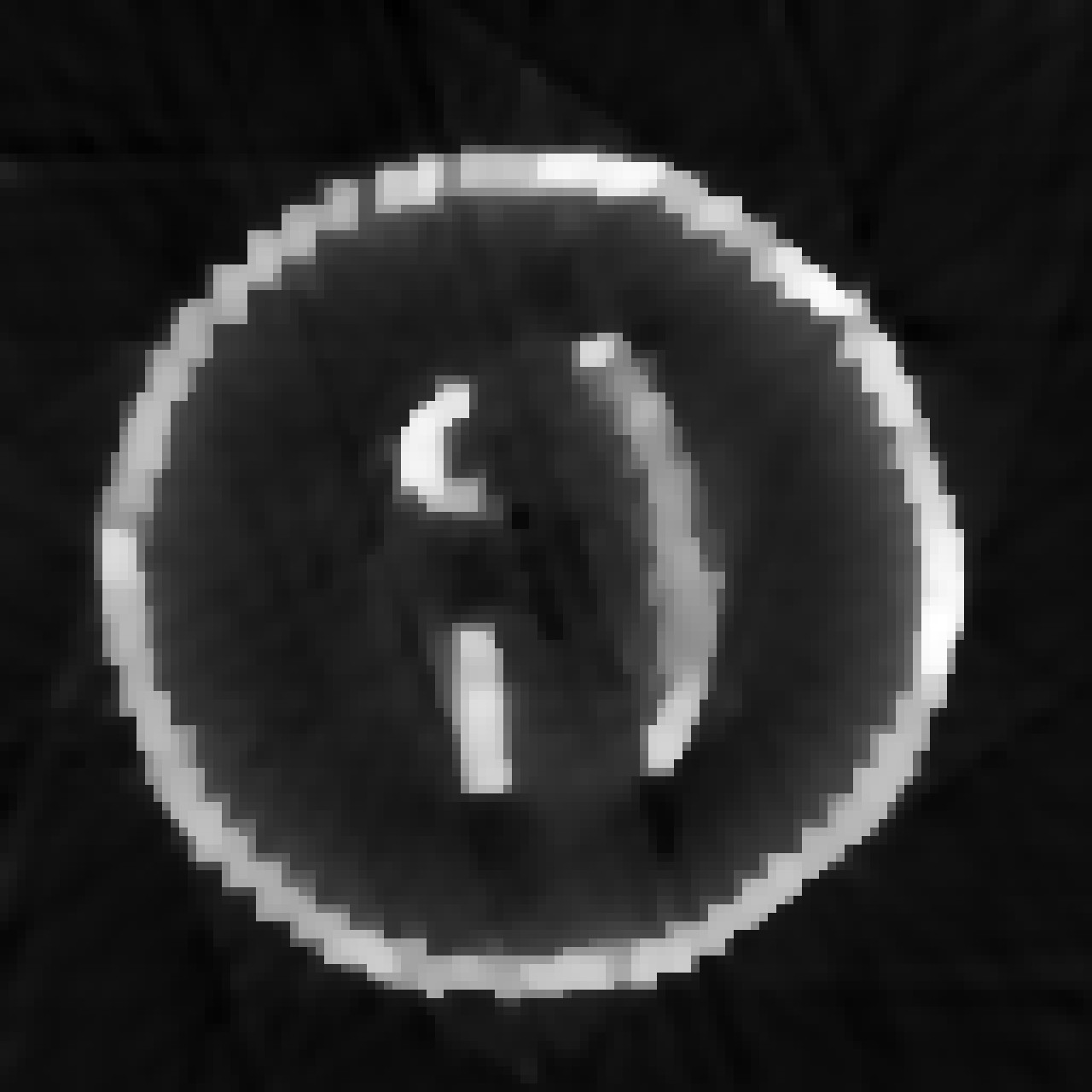}
	\end{minipage}
	\begin{minipage}{0.18\textwidth}
		\includegraphics[width=\textwidth,angle =-90]{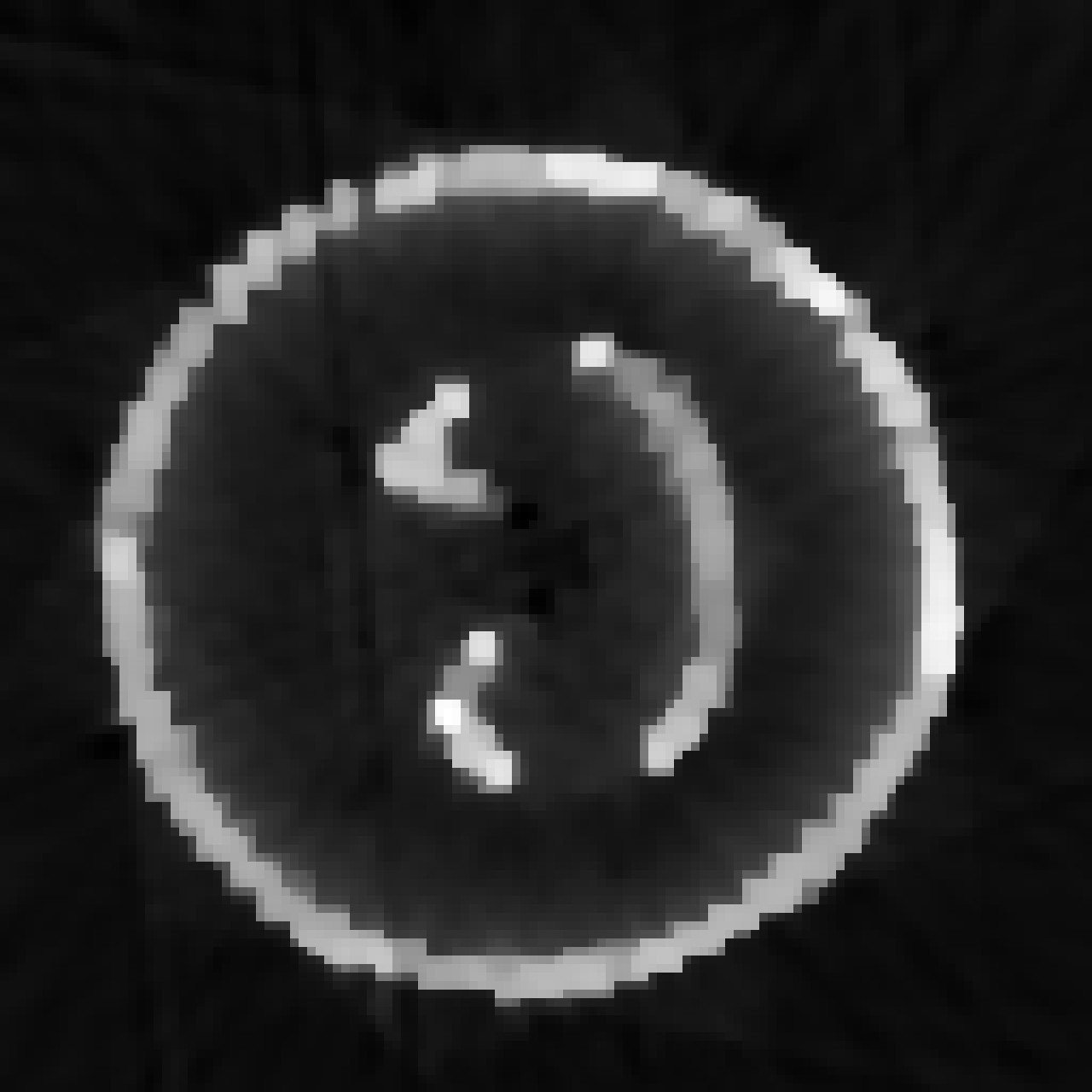}
	\end{minipage}
	\begin{minipage}{0.18\textwidth}
		\includegraphics[width=\textwidth,angle =-90]{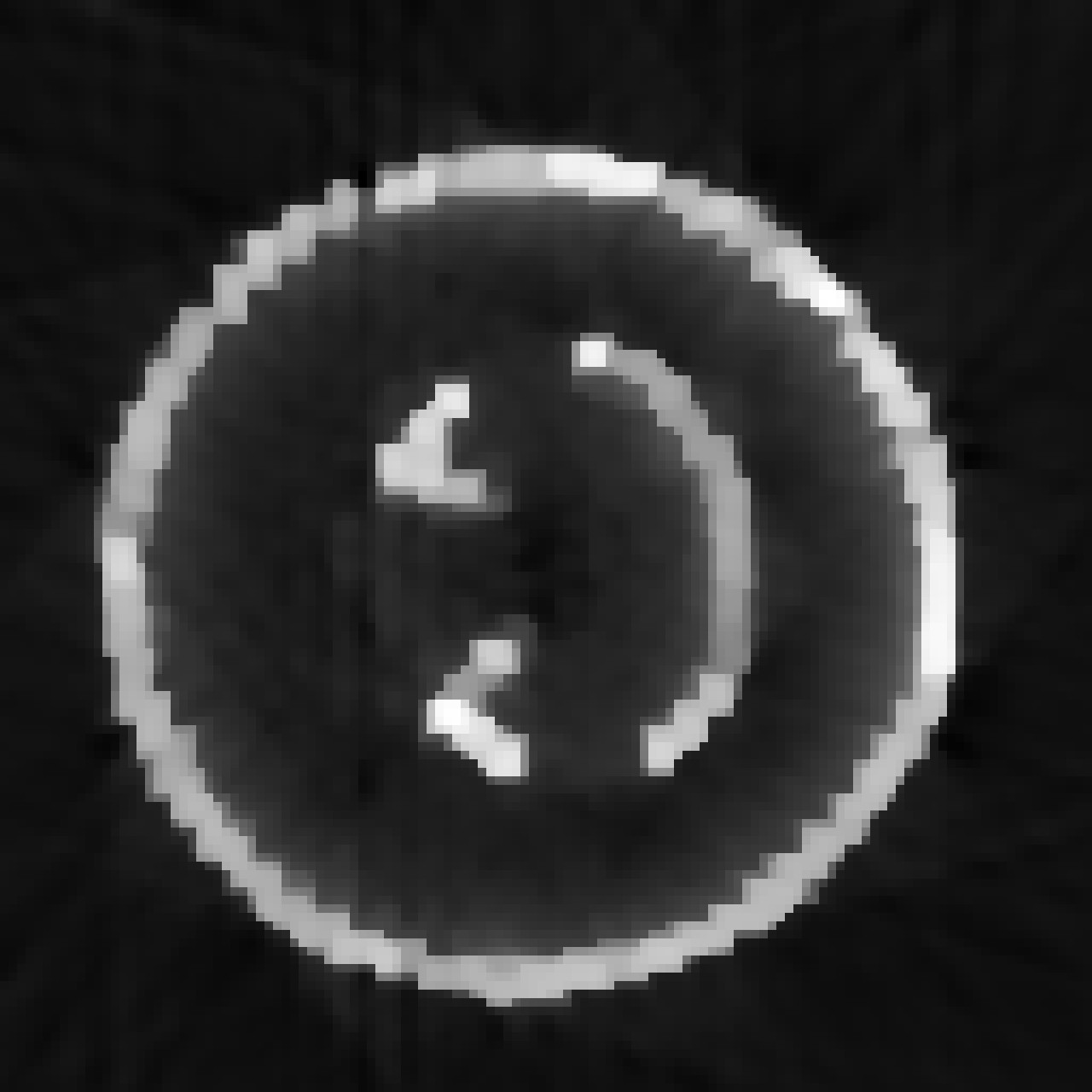}
	\end{minipage}
 	\begin{minipage}{0.18\textwidth}
		\includegraphics[width=\textwidth, angle =-90]{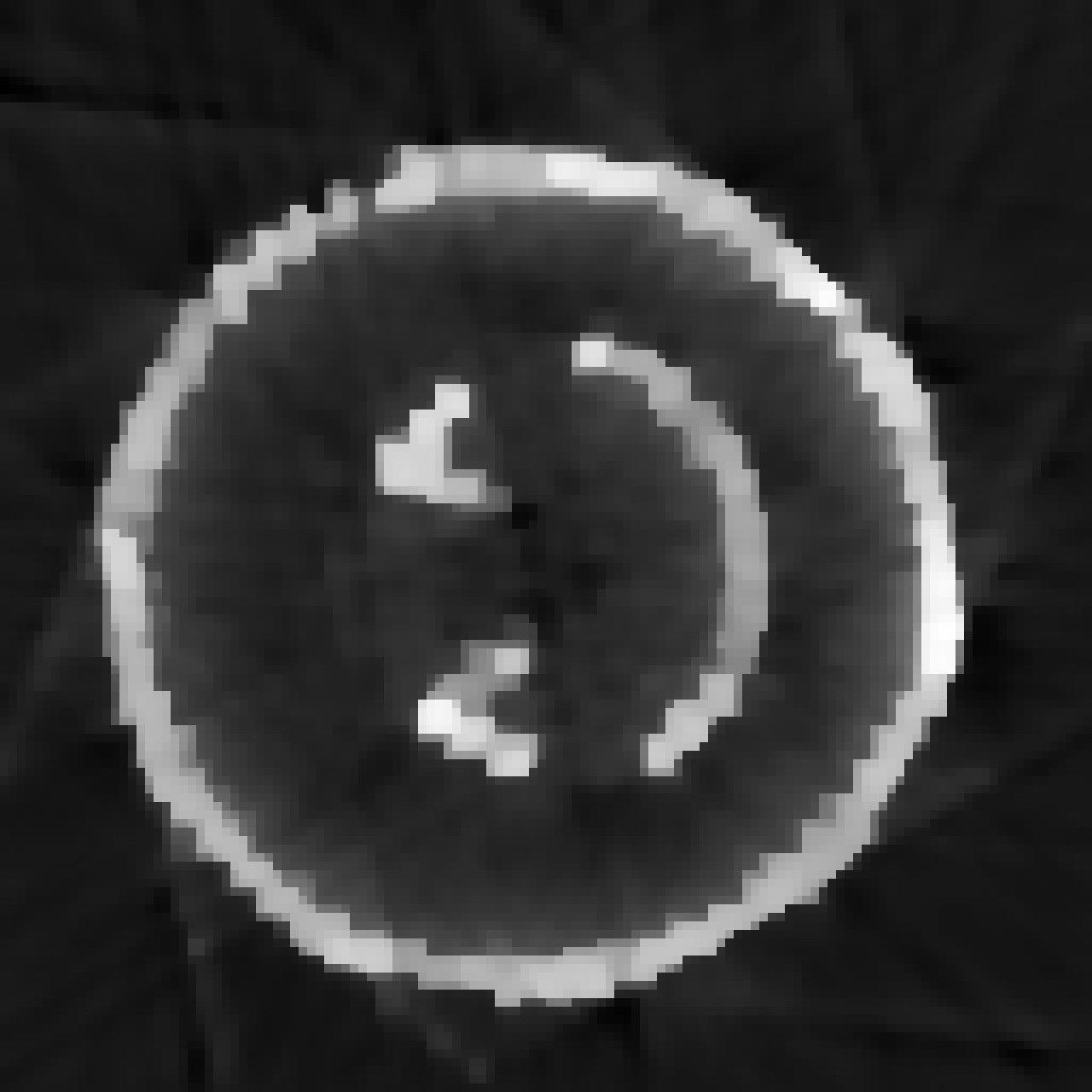}
	\end{minipage}\\
\end{centering}
	\caption{Reconstruction results for the emoji test problem with $n_a =10$. The first row shows reconstructions with MM-GKS at 22 iterations. The second row represents reconstructions with RMM-GKS until convergence when the memory limit is kept $k_{min} = 22$ at time instances $t = 2, 10, 18, 31$ (from left to right).}
	\label{Fig: smile_10angles}
\end{figure}

In this example we consider a dynamic computerized tomogrpahy example with real data. The true solution is not available, hence we only use visual inspection to compare reconstruction qualities.
We consider the dataset \verb|DataDynamic_128x30.mat| that contains the emoji phantoms measured at the University of Helsinki \cite{meaney2018tomographic}. The available data represents $n_t=33$ time steps of a series of the X-ray sinogram of emojis created of small ceramic stones. The observations are obtained by shining $217$ projections from $n_a =30$ angles from which we extract the information collected from $n_a = 10$ projection angles for our analysis. The underdetermined problems $\bA^{(t)}\bu^{(t)} + \be^{(t)} = \bd^{(t)}$, $t = 1,2,\dots, 33$ where $\bA^{(t)} \in \R^{2,170\times 16,384}$ yield a large $\bA \in \R^{71,610 \times 540,672}$ and the measurement vector $\bd \in \R^{71,610}$ containing the measured sinograms $\bd^{(t)} \in \R^{2,170}$ obtained from 217 projections around 10 equidistant angles. This example represents a limited angle computerized dynamic inverse problem. In this current work we omit explaining the need to solving the large-scale problem as it has been shown in \cite{pasha2021efficient}. In this work we illustrate the benefits of using our RMM-GKS in real data when the memory requirements are limited. In particular, we consider the above setup with 33 images of size $128\times 128$. We assume that we can only store $22$ solution subspace vectors. We run MM-GKS for $22$ iterations since 22 iterations will produce $22$ basis vectors. If we have limited memory, we can not perform anymore iterations of MM-GKS. On the other hand, since RMM-GKS alternates between enlarging and compressing, we can run as many iterations as we need to receive a high quality reconstruction while the largest number of subspace vectors remains fixed. We set the number of $k_{\rm max} = 22$ and $k_{\rm min} = 3$, meaning that at after each compression of the subspace, we enlarge by adding $19$ new solution subspace vectors in the subspace. By visual inspection of the reconstructions shown in Figure \ref{Fig: smile_10angles} we observe that if the memory limit is reached without MM-GKS having converged (see for instance reconstructions at time steps $t=2, 10, 18, 31$ in the first row of Figure \ref{Fig: smile_10angles}), we risk to reconstruct a low quality solution. Nevertheless, RMM-GKS (reconstructions shown in the second row of Figure \ref{Fig: smile_10angles}) continues to improve the reconstruction quality as we can run the method for more iterations while the memory requirement is low. 
\bibliographystyle{plain}
\bibliography{arxiv_September_27_2023}

\begin{thebibliography}{10}

\bibitem{ahuja2015recycling}
Kapil Ahuja, Peter Benner, Eric de~Sturler, and Lihong Feng.
\newblock Recycling {BiCGSTAB} with an application to parametric model order reduction.
\newblock {\em SIAM Journal on Scientific Computing}, 37(5):S429--S446, 2015.

\bibitem{buccini2020modulus}
A~Buccini, M~Pasha, and L~Reichel.
\newblock Modulus-based iterative methods for constrained $\ell_p-\ell_q$ minimization.
\newblock {\em Inverse Problems}, 36(8):084001, 2020.

\bibitem{buccini2021linearized}
Alessandro Buccini, Mirjeta Pasha, and Lothar Reichel.
\newblock Linearized {K}rylov subspace {B}regman iteration with nonnegativity constraint.
\newblock {\em Numerical Algorithms}, 87:1177--1200, 2021.

\bibitem{buccini2023limited}
Alessandro Buccini and Lothar Reichel.
\newblock Limited memory restarted $\ell_p-\ell_q$ minimization methods using generalized {K}rylov subspaces.
\newblock {\em Advances in Computational Mathematics}, 49(2):26, 2023.

\bibitem{chen2015reduced}
Yanlai Chen.
\newblock Reduced basis decomposition: a certified and fast lossy data compression algorithm.
\newblock {\em Computers \& Mathematics with Applications}, 70(10):2566--2574, 2015.

\bibitem{chung2022}
Julianne Chung, Matthias Chung, Silvia Gazzola, and Mirjeta Pasha.
\newblock Efficient learning methods for large-scale optimal inversion design.
\newblock {\em Numerical Algebra, Control and Optimization}, 2022.

\bibitem{chung2019flexible}
Julianne Chung and Silvia Gazzola.
\newblock Flexible {K}rylov methods for $\ell_p$ regularization.
\newblock {\em SIAM Journal on Scientific Computing}, 41(5):S149--S171, 2019.

\bibitem{chung2017motion}
Julianne Chung and Linh Nguyen.
\newblock Motion estimation and correction in photoacoustic tomographic reconstruction.
\newblock {\em SIAM J. Imaging Sci.}, 10(1):216--242, 2017.

\bibitem{chung2018efficient}
Julianne Chung, Arvind~K Saibaba, Matthew Brown, and Erik Westman.
\newblock Efficient generalized {G}olub--{K}ahan based methods for dynamic inverse problems.
\newblock {\em Inverse Problems}, 34(2):024005, 2018.

\bibitem{gazzola2018ir}
Silvia Gazzola, Per~Christian Hansen, and James~G Nagy.
\newblock {I}{R} {T}ools: {A} {MATLAB} package of iterative regularization methods and large-scale test problems.
\newblock {\em Numerical Algorithms}, 81:773--811, 2019.

\bibitem{gazzola2020inner}
Silvia Gazzola, Misha~E Kilmer, James~G Nagy, Oguz Semerci, and Eric~L Miller.
\newblock An inner--outer iterative method for edge preservation in image restoration and reconstruction.
\newblock {\em Inverse Problems}, 36(12):124004, 2020.

\bibitem{golub2013matrix}
GH~Golub and CF~Van~Loan.
\newblock Matrix {C}omputations.
\newblock {\em 4th edition {T}he {J}ohns {H}opkins {U}niversity {P}ress, Baltimore, MD}, 2013.

\bibitem{huang2017majorization}
G~Huang, A~Lanza, S~Morigi, L~Reichel, and F~Sgallari.
\newblock Majorization--minimization generalized {K}rylov subspace methods for $\ell_p-\ell_q$ optimization applied to image restoration.
\newblock {\em BIT Numerical Mathematics}, 57(2):351--378, 2017.

\bibitem{hunter2004tutorial}
David~R Hunter and Kenneth Lange.
\newblock A tutorial on {MM} algorithms.
\newblock {\em The American Statistician}, 58(1):30--37, 2004.

\bibitem{jiang2021hybrid}
Jiahua Jiang, Julianne Chung, and Eric de~Sturler.
\newblock Hybrid projection methods with recycling for inverse problems.
\newblock {\em SIAM Journal on Scientific Computing}, 43(5):S146--S172, 2021.

\bibitem{keuchel2016combination}
S{\"o}ren Keuchel, Jan Biermann, and Otto von Estorff.
\newblock A combination of the fast multipole boundary element method and {K}rylov subspace recycling solvers.
\newblock {\em Engineering Analysis with Boundary Elements}, 65:136--146, 2016.

\bibitem{kilmer2006recycling}
Misha~E Kilmer and Eric de~Sturler.
\newblock Recycling subspace information for diffuse optical tomography.
\newblock {\em SIAM Journal on Scientific Computing}, 27(6):2140--2166, 2006.

\bibitem{lampe2012large}
J{\"o}rg Lampe, Lothar Reichel, and Heinrich Voss.
\newblock Large-scale {T}ikhonov regularization via reduction by orthogonal projection.
\newblock {\em Linear algebra and its applications}, 436(8):2845--2865, 2012.

\bibitem{lan2023spatiotemporal}
Shiwei Lan, Mirjeta Pasha, and Shuyi Li.
\newblock Spatiotemporal {B}esov priors for {B}ayesian inverse problems.
\newblock {\em arXiv preprint arXiv:2306.16378}, 2023.

\bibitem{lange2016mm}
Kenneth Lange.
\newblock {\em {MM} optimization algorithms}.
\newblock SIAM, 2016.

\bibitem{lanza2015generalized}
Alessandro Lanza, Serena Morigi, Lothar Reichel, and Fiorella Sgallari.
\newblock A generalized {K}rylov subspace method for $\ell_p-\ell_q$ minimization.
\newblock {\em SIAM Journal on Scientific Computing}, 37(5):S30--S50, 2015.

\bibitem{lucka2018enhancing}
Felix Lucka, Nam Huynh, Marta Betcke, Edward Zhang, Paul Beard, Ben Cox, and Simon Arridge.
\newblock Enhancing compressed sensing 4{D} photoacoustic tomography by simultaneous motion estimation.
\newblock {\em SIAM Journal on Imaging Sciences}, 11(4):2224--2253, 2018.

\bibitem{meaney2018tomographic}
Alexander Meaney, Zenith Purisha, and Samuli Siltanen.
\newblock Tomographic {X}-ray data of 3{D} emoji.
\newblock {\em arXiv preprint arXiv:1802.09397}, 2018.

\bibitem{parks2006recycling}
Michael~L Parks, Eric de~Sturler, Greg Mackey, Duane~D Johnson, and Spandan Maiti.
\newblock Recycling {K}rylov subspaces for sequences of linear systems.
\newblock {\em SIAM Journal on Scientific Computing}, 28(5):1651--1674, 2006.

\bibitem{pasha2021efficient}
Mirjeta Pasha, Arvind~K Saibaba, Silvia Gazzola, Malena~I Espanol, and Eric de~Sturler.
\newblock A computational framework for edge-preserving regularization in dynamic inverse problems.
\newblock {\em Electronic Transactions on Numerical Analysis}, 58:486--516, 2023.

\bibitem{reisenhofer2018haar}
Rafael Reisenhofer, Sebastian Bosse, Gitta Kutyniok, and Thomas Wiegand.
\newblock A {H}aar wavelet-based perceptual similarity index for image quality assessment.
\newblock {\em Signal Processing: Image Communication}, 61:33--43, 2018.

\bibitem{rodriguez2008efficient}
Paul Rodr{\'\i}guez and Brendt Wohlberg.
\newblock Efficient minimization method for a generalized total variation functional.
\newblock {\em IEEE Transactions on Image Processing}, 18(2):322--332, 2008.

\bibitem{soodhalter2016block}
Kirk~M Soodhalter.
\newblock Block {K}rylov subspace recycling for shifted systems with unrelated right-hand sides.
\newblock {\em SIAM Journal on Scientific Computing}, 38(1):A302--A324, 2016.

\bibitem{soodhalter2014krylov}
Kirk~M Soodhalter, Daniel~B Szyld, and Fei Xue.
\newblock {K}rylov subspace recycling for sequences of shifted linear systems.
\newblock {\em Applied Numerical Mathematics}, 81:105--118, 2014.

\bibitem{wang2007large}
Shun Wang, Eric de~Sturler, and Glaucio~H Paulino.
\newblock Large-scale topology optimization using preconditioned {K}rylov subspace methods with recycling.
\newblock {\em International journal for numerical methods in engineering}, 69(12):2441--2468, 2007.

\bibitem{SSIM}
Zhou Wang, Alan~Conrad Bovik, Hamid~Rahim Sheikh, and Eero~P Simoncelli.
\newblock Image quality assessment: {F}rom error visibility to structural similarity.
\newblock {\em IEEE Transactions on Image Processing}, 13(4):600--612, 2004.

\end{thebibliography}
\end{document}